\documentclass[11pt, letterpaper]{elsarticle}
\usepackage{pdfpages}
\usepackage{comment}
\usepackage{amsmath}
\usepackage{amssymb}
\usepackage{graphicx}
\usepackage{multirow}
\usepackage{latexsym}
\usepackage{epic}
\usepackage{url}
\usepackage{soul}
\usepackage{afterpage}
\usepackage{subcaption}
\usepackage[unicode, pdftex]{hyperref}
\usepackage{multicol}
\usepackage{bbm}
\usepackage{setspace}
\usepackage{enumitem}
\usepackage{cleveref}
\usepackage[toc]{appendix}
\usepackage[norelsize,ruled, lined,linesnumbered]{algorithm2e}
\usepackage{tikz}
\usepackage{tikz,fullpage}
\usepackage{pgf}
\usepackage{pgfplots}
\numberwithin{equation}{section} 
\usetikzlibrary{
	pgfplots.fillbetween,
}
\usetikzlibrary{arrows,automata}
\usepackage{tkz-berge}
\usepackage{amsthm}
\usepackage[symbol]{footmisc}
\usepackage[mathscr]{euscript}
\usepackage{pifont}
\usepackage{empheq}
\usepackage{pgfplotstable}
\definecolor{forestgreen}{rgb}{0.13, 0.55, 0.13}
\definecolor{limegreen}{rgb}{0.7, 0.75, 0.0}
\definecolor{darkbluegray}{HTML}{404a69}
\pgfplotsset{compat=newest}

\DeclareMathOperator{\argmax}{arg\,max}

\DeclareMathOperator{\diag}{diag}

\makeatletter
\newtheorem*{rep@theorem}{\rep@title}
\newcommand{\newreptheorem}[2]{%
\newenvironment{rep#1}[1]{%
 \def\rep@title{#2 \ref{##1}}%
 \begin{rep@theorem}}%
 {\end{rep@theorem}}}
\makeatother

\newreptheorem{example}{Example}
\newtheorem{theorem}{Theorem}

\newtheorem{lemma}{Lemma}

\newtheorem{corollary}{Corollary}

\newtheorem{example}{Example}

\usepackage{xcolor}
\usepackage[colorinlistoftodos,prependcaption,textsize=tiny]{todonotes}
\usepackage{xargs}

\usepackage{caption}
\usepackage{subcaption}
\usepackage{multirow}

\usepackage{float}
\usepackage[section]{placeins}
\usepackage{arydshln}
\usepackage{latexsym}
\usepackage{epic}
\usepackage{amsmath,amsthm,amssymb}
\usepackage{graphicx}
\usepackage{url}
\usepackage{pgfplots}\usetikzlibrary{plotmarks}
\usepackage{soul}
\usepackage{afterpage}
\usepackage{cleveref}
\usepackage{bbm}
\usetikzlibrary{decorations.markings}
\usepackage[left=0.9in,top=0.9in,right=0.9in,bottom=0.9in,nohead]{geometry}
\usepackage{setspace}
\usepackage{appendix}
\usepackage{lineno}
\DeclareGraphicsExtensions{.pdf,.png,.jpg,.bmp}

\allowdisplaybreaks

\newcommand{\xmark}{\ding{55}}

\begin{document}
\begin{frontmatter}
\title{Data-driven interdiction with asymmetric cost uncertainty: a distributionally robust optimization approach}

\author[label2]{Sergey S.~Ketkov\footnote[2]{Corresponding author. Email: sergei.ketkov@business.uzh.ch; phone: +41-78-301-85-21}} 
\author[label2]{Oleg A. Prokopyev}

\address[label2]{Department of Business Administration, University of Zurich, Zurich, 8032, Switzerland}
\begin{abstract} \onehalfspacing 
\looseness-1 We consider a class of stochastic interdiction games between an upper-level decision-maker~(the \textit{leader}) and a lower-level decision-maker~(the \textit{follower}), where uncertainty lies in the follower's objective function coefficients.
Specifically, the follower's profits (or costs) in our model comprise a random vector, whose probability distribution is estimated independently by the leader and the follower, based on their own data.
To address the distributional uncertainty, we formulate a distributionally robust interdiction~(DRI) model, where both decision-makers solve conventional distributionally robust optimization problems based on the Wasserstein metric. For this model, we prove asymptotic consistency and derive a polynomial-size mixed-integer linear programming~(MILP) reformulation. Furthermore, in our bilevel optimization context, the leader may face uncertainty due to its incomplete knowledge of the follower's data. In this regard, we propose two distinct approximations of the true DRI model, where the leader has incomplete or no information about the follower's data. The first approach employs a pessimistic approximation, which turns out to be computationally challenging and requires a specialized reformulation amenable to a Benders-type decomposition algorithm. The second approach leverages a robust optimization approach from the leader's perspective. To address the resulting problem, we propose a scenario-based approximation that admits a potentially large single-level~MILP reformulation and satisfies asymptotic robustness guarantees. 
Finally, for a class of randomly generated instances of the packing interdiction problem, we evaluate numerically how the information asymmetry and the decision-makers' risk preferences affect the models' out-of-sample performance.
\end{abstract}

\begin{keyword} 
Bilevel optimization; Distributionally robust optimization; Interdiction; Mixed-integer programming; Wasserstein distance
\end{keyword}

\end{frontmatter}
\onehalfspacing 
\section{Introduction} \label{sec: intro}
Interdiction forms a broad class of deterministic and stochastic optimization problems, arising, for example, in the military, law-enforcement and infectious disease control contexts; see, e.g., the surveys in \citep{Smith2013, Smith2020} and the references therein. A classical \textit{attacker-defender} interdiction model can be viewed as a zero-sum game between two decision-makers: an upper-level decision-maker, commonly referred to as the \textit{leader} (or attacker), and a lower-level decision-maker, referred to as the \textit{follower}~(or defender). The leader initiates the game by allocating limited interdiction resources to target components that are crucial to the follower's operations. In response, the follower observes the leader's actions and aims to maximize its own profit within the impacted environment. Notably, the leader is informed about the follower's objective criterion and therefore chooses attacks that minimize the potential profit gained by the follower; see, e.g., the survey in \cite{Smith2020} and the references~therein. 

 In this study, we focus on a class of linear min-max interdiction problems of the form:
\begin{equation} \label{bilevel interdiction problem}
	\mbox{[\textbf{DI}]: \quad } \min_{\mathbf{x} \in X} \; \max_{\mathbf{y} \in Y(\mathbf{x}) } \;	\mathbf{c}^\top \mathbf{y},
\end{equation}	
where
\begin{equation} \label{eq: feasible sets}
	X = \big\{\mathbf{x} \in \{0, 1\}^m: \mathbf{H} \mathbf{x} \leq \mathbf{h}\big\} \; \, \mbox{ and } \; \, Y(\mathbf{x}) = \big\{\mathbf{y} \in \mathbb{R}_+^n: \mathbf{F}\mathbf{y} + \mathbf{L} \mathbf{x} \leq \mathbf{f}\big\}
\end{equation} are, respectively, the leader's and the follower's feasible sets. The leader's decision variables~$\mathbf{x}$ are restricted to be \textit{binary}, while the follower's decision variables~$\mathbf{y}$ are \textit{continuous}. To guarantee the existence of an optimal solution in [\textbf{DI}], we make the following standard assumption (see, e.g., \cite{Zare2019}):
\begin{itemize}
	\item[\textbf{A1}.] The feasible sets $X$ and $Y(\mathbf{x})$ for each $\mathbf{x} \in X$ are non-empty and bounded.
\end{itemize}

The vector $\mathbf{c} \in \mathbb{R}^n$ represents deterministic \textit{profits} (or \textit{costs}) associated with the follower's decision~$\mathbf{y}$. In particular, as there are no restrictions on the sign of $\mathbf{c}$, [\textbf{DI}] describes a broad class of interdiction problems, including, for example, the min-cost flow and the multicommodity flow interdiction~\cite{Lim2007, Smith2008,Smith2013}, packing interdiction~\cite{Dinitz2013}, as well as electric grid security problems~\cite{Salmeron2004}. It is well known that [\textbf{DI}] remains $NP$-hard even when the integrality constraints for $\mathbf{x}$ are relaxed~\cite{Hansen1992, Jeroslow1985}.


The aim of this study is to examine a \textit{stochastic} version of [\textbf{DI}], where the follower's objective function coefficients are subject to uncertainty. First, in line with the standard framework of data-driven stochastic optimization \cite{Bertsimas2018, Esfahani2018}, we assume that $\mathbf{c}$ is a random vector, whose probability distribution is not known \textit{a priori} and can only be observed through a finite training data set. Specifically, in our bilevel optimization context, both the leader and the follower are supposed to use their own data sets to estimate the distribution of~$\mathbf{c}$. Second, we argue that the leader in our model may have limited access to the follower’s private data.  For example, in practice, the follower may deliberately withhold its private data since the decision-makers in [\textbf{DI}] estimate the same profit vector $\mathbf{c}$ and pursue conflicting objectives. Motivated by these considerations, we formulate three distinct \textit{distributionally robust interdiction}~(DRI) models that address the distributional uncertainty and reflect different levels of the leader's knowledge about the follower's~data. 


\subsection{Related literature}
Uncertainty in interdiction models may arise due to various different factors. For instance, interdiction actions may fail with a certain probability, or key problem parameters, such as arc capacities or arc costs in network interdiction models, could be subject to uncertainty; see, e.g., the surveys in~\cite{Beck2023, Cormican1998, Smith2020} and the references therein. Furthermore, as argued in \cite{Beck2023}, the sources of uncertainty in bilevel optimization are more complex than those in standard single-level problems. That is, not only the problem's parameters may be uncertain but also decisions of one of the decision-makers may only be observed partially by the other decision-maker.

In this literature review, we categorize the most relevant studies according to their approach to uncertainty, with the focus on \textit{robust optimization}, \textit{stochastic programming}, and \textit{distributionally robust optimization} models. In particular, we discuss whether the leader's and the follower's decisions are made \textit{before} or \textit{after} the realization of uncertainty and, where appropriate, review the computational complexity and solution approaches for the resulting optimization models.

\looseness-1 \textbf{Robust interdiction models.} In robust optimization, uncertain problem parameters are commonly represented using deterministic \textit{uncertainty sets}, and the decision-maker chooses an optimal solution based on the worst-case possible realization of the uncertain parameters; see, e.g., \cite{Bental2009, Bental2002}. 
While robust optimization problems are extensively studied in the context of single-level optimization problems, relatively few studies consider a robust optimization approach in the context of interdiction or more general bilevel optimization problems; see, e.g., \cite{Beck2023} and the references therein. 

First, we refer to Beck~et~al.~\cite{Beck2023b}, who consider a min-max knapsack interdiction problem, where both the leader's and the follower's decision variables are integer, while the follower does not have full information about its objective function or constraint coefficients. The authors model the uncertainty using budget-constrained uncertainty sets \cite{Bertsimas2003} and propose a general branch-and-cut framework for solving the resulting problem. Notably, the leader in \cite{Beck2023b} has complete information about the follower's uncertainty set, whereas both decision-makers implement their decisions \textit{here-and-now}, prior to the realization of uncertainty.

Another relevant study is by Buchheim et al.~\cite{Buchheim2021}, which examines 
 mixed-integer bilevel linear programs with binary decision variables for the leader and continuous decision variables for the follower. Specifically, the only difference between [\textbf{DI}] and the model in \cite{Buchheim2021} is that the latter allows different objective function coefficients for the leader and the follower. Buchheim et al.~\cite{Buchheim2021} assume that the leader has incomplete information about the follower's objective function coefficients and estimates them using either interval or discrete uncertainty sets. As a result, the leader in \cite{Buchheim2021} makes a here-and-now decision, while the follower makes a \textit{wait-and-see} decision, affected by the realization of uncertainty. Buchheim et al.~\cite{Buchheim2021} show that, under interval uncertainty, the resulting robust bilevel problem is~$\Sigma^p_2$-hard. In other words, this problem is located at the second level of the polynomial hierarchy and, under some reasonable assumptions, there is no way of formulating it as a single-level mixed-integer linear programming~(MILP) problem of polynomial size; see, e.g.,~\citep{Jeroslow1985}. On a positive note, it is established in \cite{Buchheim2021} that, under discrete uncertainty, the problem 
is at most one level harder than the follower's~problem. 

Then, several recent studies consider robust formulations in the context of network interdiction problems. For example, Azizi and Seifi \cite{Azizi2024} consider the shortest path network interdiction problem, where the follower is uncertain about the additional costs imposed for traversing interdicted arcs. The authors consider here-and-now decision-makers, represent the uncertain parameters using a budget-constrained uncertainty set, and solve the resulting robust problem by using its single-level dual reformulation. In contrast, Chauhan~et~al.~\cite{Chauhan2024} study a robust maximum flow interdiction problem with a wait-and-see follower, considering uncertain arc capacities and interdiction resources. The problem is approximated by a single-level MILP reformulation and then addressed using three heuristics designed for this problem and its underlying network structure. Finally, we refer to the studies in \cite{Borrero2016, Borrero2019} for multi-stage interdiction models with incomplete information, in which the leader refines its information about the follower's costs and the structure of the follower's problem by observing the follower's~actions.


\textbf{Stochastic interdiction models.} In stochastic optimization, uncertain problem parameters are typically assumed to follow a fixed and known probability distribution, and the goal is to optimize a risk measure such as the expected value or the conditional value-at-risk (CVaR) based on this distribution; see, e.g.,~\cite{Birge2011, Rockafellar2000,Shapiro2021}. 
Most of the existing literature considers stochastic interdiction problems involving a \textit{risk-neutral} leader, whose objective is to minimize the follower's expected utility; 
see, e.g.,~\cite{Cormican1998, Janjarassuk2008, Nguyen2022} and the survey in \cite{Smith2020}. In contrast, 
we focus on \textit{risk-averse} stochastic interdiction models, where the decision-makers account for the risks associated with the realization of uncertain parameters and aim to mitigate the impact of extremely large losses.

Of particular relevance to our model is the study by Lei~et~al.~\cite{Lei2018}, which explores a stochastic maximum flow interdiction model with here-and-now decision-makers. In their framework, the arc capacities and the effects of interdiction actions follow a discrete probability distribution defined by a finite set of \textit{scenarios} and associated probabilities. Furthermore, the follower can introduce additional arc capacities for mitigating flow losses, after knowing the leader's interdiction plan but before realizing the uncertainty. The leader and the follower in \cite{Lei2018} seek to minimize and maximize, respectively, the right- and the left-tail CVaR; see, e.g.,~\cite{Rockafellar2000}. That is, since the decision-makers have conflicting objective functions, the pessimism of one decision-maker naturally translates to an optimistic viewpoint for the other. Based on the strong duality, Lei~et~al.~\cite{Lei2018} propose various MILP reformulations of the resulting problem, depending on the decision-makers' risk preferences. 

 Next, Song and Shen \cite{Song2016} and Pay~et~al.~\cite{Pay2019} consider a stochastic shortest path interdiction problem with uncertain travel times and interdiction effects. Both models assume a known discrete distribution of the uncertain problem parameters and focus on the case of a here-and-now leader and a wait-and-see follower. Song and Shen \cite{Song2016} introduce a chance constraint on the length of the follower's shortest path and address the resulting problem using a branch-and-cut algorithm. 
On the other hand, Pay~et~al.~\cite{Pay2019} assume ambiguous risk-preferences from the leader's perspective. In this regard, the authors provide several approaches to optimization over a class of utility functions. 

Finally, Atamt{\"u}rk~et~al.~\cite{Atamturk2020} consider mean-risk maximum flow interdiction problem with uncertain arc capacities. In contrast to the scenario-based approach, the model presented in~\cite{Atamturk2020} assumes that the uncertain problem parameters follow a nearly normal distribution, characterized by a known mean and a known covariance matrix. Then, the problem of minimizing the maximum flow-at-risk over a discrete set of actions is reduced to a sequence of discrete quadratic optimization problems. 

\textbf{Distributionally robust interdiction models.}
Distributionally robust optimization (DRO) models assume that the distribution of uncertain problem parameters is unknown and belongs to an \textit{ambiguity set} of plausible distributions; see, e.g.,~\cite{Delage2010,Esfahani2018,Wiesemann2014} and the review in \cite{Rahimian2022}. 
Then, similar to robust optimization, the decision-maker selects an optimal solution based on the worst-case distribution of the uncertain parameters. Overall, one-stage DRO models can often be reformulated as tractable convex optimization problems and provide, on average, a better out-of-sample performance compared to stochastic programming models constrained to a single candidate distribution \cite{Esfahani2018, Wiesemann2014}.

To the best of our knowledge, there are only a few studies that examine a DRO approach in the context of interdiction problems. Notably, Sadana and Delage \cite{Sadana2023} analyze the value of randomization for a class of distributionally robust risk-averse maximum flow interdiction problems. Specifically, it is assumed that the arc capacities are subject to distributional uncertainty and the ambiguity set is comprised of \textit{discrete} distributions concentrated around a fixed reference distribution. The leader in~\cite{Sadana2023} randomizes over the feasible interdiction plans in order to minimize the worst-case CVaR of the maximum flow with respect to both the unknown distribution of the capacities and its own randomized strategy. Meanwhile, the follower makes a wait-and-see decision by solving a maximum flow problem with given arc capacities. Sadana and Delage~\cite{Sadana2023} reformulate the leader's DRO problem as a bilinear optimization problem and solve it using a specialized spatial branch-and-bound algorithm. 

In conclusion, we refer to Kang and Bansal \cite{Kang2025}, who investigate multistage \textit{risk-neutral} stochastic mixed-integer programs. In the context of network interdiction, the authors account for uncertainty in arc capacities and interdiction effects, considering a leader who may be either ambiguity-averse or ambiguity-receptive. In these two variations of the problem, the leader focuses, respectively, on the worst-case and the best-case expected losses incurred by the follower. Similar to the work of Sadana~et~al.~\cite{Sadana2023}, both formulations involve a wait-and-see follower but allow for discrepancy-based ambiguity sets based on the Wasserstein metric. To solve both problem versions, Kang and Bansal~\cite{Kang2025} propose several exact solution methods that leverage decomposition and reformulation~techniques.

\subsection{Our approach and contributions} \label{subsec: approach and contributions}
In this paper, we analyze a \textit{data-driven distributionally robust interdiction} (DRI) problem associated with [\textbf{DI}], where the profit vector $\mathbf{c}$ is subject to distributional uncertainty. That is, the leader and the follower in our model are given two distinct training data sets drawn from the unknown \textit{true distribution} of~$\mathbf{c}$. 
 We first introduce a \textit{basic DRI model}, where the leader has \textit{complete knowledge} of the follower's private data and the only source of uncertainty is associated with the distribution~of~$\mathbf{c}$. Notably, as both data sets are used to estimate the \textit{same} probability~distribution, this model requires the assumption that the follower's data set is included in the leader's. 

To address the distributional uncertainty, we employ data-driven ambiguity sets based on the Wasserstein distance from the empirical distribution of the data; see, e.g.,~\cite{Esfahani2018, Gao2023}. 
 Thus, similar to Lei~et~al.~\cite{Lei2018}, both decision-makers in our model solve conventional one-stage Wasserstein DRO problems based on the conditional value-at-risk (CVaR) and act \textit{here-and-now}, prior to the realization of~uncertainty.
Our choice of the ambiguity sets 
is justified by the fact that one-stage Wasserstein DRO models effectively balance tractability and statistical guarantees, leading to their broad applications in machine learning, portfolio optimization and network routing problems; see,~e.g.,~\cite{Blanchet2022, Kuhn2019, Wang2020}. 

 Next, following the major part of the bilevel optimization literature, we consider an \textit{optimistic} formulation of the resulting bilevel problem. That is, if the follower has multiple optimal solutions, then the one most favorable to the leader is selected; see, e.g., the surveys in~\cite{Colson2007, Kleinert2021}. The justification for the optimistic formulation is two-fold. First, \textit{pessimistic} formulations of mixed-integer bilevel linear programs, while arguably more realistic, typically require specialized solution techniques; see, e.g., \cite{Wiesemann2013, Zeng2020}. Second, as we establish later in Section \ref{sec: basic}, the follower’s problem in our basic DRI model admits a linear programming (LP) reformulation whose constraint coefficients depend on the training data. Therefore, under mild regularity conditions, the follower’s problem has a unique optimal solution, implying that the optimistic and the pessimistic models coincide.  

As outlined earlier, the inherent conflict between the decision-makers in [\textbf{DI}] suggests that the leader may also face uncertainty in the follower's private data and, consequently, in the follower's optimal decision. From a practical perspective, similar information asymmetry arises in Stackelberg-type coordination settings, where each agent optimizes its own operational or economic objective based on its proprietary or privacy-sensitive data, and the follower may only disclose a sanitized or partial version of this information; see, e.g., the privacy-preserving Stackelberg mechanisms for sequential energy-market coordination in \cite{Fioretto2020, Mitridati2022}. 
More generally, bilevel models with incomplete information explicitly treat parts of the follower’s problem data as unknown to the leader and motivate learning/estimation from observed follower reactions; see, e.g., \cite{Borrero2019, Borrero2022}. 

 To address the leader's partial knowledge of the follower's data, we develop two alternative optimization models.
For these models, 
we first introduce a hypothetical \textit{true basic DRI model}, in which the leader has full knowledge of the follower's data set, but this set is not a subset of the leader's. We~then propose two approximations of the true basic model. In the \textit{pessimistic approximation}, the leader disregards any available information about the follower's data and optimizes its own objective function with respect to the worst-case feasible follower's policy. On the other hand, in the \textit{semi-pessimistic approximation}, we 
assume that the follower's data set is subject to component-wise \textit{interval uncertainty}. Then, similar to the one-stage DRO model with data uncertainty in \cite{Ketkov2024}, the leader optimizes its objective function assuming the worst-case possible realization of the follower's data within the uncertainty set. 

Our main theoretical contributions can be summarized as follows:
\begin{itemize} 
	\setlength{\itemsep}{0.05cm} 
	\setlength{\parskip}{0.05cm} 
	\setlength{\topsep}{2pt} 
	
	\item First, provided that the Wasserstein ambiguity sets for both decision-makers are defined in terms of $\ell_1$- or $\ell_\infty$-norm, we establish the computational complexity for all proposed optimization models. That is, while the basic DRI model is $NP$-hard, both pessimistic and semi-pessimistic approximations turn out to be $\Sigma_2^p$-hard. To prove our complexity results, we use reductions from the dominating set interdiction and the robust bilevel interdiction problem with interval uncertainty that are known to be $\Sigma_2^p$-hard \cite{Buchheim2021, Grune2025}.
	
	\item In line with the one-stage DRO models in \cite{Esfahani2018, Gao2023}, we demonstrate that, under mild assumptions, the basic DRI model is \textit{asymptotically consistent}. That is, when the decision-makers acquire more data, their optimal solutions and the respective optimal objective function values converge, in a sense, to those of the underlying stochastic programming problem. This is, to our knowledge, the first result establishing convergence of DRO models in a bilevel interdiction problem setting. 
	\item The basic DRI model is reformulated and solved as a single-level MILP problem of polynomial size. This reformulation builds on linear programming reformulations of the worst-case CVaR problems for each decision-maker~\cite{Esfahani2018} and a strong duality-based reformulation of the follower's problem; see, e.g.,~\cite{Zare2019}. 
	\item To address the pessimistic approximation, we design a Benders-type decomposition algorithm tailored to two special cases of the problem, where the leader is either \textit{risk-neutral} or \textit{ambiguity-free}. The algorithm for each case is based on the standard decomposition techniques for bilevel optimization~\cite{Israeli2002, Zeng2014} and leverages a disjoint bilinear structure of the inner optimization problem. 
	\item To address the semi-pessimistic approximation, we propose a scenario-based discretization of the leader's uncertainty set. This approach, as indicated in~\cite{Buchheim2021}, reduces the computational complexity and enables a unified MILP reformulation of the problem. To justify our approach, we show that the discretized semi-pessimistic approximation is almost surely \textit{robust} with respect to the follower's data in an asymptotic sense, i.e., when the number of scenarios tends to infinity. 
\end{itemize}

Despite the fact that some standard no-good cuts can be applied to general bilevel interdiction problems with a non-convex lower level \cite{Bomze2025}, $\Sigma_2^p$-hard problems remain notoriously difficult, even for medium-sized problem instances. For this reason, the focus of our study is on solution techniques tailored to the structure of the proposed approximations. 

\begin{table}[ht]
	\centering \footnotesize \doublespacing 
	\begin{tabular}{c c |c c c}
		& & Our paper & Lei~et~al.~\cite{Lei2018} & Sadana \& Delage \cite{Sadana2023} \\\hline
		
		\multirow{2}{*}{Follower} & here-and-now & \checkmark & \checkmark & \xmark \\
		& wait-and-see & \xmark & \xmark & \checkmark 
		\\
		\hline
		\multirow{2}{*}{Follower's decision variables} & continuous & \checkmark & \checkmark & \checkmark \\
		& integer & \xmark & \xmark & \xmark /\checkmark 
		\\
		\hline
		\multirow{4}{*}{Uncertainty} & scenario-based & \xmark & \checkmark & \checkmark \\ 
		& data-driven & \checkmark & \xmark & \xmark \\
		& objective function & \checkmark & \checkmark & \checkmark \\
		& constraints & \xmark & \checkmark & \checkmark \\\hline
		\multirow{2}{*}{Risk-aversion} & leader & \checkmark & \checkmark & \checkmark \\
		& follower & \checkmark & \checkmark & \xmark \\\hline
		\multirow{2}{*}{Ambiguity-aversion} & leader & \checkmark & \xmark & \checkmark \\
		& follower & \checkmark & \xmark & \xmark \\\hline 
		\multicolumn{2}{c|}{\multirow{2}{*}{Information asymmetry}} & \multirow{2}{*}{\checkmark} & \multirow{2}{*}{\xmark} & \multirow{2}{*}{\xmark} \\
		& & & & \\\hline
	\end{tabular}
	\caption{Comparison of our model with the related models in \cite{Lei2018} and \cite{Sadana2023}. Symbol \xmark /\checkmark indicates that the model in \cite{Sadana2023} can only be applied to the follower's problem with a unimodular constraint matrix. }
	\label{tab: contributions}
\end{table} 

\textbf{Contributions to the related literature.} In view of the discussion above, we compare our DRI model with the most related studies in Lei~et~al.~\cite{Lei2018} and Sadana and Delage \cite{Sadana2023}; see Table \ref{tab: contributions}. Our contributions can be summarized as follows:
\begin{itemize} \setlength{\itemsep}{0.05cm} 
\setlength{\parskip}{0.05cm} 
\setlength{\topsep}{2pt} 
	 \item To the best of our knowledge, our model is the first in the interdiction literature to incorporate data-driven distributional uncertainty and information asymmetry with respect to the data. 
\item Unlike the studies in~\cite{Lei2018, Sadana2023}, we incorporate risk aversion and ambiguity aversion for both decision-makers, thereby adopting the standard conventions of one-stage DRO models.
\item Beyond the standard distributional uncertainty arising from the decision-makers' incomplete knowledge of the true data-generating distribution, our model also examines uncertainty in the follower's decision caused by the leader's limited knowledge of the follower's data.
\end{itemize}

\looseness-1 Given these advantages, our model also has two major limitations. First, unlike existing interdiction models with a wait-and-see follower \cite{Pay2019, Sadana2023, Song2016}, our approach cannot address 
 interdiction problems with a unimodular structure at the lower level \cite{Israeli2002, Zenklusen2010}. 
Instead, we focus on a more specific class of interdiction problems, such as the min-cost flow or the multicommodity flow interdiction problems~\cite{Lim2007}, as well as the packing interdiction problem \cite{Dinitz2013}, where the follower's continuous decisions have a natural interpretation; see row 2 of Table \ref{tab: contributions}. Second, we account for uncertainty in the follower's objective function coefficients but not in the follower's \textit{constraints}; see row 3 of Table \ref{tab: contributions}. In this regard, we observe that distributionally robust chance-constrained problems with a Wasserstein ambiguity set do not admit LP reformulations \cite{Chen2024}, which makes even the basic DRI model computationally challenging. 
On a positive note, with a slight modification, our model can also address uncertainty in the interdiction effects, since this uncertainty appears in the follower's objective~function. 

The remainder of the paper is organized as follows. First, in Section \ref{sec: problem}, we introduce a full information stochastic programming model and our modeling assumptions. Next, in Sections \ref{sec: basic}--\ref{sec: semi-pessimistic}, we analyze the basic, pessimistic and semi-pessimistic DRI models, respectively, along with their theoretical properties and corresponding solution methods. In Section \ref{sec: comp study}, we perform a numerical analysis of the proposed DRI models using a set of randomly generated instances of the packing interdiction problem. Finally, Section \ref{sec: conclusion} provides conclusions and outlines possible directions for future research. 

\textbf{Notation.} All vectors and matrices are denoted by bold letters. A vector of all ones is represented as~$\mathbf{1}$ and a diagonal matrix with vector $\mathbf{x}$ on the diagonal is referred to as $\diag(\mathbf{x})$. Also, $\dim(\mathbf{x})$ refers to the dimension of $\mathbf{x}$. 
We use $\mathbb{R}_+$ to denote the set of nonnegative real numbers and $\mathbb{Z}_+$~($\mathbb{N}$) to denote the set of nonnegative (positive) integers. Additionally, subscripts ``$l$'' and~``$f$'' refer to parameters of the leader's and the follower's problems, respectively. For $a \in \mathbb{R}$, we use a standard notation $a^+ = \max\{a; 0\}$ and $a^- = \min\{a; 0\}$. For a given support set $S \subseteq \mathbb{R}^n$, let $\mathcal{Q}_0(S)$ denote the space of all probability distributions defined on $S$. Finally, the $k$-fold product of a distribution $\mathbb{Q} \in \mathcal{Q}_0(S)$ is denoted by $\mathbb{Q}^k$, which represents a distribution on the Cartesian product space $S^k$. 

\section{Model setup and assumptions} \label{sec: problem}
\textbf{Full information model.}
Let $\mathbb{Q}^* \in \mathcal{Q}_0(S)$ represent a true (nominal) distribution of the profit vector $\mathbf{c}$, supported on $S \subseteq \mathbb{R}^n$. Ideally, with complete knowledge of $\mathbb{Q}^*$, our goal is to solve the following stochastic programming formulation of [\textbf{DI}]:
\begin{subequations} \label{stochastic programming problem}
\begin{align} 
 \mbox{[\textbf{SI}]: \quad } z^* :=	& \min_{\mathbf{x}, \mathbf{y}} \; \rho_l(\mathbf{y}, \mathbb{Q}^*) \label{obj: stochastic programming problem}\\ 
 \mbox{s.t. } & \mathbf{x} \in X \\
	& \mathbf{y} \in \argmax_{\,\tilde{\mathbf{y}} \in Y(\mathbf{x})} \;\rho_f(\tilde{\mathbf{y}}, \mathbb{Q}^*),
\end{align}	
\end{subequations}
where $\rho_l(\mathbf{y}, \mathbb{Q}^*)$ and $\rho_f(\mathbf{y}, \mathbb{Q}^*)$ denote the leader's and the follower's risk functions, and the minimization with respect to $\mathbf{y}$ in (\ref{obj: stochastic programming problem}) implies an optimistic version of the bilevel problem; recall our discussion in Section \ref{subsec: approach and contributions}. 

Following the model of Lei~et~al.~\cite{Lei2018}, we define~$\rho_l(\mathbf{y}, \mathbb{Q}^*)$ and $\rho_f(\mathbf{y}, \mathbb{Q}^*)$ using the right-tail and the left-tail conditional values-at-risk (CVaRs), respectively. That is, given confidence levels $\alpha_l \in (0,1)$ and $\alpha_f \in (0,1)$, we have: 
\begin{subequations} \label{eq: risk functions} 
	\begin{align}
		& \rho_l(\mathbf{y}, \mathbb{Q}^*) := \min_{t_l \in \mathbb{R}} \Big\{ t_l + \frac{1}{1 - \alpha_l} \mathbb{E}_{\mathbb{Q}^*} \big\{ (\mathbf{c}^\top \mathbf{y} - t_l)^+ \big\} \Big\} \; \mbox{ and }\label{eq: right-tail CVaR} \\
		& \rho_f(\mathbf{y}, \mathbb{Q}^*) := \max_{t_f \in \mathbb{R}} \Big\{ t_f + \frac{1}{1 - \alpha_f} \mathbb{E}_{\mathbb{Q}^*} \big\{ (\mathbf{c}^\top \mathbf{y} - t_f)^- \big\} \Big\}, \label{eq: left-tail CVaR}
	\end{align}
\end{subequations}
where an equivalent representation of CVaR due to \cite{Rockafellar2000} 
is employed. In this context, (\ref{eq: right-tail CVaR}) and (\ref{eq: left-tail CVaR}) measure the expected value of extremely high losses for the leader (i.e., upper-tail risk) and extremely low profits for the follower (i.e., lower-tail risk), respectively. 
Notably, both risk measures reduce to the expected value of $\mathbf{c}^\top \mathbf{y}$ under~$\mathbb{Q}^*$ in the limiting cases of $\alpha_l = 0$ and~$\alpha_f = 0$. 
Consequently, by choosing appropriate confidence levels, $\alpha_l$ and $\alpha_f$, the decision-makers can adjust their level of risk-aversion.


\textbf{Distributional uncertainty.} To address the distributional uncertainty regarding $\mathbb{Q}^*$, we make the following assumptions about the partial distributional information available to the decision-makers:
\begin{itemize}
	\item[\textbf{A2.}] The \textit{support} of $\mathbb{Q}^*$ is known and given by a non-empty compact polyhedral set
	\begin{equation} \label{eq: support set}
		S = \Big\{\mathbf{c} \in \mathbb{R}^n: \mathbf{B} \mathbf{c} \leq \mathbf{b} \Big\}.
	\end{equation}
	\item[\textbf{A3.}] The leader and the follower have access to two independent and identically distributed (i.i.d.) training data sets, generated according to~$\mathbb{Q}^*$, 
	\begin{equation} \label{eq: data sets}
		\hat{\mathbf{C}}_l = \Big\{ \hat{\mathbf{c}}_l^{\text{\tiny (\textit{k})}} \in S, \; k \in K_l := \{1, \ldots, k_l\}\Big\} \; \mbox{ and } \; \hat{\mathbf{C}}_f = \Big\{ \hat{\mathbf{c}}_f^{\text{\tiny (\textit{k})}} \in S, \; k \in K_f := \{1, \ldots, k_f\} \Big\}.
	\end{equation}
	Furthermore, $\hat{\mathbf{C}}_f \subseteq \hat{\mathbf{C}}_l$ or, equivalently, $k_f \leq k_l$ and
	$\hat{\mathbf{c}}_f^{\text{\tiny (\textit{k})}} = \hat{\mathbf{c}}_l^{\text{\tiny (\textit{k})}}$ for each $k \in K_f$.	
\end{itemize}

 Both Assumptions \textbf{A2} and \textbf{A3} align with the standard one-stage DRO models in \cite{Delage2010, Esfahani2018, Wiesemann2014} and offer a natural, myopic extension to our bilevel interdiction problem setting. Notably, Assumption~\textbf{A2} also requires compactness of the support set $S$. 
This assumption simplifies our convergence analysis and naturally fits interdiction problems, where uncertain problem parameters usually come from a finite set of scenarios~\cite{Lei2018, Sadana2023} or satisfy specified interval constraints~\cite{Borrero2016, Buchheim2021}. Finally, the second part of Assumption~\textbf{A3} indicates that solving the leader’s problem requires full knowledge of the follower’s private data;  recall our discussion of the basic DRI model in Section \ref{subsec: approach and contributions}. 




Next, following conventional Wasserstein DRO models in \cite{Esfahani2018, Gao2023}, we define each decision-maker's ambiguity set as a type-1 Wasserstein ball, centered at their empirical data distribution. Specifically, for $i \in \{l, f\}$, the ambiguity set is given by: 
\begin{equation} \label{ambiguity set}
	\mathcal{Q}_i := \left\{ \mathbb{Q} \in \mathcal{Q}_0(S) : W^p(\hat{\mathbb{Q}}_i, \mathbb{Q}) \leq \varepsilon_i \right\},
\end{equation}
where 
\begin{equation} \label{eq: empirical distribution} \nonumber
	\hat{\mathbb{Q}}_i (\hat{\mathbf{C}}_i) := \frac{1}{k_i}\sum_{k \in K_i} \; \delta(\hat{\mathbf{c}}_i^{\text{\tiny (\textit{k})}}), \quad i \in \{l, f\},
\end{equation}
is the empirical distribution of the observed data, with $\delta: \mathbb{R}^{n} \rightarrow \mathbb{R}$ representing the Dirac delta function, and $\varepsilon_i \in \mathbb{R}_+$ is the Wasserstein radius. Furthermore, given that $\Pi(\hat{\mathbb{Q}}_i, \mathbb{Q})$ is the set of all couplings with marginals $\hat{\mathbb{Q}}_i$ and $\mathbb{Q}$, 
\begin{equation} \nonumber
	W^p(\hat{\mathbb{Q}}_i, \mathbb{Q}):= \inf_{\pi \in \Pi(\hat{\mathbb{Q}}_i, \mathbb{Q})} \int_{S \times S} \Vert \hat{\mathbf{c}}_i - \mathbf{c} \Vert_p \pi( \mathrm{d}\hat{\mathbf{c}}_i, \mathrm{d}\mathbf{c}), \quad i \in \{l, f\},
\end{equation} 
denotes the type-1 Wasserstein distance between $\hat{\mathbb{Q}}_i$ and $\mathbb{Q}$ with respect to the $\ell_p$-norm \cite{Kantorovich1958}. 

 In the following, we focus on the Wasserstein distance with respect to $\ell_1$- or~$\ell_{\infty}$-norm, since computing the worst-case CVaR within such Wasserstein balls is equivalent to solving a finite linear programming~(LP) problem; see Corollary~5.1 in~\cite{Esfahani2018}. Finally, $\varepsilon_{i}$, $i \in \{l, f\}$, are typically selected to ensure that, with high probability, the true distribution~$\mathbb{Q}^*$ belongs to the ambiguity sets~(\ref{ambiguity set}); see,~e.g.,~\cite{Esfahani2018, Gao2023}. We discuss our choice of the Wasserstein radii in more detail within~Section \ref{sec: comp study}. 

To motivate our basic DRI model, as well as the pessimistic and the semi-pessimistic approximations discussed in Section \ref{subsec: approach and contributions}, we present the following illustrative example.

		
		
		

\begin{table}[h!]
	\centering
	\doublespacing
	\scriptsize
	\begin{tabular}{cccc} 
		\hline 	
		Model & Optimal leader's decision & Leader's estimated objective & Leader's true objective \\\hline
		Full information & $(0, 0, 1, 1)^\top$ & 6 & 6 \\\
		True basic & $(0,0,0,1)^\top$ & 7 & 9 \\
		Pessimistic & $(1,1,0,0)^\top$ & 8 & 10 \\
		Semi-pessimistic & $(0,0,0,1)^\top$ & 7 & 9 \\
		Augmented basic & $(1,1,0,0)^\top$ & 8 & 10 \\
		\hline
	\end{tabular}
	\caption{ Optimal solutions and the respective leader's estimated and true objective function values under optimization models (\ref{ex: full information model})--(\ref{ex: basic model}), respectively, discussed in Example~\ref{ex: example 1}.}
	\label{tab: example 1}
\end{table}

\begin{example} \upshape \label{ex: example 1}
	Consider an instance of the packing interdiction problem \cite{Dinitz2013} with uncertain profits given~by:
	\[\min_{\mathbf{x} \in X} \max_{\;\mathbf{y} \in Y(\mathbf{x})} \mathbf{c}^\top \mathbf{y},\]
	where
	\[X = \Big\{\mathbf{x} \in \{0, 1\}^4: \, \sum_{i = 1}^4 x_i \leq 2 \Big\} \, \, \mbox{ and } \, \, Y(\mathbf{x}) = \Big\{\mathbf{y} \in [0, 1]^4: \, \mathbf{y} \leq \mathbf{1} - \mathbf{x}, \, \sum_{i = 1}^4 y_i \leq 1 \Big\}\]
	are, respectively, the leader's and the follower's feasible sets. Assume that the follower's profit vector~$\mathbf{c}$ belongs to a support set 
	\[S = \Big\{\mathbf{c} \in \mathbb{R}^4: \, c_1 \in [0, 10], \, c_2 \in [1, 11], \, c_3 \in [6, 12], \, c_4 \in [7, 13] \Big\},\]
	with the true distribution $\mathbb{Q}^*$ being a \textit{discrete uniform distribution} on $S$. 

	Given that both decision-makers are risk-neutral, i.e., $\alpha_l = 0$ and $\alpha_f = 0$, the full information model [\textbf{SI}] reads as:
	\begin{equation} \label{ex: full information model}
	z^* = \min_{\mathbf{x} \in X} \max_{\;\mathbf{y} \in Y(\mathbf{x})} \; \bar{\mathbf{c}}^{*\top} \mathbf{y}, 
	\end{equation}
where $\bar{\mathbf{c}}^{*} = (5, 6, 9, 10)^\top$ is a vector of the true expected costs.
The optimal solution of (\ref{ex: full information model}) is given by $\mathbf{x}^* = (0, 0, 1, 1)^\top$, 
with the true expected profit of the follower $z^* = 6$; see Table \ref{tab: example 1}. 

Next, for simplicity, we assume that $\varepsilon_l = \varepsilon_f = 0$, i.e., both the leader and the follower implement 
a myopic \textit{sample average approximation} based on $k_l = k_f = 5$ samples given by: 
\[
\renewcommand{\arraystretch}{0.8}
\hat{\mathbf{C}}_l = \left(
\begin{array}{cccc}
	10 & 11 & 8 & 7 \\
	6 & 7 & 7 & 9 \\
	10 & 10 & 7 & 9 \\
	9 & 11 & 6 & 7 \\
	10 & 11 & 7 & 8
\end{array}\right) \, \, \mbox{ and } \, \, \hat{\mathbf{C}}_f = \left(
\begin{array}{cccc}
3 & 3 & 8 & 13 \\
1 & 2 & 9 & 10 \\
3 & 11 & 11 & 13 \\
10 & 3 & 10 & 12 \\
4 & 1 & 12 & 7
\end{array}
\right),
\]
with average profits \[\bar{\mathbf{c}}_l := (9, 10, 7, 8)^\top \, \, \mbox{ and } \, \, \bar{\mathbf{c}}_f := (4, 5, 10, 11)^\top,\] respectively. 
Clearly, in this setting, the second part of Assumption \textbf{A3} does not hold, i.e., $\hat{\mathbf{C}}_f \nsubseteq \hat{\mathbf{C}}_l$.

Therefore, we first formulate a hypothetical \textit{true basic model} 
	\begin{subequations} \label{ex: true basic model} 
	\begin{align}
		z_b^* := & \min_{\mathbf{x}, \mathbf{y}} \; \bar{\mathbf{c}}_l^{\top} \mathbf{y} \\
		\mbox{s.t. } & \mathbf{x} \in X \\
		& \mathbf{y} \in \argmax_{\,\tilde{\mathbf{y}} \in Y(\mathbf{x})} \; \bar{\mathbf{c}}_f^\top \tilde{\mathbf{y}},
	\end{align}
\end{subequations} 
where the leader correctly estimates the follower's optimal decision. 
To address uncertainty in the follower's data, we consider the following three possible approaches:
	\begin{itemize}	
	\setlength{\itemsep}{0.05cm} 
	\setlength{\parskip}{0.05cm} 
	\setlength{\topsep}{2pt}	
	\item[(\textit{i})] The leader may solve a \textit{pessimistic approximation} of (\ref{ex: true basic model}) given by:
	\begin{equation} 
		\label{ex: pessimistic approximation} \hat{z}^*_p = \min_{\mathbf{x} \in X} \max_{\mathbf{y} \in Y(\mathbf{x})} \; \bar{\mathbf{c}}_l^{\top} \mathbf{y}. 
	\end{equation}
	The leader in (\ref{ex: pessimistic approximation}) disregards any information about the follower's data and minimizes its objective function value assuming the worst-case feasible follower's policy. 
	
	\item[(\textit{ii})] Alternatively, assume that the leader has partial information about the follower's data, e.g., it is aware of the first two components in $\hat{\mathbf{C}}_f$. In this case, the leader may solve a \textit{semi-pessimistic approximation} of (\ref{ex: true basic model}) given~by:
	\begin{subequations} \label{ex: semi-pessimistic approximation} 
		\begin{align} \hat{z}^*_{sp} := \min_{\mathbf{x} \in X} \; \max_{\bar{\mathbf{c}}_f \in \bar{S}_l} \; & \min_{\mathbf{y}}
			\; \bar{\mathbf{c}}_l^{\top} \mathbf{y} \\ 
			\mbox{s.t. } & \mathbf{y} \in \argmax_{\,\tilde{\mathbf{y}} \in Y(\mathbf{x})} \; \bar{\mathbf{c}}_f^\top \tilde{\mathbf{y} }, 
		\end{align}
	\end{subequations}
	where \[\bar{S}_l = \Big\{\bar{\mathbf{c}} \in \mathbb{R}^4: \; \bar{c}_1 = 4, \, \bar{c}_2 = 5,\, \bar{c}_3 \in [6, 12], \bar{c}_4 \in [7, 13] \Big\}. \] 
 In other words, based on the partial information about $\hat{\mathbf{C}}_f$, the leader in (\ref{ex: semi-pessimistic approximation}) may construct the uncertainty set $\bar{S}_l$ for the follower's average profits $\bar{\mathbf{c}}_f$ and account for the worst-case scenario within this set. 
	\item[(\textit{iii})] To enforce Assumption \textbf{A3}, one may assume that the leader generates the missing samples, e.g., by resampling from its own data set. In our case, since the decision-makers' data sets have the same size, we assume for simplicity that the leader's data set is used \textit{instead} of the follower's data set. 
	In other words, the leader may solve an \textit{augmented basic model} of the form:
	\begin{subequations} \label{ex: basic model} 
		\begin{align}
			\hat{z}_{ab}^* := & \min_{\mathbf{x}, \mathbf{y}} \; \bar{\mathbf{c}}_l^{\top} \mathbf{y} \\
			\mbox{s.t. } & \mathbf{x} \in X \\
			& \mathbf{y} \in \argmax_{\,\tilde{\mathbf{y}} \in Y(\mathbf{x})} \; \bar{\mathbf{c}}_l^{\top} \tilde{\mathbf{y}}.
		\end{align}
	\end{subequations}
	Within this example, (\ref{ex: basic model}) coincides with the pessimistic approximation (\ref{ex: pessimistic approximation}). 
	However, in general, these are two distinct models, as the decision-makers' objective functions may~vary.
\end{itemize}
 For each of the outlined models, we provide the leader's optimal solution along with the corresponding estimated and true objective function values in Table \ref{tab: example 1}. For all models except the full information model, the true value is computed under the true follower's policy, as defined in (\ref{ex: true basic model}), and the true~distribution~$\mathbb{Q}^*$. \hfill$\square$ 
\end{example}

We draw several important observations from Example \ref{ex: example 1}. First, both the leader and the follower may experience out-of-sample errors due to inaccuracies in their perception of the profit vector $\mathbf{c}$. For example, replacing the expected profit vector $\bar{\mathbf{c}}^*$ with the leader's estimate $\bar{\mathbf{c}}_l$ in the true basic model~(\ref{ex: true basic model}) increases the follower's true expected profit from $6$ to $9$. Another important source of error arises from the leader's incorrect estimate of the follower's policy. For instance, the augmented basic model~(\ref{ex: basic model}) and the pessimistic approximation~(\ref{ex: pessimistic approximation}) further increase the true follower's profit from~$9$ to $10$, which is not the case for the semi-pessimistic approximation (\ref{ex: semi-pessimistic approximation}). 
To provide additional practical insights, we investigate the effect of these errors numerically in Section~\ref{sec: comp study}.

Second, we observe that the optimal objective function values in (\ref{ex: true basic model}), (\ref{ex: pessimistic approximation}) and (\ref{ex: semi-pessimistic approximation}) satisfy the following hierarchy:
\begin{equation} \label{eq: hierarchy of bounds} \nonumber
z^*_{b} \leq \hat{z}^*_{sp} \leq \hat{z}^*_{p}.
\end{equation}
Specifically, in~Example \ref{ex: example 1}, we have $z^*_{b} = \hat{z}^*_{sp} = 7$ and $\hat{z}^*_{p} = 8$. The inequality $\hat{z}^*_{sp} < \hat{z}^*_{p}$ suggests that, even with partial information about the follower's data, the leader may potentially refine its estimate of the true follower's policy and thereby improve its in-sample performance. 
Finally, the optimal objective function value of (\ref{ex: basic model}), $\hat{z}^*_{ab}$, does not have any specific relationship with $z^*_b$, and therefore the augmented basic model (\ref{ex: basic model}) is only discussed within our numerical study in Section \ref{sec: comp study}. 

\section{Basic DRI model} \label{sec: basic}
\subsection{Problem formulation and solution approach}
In contrast to the full information model [\textbf{SI}], we assume that the true distribution $\mathbb{Q}^*$ is known to the leader and the follower only partially, as specified by Assumptions \textbf{A2} and \textbf{A3}. 
Then, the basic distributionally robust interdiction (DRI) model is formulated as: 
\begin{subequations} \label{bilevel distributionally robust problem} 
	\begin{align}
		\mbox{[\textbf{DRI}]: \quad } \hat{z}^*_b := & \min_{\mathbf{x}, \mathbf{y}} \; \Big\{ \max_{\mathbb{Q}_l \in \mathcal{Q}_l} \;\rho_l(\mathbf{y}, \mathbb{Q}_l) \Big\}\label{obj: bilevel DR}\\
		\mbox{s.t. } & \mathbf{x} \in X \label{cons: bilevel DR 1} \\
		& \mathbf{y} \in \argmax_{\,\tilde{\mathbf{y}} \in Y(\mathbf{x})} \; \Big\{ \min_{\mathbb{Q}_f \in \mathcal{Q}_f} \rho_f(\tilde{\mathbf{y}}, \mathbb{Q}_f) \Big\}, \label{cons: bilevel DR 2}
	\end{align}
\end{subequations}
where the data-driven ambiguity sets $\mathcal{Q}_l$ and $\mathcal{Q}_f$ are defined by equation (\ref{ambiguity set}). Additionally, the risk functions $\rho_l(\mathbf{y}, \mathbb{Q}_l)$ and $\rho_f(\mathbf{y}, \mathbb{Q}_f)$ are defined similar to equations (\ref{eq: right-tail CVaR}) and (\ref{eq: left-tail CVaR}), i.e.,
\begin{subequations} \label{eq: risk functions basic} 
	\begin{align}
		& \rho_l(\mathbf{y}, \mathbb{Q}_l) := \min_{t_l \in \mathbb{R}} \Big\{ t_l + \frac{1}{1 - \alpha_l} \mathbb{E}_{\mathbb{Q}_l} \big\{ (\mathbf{c}^\top \mathbf{y} - t_l)^+ \big\} \Big\} \; \mbox{ and }\label{eq: right-tail CVaR basic} \\
		& \rho_f(\mathbf{y}, \mathbb{Q}_f) := \max_{t_f \in \mathbb{R}} \Big\{ t_f + \frac{1}{1 - \alpha_f} \mathbb{E}_{\mathbb{Q}_f} \big\{ (\mathbf{c}^\top \mathbf{y} - t_f)^- \big\} \Big\}, \label{eq: left-tail CVaR basic}
	\end{align}
\end{subequations}
where $\alpha_l \in (0, 1)$ and $\alpha_f \in (0, 1)$ are the prespecified confidence levels. 

Importantly, the compactness of the support set $S$ (Assumption \textbf{A2}) ensures the existence of the worst-case distributions in (\ref{obj: bilevel DR}) and~(\ref{cons: bilevel DR 2}), allowing the supremum (infimum) over the ambiguity sets to be replaced by the maximum~(minimum); see Corollary 4.6 in \cite{Esfahani2018}. Moreover, we implicitly assume that the leader in [\textbf{DRI}] has full information about the follower's feasible set $Y(\mathbf{x})$ and two global parameters, $\alpha_f$ and $\varepsilon_{f}$; possible relaxations of this assumption are briefly discussed in Section~\ref{subsec: summary}. 

The following results establish a polynomial-size single-level MILP reformulation of [\textbf{DRI}]. 
\begin{lemma} \label{lemma 1}
	Let Assumptions \textbf{A1}-\textbf{A3} hold, and the ambiguity sets in \upshape (\ref{ambiguity set}) \itshape be defined in terms of $p$-norm with $p \in \{1, \infty\}$. Then, for any fixed $\mathbf{x} \in X$, the follower's problem 
	\begin{equation} \label{eq: follower's problem}
		\max_{\mathbf{y} \in Y(\mathbf{x})} \, \min_{\mathbb{Q}_f \in \mathcal{Q}_f} \rho_f(\mathbf{y}, \mathbb{Q}_f)
	\end{equation}
	can be equivalently reformulated as a linear programming problem of the form: \upshape
	\begin{subequations} \label{follower's LP}
		\begin{align} 
			&\max_{\mathbf{y},\, \boldsymbol{\nu}_f, \,\mathbf{s}_f,\, \lambda_f,\, t_f} \Big\{t_f - \frac{1}{1 - \alpha_f}\Big(\varepsilon_f \lambda_f + \frac{1}{k_f} \; \sum_{k \in K_f} s_f^{\text{\tiny (\textit{k})}}\Big)\Big\} \label{obj: follower's LP} \\
			\mbox{\upshape s.t. } 	& \mathbf{F}\mathbf{y} \leq \mathbf{f} - \mathbf{L}\mathbf{x} \label{cons: follower's LP 1}\\
			& \begin{rcases} 
				 - \hat{\mathbf{c}}_f^{\text{\tiny (\textit{k})}\top} \mathbf{y} 	 + \; t_f + \boldsymbol{\Delta}_f^{\text{\tiny (\textit{k})}\top} \boldsymbol{\nu}_f^{\text{\tiny (\textit{k})}} \leq s^{\text{\tiny (\textit{k})}}_{f} \quad \\
				\Vert \mathbf{B}^\top \boldsymbol{\nu}_f^{\text{\tiny (\textit{k})}} + \mathbf{y} \Vert_{q} \leq \lambda_f \\
				\boldsymbol{\nu}_f^{\text{\tiny (\textit{k})}} \geq \mathbf{0}, \;
				s_f^{\text{\tiny (\textit{k})}} \geq 0 
			\end{rcases} \forall k \in K_f \label{cons: follower's LP 2-5} \\
			& \mathbf{y} \geq \mathbf{0}, \; \lambda_f \geq 0. \label{cons: follower's LP 6} \end{align}
	\end{subequations}
\itshape where \upshape $\boldsymbol{\Delta}_f^{\text{\tiny (\textit{k})}} = \mathbf{b} - \mathbf{B}
\, \hat{\mathbf{c}}_f^{\text{\tiny (\textit{k})}} \geq \mathbf{0}$ \itshape for each $k \in K_f$ and $\frac{1}{p} + \frac{1}{q} = 1$. 
	Furthermore, \upshape (\ref{follower's LP}) \itshape admits the following dual reformulation: \upshape
	\begin{subequations} \label{follower's dual LP}
		\begin{align} 
			& \min_{\overline{\boldsymbol{\mu}}_{f},\, \underline{\boldsymbol{\mu}}_{f}, \,\boldsymbol{\beta}_f,\, \boldsymbol{\gamma}_f} \Big\{( - \mathbf{L}\mathbf{x} + \mathbf{f})^\top \boldsymbol{\beta}_f\Big\} \\
			\mbox{\upshape s.t. } 
			& \mathbf{F}^\top \boldsymbol{\beta}_f - \sum_{k \in K_f} \gamma^{\text{\tiny (\textit{k})}}_{f} \hat{\mathbf{c}}_f^{\text{\tiny (\textit{k})}} + \sum_{k \in K_f} \big(\overline{\boldsymbol{\mu}}_{f}^{\text{\tiny (\textit{k})}} - \underline{\boldsymbol{\mu}}_{f}^{\text{\tiny (\textit{k})}} \big) \geq \mathbf{0} \label{cons: follower's LP dual 1} \\
			& \begin{rcases} \gamma^{\text{\tiny (\textit{k})}}_{f} \boldsymbol{\Delta}_f^{\text{\tiny (\textit{k})}} + \mathbf{B} (\overline{\boldsymbol{\mu}}_{f}^{\text{\tiny (\textit{k})}} - \underline{\boldsymbol{\mu}}_{f}^{\text{\tiny (\textit{k})}}) \geq \mathbf{0} \quad \\
				0 \leq \gamma^{\text{\tiny (\textit{k})}}_{f} \leq \frac{1}{(1 - \alpha_f)k_f} \\
				\overline{\boldsymbol{\mu}}_{f}^{\text{\tiny (\textit{k})}} \geq \mathbf{0}, \;
				\underline{\boldsymbol{\mu}}_{f}^{\text{\tiny (\textit{k})}} \geq \mathbf{0} \quad \end{rcases} \forall k \in K_f \label{cons: follower's LP dual 2-4}\\
			& \Big\Vert \sum_{k \in K_f}\big(\overline{\boldsymbol{\mu}}_{f}^{\text{\tiny (\textit{k})}} + \underline{\boldsymbol{\mu}}_{f}^{\text{\tiny (\textit{k})}}\big) \Big\Vert_p \leq \frac{1}{1 - \alpha_f}\varepsilon_f \label{cons: follower's LP dual 5}\\
			& \sum_{k \in K_f} \gamma^{\text{\tiny (\textit{k})}}_{f} = 1 \label{cons: follower's LP dual 6} \\
			& \boldsymbol{\beta}_f \geq \mathbf{0} \label{cons: follower's LP dual 7}.
		\end{align}
	\end{subequations}
	\begin{proof} See Supplementary Material A.
	\end{proof}	
\end{lemma}

\begin{theorem} \label{theorem 1} 
	If the conditions of Lemma \ref{lemma 1} are satisfied, then the basic DRI problem \upshape [\textbf{DRI}] \itshape admits the following single-level reformulation: \upshape
	\begin{subequations} \label{MILP reformulation}
		\begin{align} 
			\hat{z}_b^* = & \min_{\mathbf{x}, \mathbf{y}, \boldsymbol{\nu}, \boldsymbol{s}, \boldsymbol{\mu} ,\boldsymbol{\beta}, \boldsymbol{\gamma}, \lambda, t} \Big\{ t_l + \frac{1}{1 - \alpha_l} \big(\varepsilon_l \lambda_l + \frac{1}{k_l} \sum_{k \in K_l} s^{\text{\tiny (\textit{k})}}_{l} \big)\Big\} \label{obj: MILP}\\
			\mbox{\upshape s.t. } & \mathbf{x} \in X \label{cons: MILP 1} \\
			& \begin{rcases} 
				 \hat{\mathbf{c}}_l^{\text{\tiny (\textit{k})}\top} \mathbf{y} - t_l + \boldsymbol{\Delta}_l^{\text{\tiny (\textit{k})}\top} \boldsymbol{\nu}_l^{\text{\tiny (\textit{k})}} \leq s^{\text{\tiny (\textit{k})}}_{l} \quad\\	
				\Vert \mathbf{B}^\top \boldsymbol{\nu}_l^{\text{\tiny (\textit{k})}} - \mathbf{y} \Vert_{q} \leq \lambda_l \\
				\boldsymbol{\nu}_l^{\text{\tiny (\textit{k})}} \geq \mathbf{0}, \; s^{\text{\tiny (\textit{k})}}_{l} \geq 0
			\end{rcases} \forall{k} \in K_l \label{cons: MILP 2-4}\\
			&	\mbox{\upshape (\ref{cons: follower's LP 1})--(\ref{cons: follower's LP 6}), (\ref{cons: follower's LP dual 1})--(\ref{cons: follower's LP dual 7})} \\
			& ( - \mathbf{L}\mathbf{x} + \mathbf{f} )^\top \boldsymbol{\beta}_f = t_f - \frac{1}{1 - \alpha_f}\big(\varepsilon_f \lambda_f + \frac{1}{k_f} \; \sum_{k \in K_f} s^{\text{\tiny (\textit{k})}}_{f}\big), \label{cons: MILP 5}
		\end{align}
	\end{subequations}
 \itshape where \upshape $\boldsymbol{\Delta}_l^{\text{\tiny (\textit{k})}} = \mathbf{b} - \mathbf{B} \hat{\mathbf{c}}_l^{\text{\tiny (\textit{k})}} \geq \mathbf{0}$ \itshape for each $k \in K_l$ and $\frac{1}{p} + \frac{1}{q} = 1$. 
	\begin{proof} See Supplementary Material A.
	\end{proof}
\end{theorem}

 Notably, the single-level reformulation (\ref{MILP reformulation}) contains a product of continuous dual variables,~$\boldsymbol{\beta}_f$, and binary decision variables, $\mathbf{x}$. This product can be linearized using standard Big-$M$ constraints, yielding a single-level MILP reformulation of [\textbf{DRI}]; see, e.g., \cite{Zare2019}. As the choice of Big-$M$ depends on the underlying problem's structure, we provide a detailed linearization approach along with specific test instances in Section \ref{sec: comp study}. Finally, the existence of a polynomial-size MILP reformulation for [\textbf{DRI}] rules out $\Sigma^p_2$-hardness and implies that [\textbf{DRI}] is $NP$-hard, matching the complexity of the deterministic problem [\textbf{DI}]. Although not stated explicitly, analogous reformulation-based arguments suggest that the related interdiction models in \cite{Lei2018, Sadana2023} exhibit comparable computational complexity. 

\subsection{Convergence analysis}
Next, we analyze a relation between [\textbf{DRI}] and the respective full information model [\textbf{SI}], as the sample sizes of the leader and the follower tend to infinity. Our results are consistent with the standard asymptotic findings for one-stage Wasserstein DRO problems \cite{Esfahani2018}, reinforcing that [\textbf{DRI}] serves as a valid approximation of [\textbf{SI}]. To establish our convergence results, we introduce the following additional~assumption:
\begin{itemize}
	\item[\textbf{A4}.] For given numbers of samples, $k_l$ and $k_f$, and significance levels, $\theta_{k_l} \in (0, 1)$ and $\theta_{k_f} \in (0, 1)$, the Wasserstein radii, $\varepsilon_l(\theta_{k_l})$ and $\varepsilon_f(\theta_{k_f})$, can be selected so that the true distribution \(\mathbb{Q}^*\) belongs to~$\mathcal{Q}_l$ and $\mathcal{Q}_f$, respectively, with probabilities of at least $1 - \theta_{k_l}$ and $1 - \theta_{k_f}$. Moreover, let
	\[\sum_{k_i = 1}^{\infty} \theta_{k_i} < \infty \, \, \mbox{ and } \, \, \lim_{k_i \to \infty} \varepsilon_i(\theta_{k_i}) = 0 \quad \forall i \in \{l, f\}.\] 
\end{itemize}

Assumption \textbf{A4} is supported by standard measure concentration results for the Wasserstein distance and is essential for the convergence of one-stage Wasserstein DRO models; see Theorem 3.6 in~\cite{Esfahani2018}.
In our setting, due to a similar structure of the leader's and the follower's ambiguity sets, the conditions of Assumption \textbf{A4} are simply reiterated for both decision-makers. Also, in the remainder of this section, $\mathbb{Q}^{* \infty}$ denotes the limiting distribution of the sample's behavior, as the sample size approaches infinity; recall our notations at the end of Section \ref{subsec: approach and contributions}.

\subsubsection{Convergence for the follower}
We begin with analyzing the follower's problem in [\textbf{DRI}] for a fixed leader's decision $\mathbf{x} \in X$. Thus, by leveraging the minimax theorem (see, e.g., Example 3 in \cite{Delage2010}), the follower's problem in~(\ref{cons: bilevel DR 2}) can be expressed as:
\begin{equation} \label{eq: follower's problem convergence}
	\hat{z}^*_f(\mathbf{x}, \hat{\mathbf{C}}_f) := \max_{\bar{\mathbf{y}} \in \bar{Y}(\mathbf{x})\,} \, \min_{\mathbb{Q}_f \in \mathcal{Q}_f} \mathbb{E}_{\mathbb{Q}_f} \big\{ h_f( \bar{\mathbf{y}}, \mathbf{c})\big\},
\end{equation}
where $\bar{\mathbf{y}} := (\mathbf{y}, t_f)$, \[h_f( \bar{\mathbf{y}}, \mathbf{c}) := t_f + \frac{1}{1 - \alpha_f} (\mathbf{c}^\top \mathbf{y} - t_f)^- \; \, \mbox{ and } \, \; \bar{Y}(\mathbf{x}) := Y(\mathbf{x}) \times [\,\underline{t}_f, \overline{t}_f]\] for some appropriate $\underline{t}_f$ and $\overline{t}_f$. In particular, compactness of the sets $X$, $Y(\mathbf{x})$ and $S$, which is guaranteed by Assumptions~\textbf{A1} and \textbf{A2}, yields that any optimal $t^*_f$ in (\ref{eq: left-tail CVaR basic}) is such that: 
\[
-\infty < \underline{t}_f := \min_{\mathbf{x} \in X} \; \min_{(\mathbf{c}, \, \mathbf{y}) \, \in \, S \times Y(\mathbf{x})} \mathbf{c}^\top \mathbf{y} \leq t^*_f \leq \overline{t}_f := \max_{\mathbf{x} \in X} \; \max_{(\mathbf{c}, \,\mathbf{y}) \, \in \, S \times Y(\mathbf{x})} \mathbf{c}^\top \mathbf{y} < +\infty.
\]
The following result establishes convergence for the follower.

\begin{lemma} \label{lemma 2}
Let $\mathbf{x} \in X$ be fixed and Assumptions \textbf{A1}-\textbf{A4} hold. Then, $\mathbb{Q}^{* \infty}$-almost surely the follower's optimal objective function value in \upshape (\ref{eq: follower's problem convergence}) \itshape satisfies
	\begin{equation} \label{eq: follower's problem true}
	\liminf_{k_f \to \infty} \,	\hat{z}^*_f(\mathbf{x}, \hat{\mathbf{C}}_f) = z_f^*(\mathbf{x}) := \max_{\bar{\mathbf{y}} \in \bar{Y}(\mathbf{x})\,} \mathbb{E}_{\mathbb{Q}^*} \big\{ h_f( \bar{\mathbf{y}}, \mathbf{c})\big\}.
	\end{equation}
Furthermore, if \upshape $\bar{\mathbf{y}}^*(k_f)$ \itshape is an optimal solution of the follower's problem in \upshape (\ref{eq: follower's problem convergence}), \itshape then $\mathbb{Q}^{* \infty}$-almost surely any accumulation point of \upshape $\big\{\bar{\mathbf{y}}^*(k_f)\big\}_{k_f \in \mathbb{N}}$ \itshape is an optimal solution of the full information problem~in~\upshape(\ref{eq: follower's problem true}). \itshape	
\begin{proof}
Since $\bar{Y}(\mathbf{x})$ and $S$ are compact, we conclude that $h_f( \bar{\mathbf{y}}, \mathbf{c})$ is continuous and bounded uniformly for all feasible $\bar{\mathbf{y}}$ and $\mathbf{c}$. As $Y(\mathbf{x})$ is also closed, the result follows from Theorem 3.6~in~\cite{Esfahani2018}. 
\end{proof}	
\end{lemma}

\subsubsection{Convergence for the leader}
Analogously to the follower's problem, we formulate the leader's problem in [\textbf{DRI}] as:
\begin{equation} \label{eq: leader's problem convergence}
	\hat{z}_l^*(\hat{\mathbf{C}}_l, \hat{\mathbf{C}}_f) := \min_{\bar{\mathbf{x}} \in 
	\bar{X}(\hat{\mathbf{C}}_f)} \; \max_{\mathbb{Q}_l \in \mathcal{Q}_l} \mathbb{E}_{\mathbb{Q}_l} \big\{ h_l( \bar{\mathbf{x}}, \mathbf{c})\big\},
\end{equation}
where $\bar{\mathbf{x}} := (\mathbf{x}, \bar{\mathbf{y}}, t_l)$, $\bar{\mathbf{y}} := (\mathbf{y}, t_f)$, $h_l( \bar{\mathbf{x}}, \mathbf{c}) := t_l + \frac{1}{1 - \alpha_l} (\mathbf{c}^\top \mathbf{y} - t_l)^+$ and 
\begin{equation} \label{eq: leader's problem convergence feasible set}
\bar{X}(\hat{\mathbf{C}}_f) := \Big\{(\mathbf{x}, \bar{\mathbf{y}}, t_l): \; \mathbf{x} \in X, \; \bar{\mathbf{y}} \in 	\hat{V}_f(\mathbf{x}, \hat{\mathbf{C}}_f), \; t_l \in [\underline{t}_l, \overline{t}_l] \Big\}.
\end{equation}
Specifically, we define $\underline{t}_l := \underline{t}_f$, $\overline{t}_l := \overline{t}_f$ and
\begin{equation} \label{eq: follower's optimal set convergence} 
	\hat{V}_f(\mathbf{x}, \hat{\mathbf{C}}_f) := \argmax_{\, \bar{\mathbf{y}} \, \in \, \bar{Y}(\mathbf{x})} \Big\{ \min_{\mathbb{Q}_f \in \mathcal{Q}_f} \mathbb{E}_{\mathbb{Q}_f} \{ h_f( \bar{\mathbf{y}}, \mathbf{c}) \}\Big\}.
\end{equation}

To establish convergence for the leader, we need to make the following additional assumption: 
\begin{itemize}
\item[\textbf{A5}.] For every fixed $\mathbf{x} \in X$, the full information follower's problem 
\begin{equation} \label{eq: full information follower's problem}
\max_{\,\bar{\mathbf{y}} \in \bar{Y}(\mathbf{x})} \mathbb{E}_{\mathbb{Q}^*} \{ h_f( \bar{\mathbf{y}}, \mathbf{c}) \},
\end{equation}
has a unique optimal solution. 
\end{itemize}
Assumption \textbf{A5} is often used for deriving theoretical properties of bilevel optimization problems; see, e.g., \cite{Dempe2002}. Some additional limitations related to Assumption \textbf{A5} are discussed later, after presenting the proof of our convergence~result.

\begin{lemma} \label{lemma 3} Let 
	\begin{equation} \nonumber
	v_f( \bar{\mathbf{y}}, \hat{\mathbf{C}}_f) := \min_{\mathbb{Q}_f \in \mathcal{Q}_f(\hat{\mathbf{C}}_f)} \mathbb{E}_{\mathbb{Q}_f} \{ h_f(\bar{\mathbf{y}}, \mathbf{c})\} 
	\end{equation} 
	be the follower's optimal objective function value 
	in \upshape (\ref{eq: follower's optimal set convergence}). \itshape Then, $v_f( \bar{\mathbf{y}}, \hat{\mathbf{C}}_f)$ is continuous in both $ \bar{\mathbf{y}}$ and $\hat{\mathbf{C}}_f$, and the follower's optimal solution set $\hat{V}_f(\mathbf{x}, \hat{\mathbf{C}}_f)$ given by equation \upshape (\ref{eq: follower's optimal set convergence}) \itshape is non-empty and upper semicontinuous in $\hat{\mathbf{C}}_f$ for any fixed $\mathbf{x} \in X$. 
\begin{proof} See Supplementary Material A.
\end{proof}	
\end{lemma}

\begin{theorem} \label{theorem 2}
Let Assumptions \textbf{A1}-\textbf{A5} hold. Then, $\mathbb{Q}^{* \infty}$-almost surely the leader's optimal objective function value in \upshape (\ref{eq: leader's problem convergence}) \itshape satisfies
\begin{equation} \label{eq: leader's problem true}
		\limsup_{k_f \to \infty} \, \limsup_{k_l \to \infty} \,	\hat{z}_l^*(\hat{\mathbf{C}}_l, \hat{\mathbf{C}}_f) = z^* := \min_{\bar{\mathbf{x}} \in 
			\bar{X}^{*}} \mathbb{E}_{\mathbb{Q}^*} \big\{ h_l( \bar{\mathbf{x}}, \mathbf{c})\big\},
\end{equation}
where $z^*$ is the optimal objective function value of the stochastic programming problem \upshape [\textbf{SP}]\itshape, with 
\begin{align} 
& \bar{X}^* := \Big\{(\mathbf{x}, \bar{\mathbf{y}}, t_l): \; \mathbf{x} \in X, \; \bar{\mathbf{y}} \in V^*_f(\mathbf{x}), \; t_l \in [\underline{t}_l, \overline{t}_l] \Big\} \; \mbox{ and } \label{eq: leader's feasible set convergence} \\
&	V^*_f(\mathbf{x}) := \argmax_{\, \bar{\mathbf{y}} \, \in \, \bar{Y}(\mathbf{x})} \mathbb{E}_{\mathbb{Q}^*} \{ h_f( \bar{\mathbf{y}} , \mathbf{c}) \}. \label{eq: follower's optimal set true convergence}
\end{align}
Furthermore, if \upshape $\bar{\mathbf{x}}^*(k_l, k_f)$ \itshape 
is an optimal solution of the leader's problem in \upshape (\ref{eq: leader's problem convergence})\itshape, then~$\mathbb{Q}^{* \infty}$-almost surely any resulting accumulation point of \upshape $\big\{\bar{\mathbf{x}}^*(k_l, k_f) \big\}_{k_l,k_f \in \mathbb{N}}$ \itshape is an optimal solution of the full information problem in the right hand side of \upshape (\ref{eq: leader's problem true}). \itshape	
\begin{proof}
Let the follower's sample size, $k_f$, and its data set, $\hat{\mathbf{C}}_f$, be fixed. Additionally, based on Assumption~\textbf{A3}, we also fix the first $k_f$ samples in the leader's data set, $\hat{\mathbf{C}}_l$. It is important to note, however, that fixing a finite number of samples in the leader's data set does not impact our convergence results, since we focus only on asymptotic performance guarantees for the leader.

\textit{Closedness of $\bar{X}(\hat{\mathbf{C}}_f)$.}
Initially, we demonstrate that Theorem 3.6 in \cite{Esfahani2018} can also be applied to the leader's problem in~(\ref{eq: leader's problem convergence}). Similar to the proof of Lemma \ref{lemma 2}, we observe that $h_l( \bar{\mathbf{x}} , \mathbf{c})$ is continuous and bounded uniformly for all feasible~ $\bar{\mathbf{x}}$ and $\mathbf{c}$. Thus, it suffices to show that the leader's feasible set~$\bar{X}(\hat{\mathbf{C}}_f)$ given by equation (\ref{eq: leader's problem convergence feasible set}) is closed. Let $\bar{\mathbf{x}}(k) = (\mathbf{x}, \mathbf{y}, t_f, t_l) \in \bar{X}(\hat{\mathbf{C}}_f)$
be a sequence of points such that
\[\lim_{k \to \infty} \bar{\mathbf{x}}(k) = \bar{\mathbf{x}}^{\text{\tiny lim}} = (\mathbf{x}^{\text{\tiny lim}}, \mathbf{y}^{\text{\tiny lim}}, t_f^{\text{\tiny lim}}, t_l^{\text{\tiny lim}});\]
for simplicity, the dependence on $k$ is omitted in the notation whenever it is clear from the context. Since $X$ is a discrete set, it can be observed that $\mathbf{x}^{\text{\tiny lim}} \in X$ and furthermore there exists $k_0 \in \mathbb{N}$ such that, for all $k \geq k_0$, we have $\mathbf{x}(k) = \mathbf{x}^{\text{\tiny lim}}$. Also, $t_l^{\text{\tiny lim}} \in [\underline{t}_l, \overline{t}_l]$ and, to verify that $\bar{X}(\hat{\mathbf{C}}_f)$ is closed, we only need to check that 
\begin{equation} \label{eq: closedeness convergence}
\bar{\mathbf{y}}^{\text{\tiny lim}} = (\mathbf{y}^{\text{\tiny lim}}, t_f^{\text{\tiny lim}}) \in \hat{V}_f(\mathbf{x}^{\text{\tiny lim}}, \hat{\mathbf{C}}_f),
\end{equation}
given that \[\bar{\mathbf{y}}(k) =\big(\mathbf{y}, t_f \big) \in \hat{V}_f(\mathbf{x}^{\text{\tiny lim}}, \hat{\mathbf{C}}_f) \quad \forall \, k \geq k_0.\] 

The latter observation follows from the closedness of the optimal solution set $\hat{V}_f(\mathbf{x}^{\text{\tiny lim}}, \hat{\mathbf{C}}_f)$. Indeed, by Lemma \ref{lemma 3}, the follower's objective function in (\ref{eq: follower's optimal set convergence}) is continuous and $\hat{V}_f(\mathbf{x}^{\text{\tiny lim}}, \hat{\mathbf{C}}_f)$ is non-empty. Since $\bar{Y}(\mathbf{x}^{\text{\tiny lim}})$ is also compact, $\bar{\mathbf{y}}(k)$, $k \geq k_0$, and $\bar{\mathbf{y}}^{\text{\tiny lim}}$ provide the same optimal objective function value for the follower and the relation (\ref{eq: closedeness convergence}) holds. 
	
As a result, Theorem 3.6 in \cite{Esfahani2018} implies that $\mathbb{Q}^{*\infty}$-almost surely
\begin{equation} \label{eq: leader's problem convergence 3} 
\limsup_{k_l \to \infty} \hat{z}_{ l }^*(\hat{\mathbf{C}}_l, \hat{\mathbf{C}}_f) = \hat{z}_{ l }^*(\hat{\mathbf{C}}_f) = \min_{ \bar{\mathbf{x}} \in \bar{X}(\hat{\mathbf{C}}_f)} \mathbb{E}_{\mathbb{Q}^*} \left\{ h_l( \bar{\mathbf{x}} , \mathbf{c}) \right\}.
\end{equation}
Also, given an optimal solution 
$\bar{\mathbf{x}}^*(k_l, k_f)$ 
of the leader's problem (\ref{eq: leader's problem convergence}), any accumulation point 
\begin{equation} \label{eq: partial full optimal solution} 
\bar{\mathbf{x}}^*(k_f) := (\tilde{\mathbf{x}}^*, \tilde{\mathbf{y}}^*, \tilde{t}^*_f, \tilde{t}_l^*)
\end{equation} of the sequence $\big\{\bar{\mathbf{x}}^*(k_l, k_f) \big\}_{k_l \in \mathbb{N}}$ is an optimal solution of (\ref{eq: leader's problem convergence 3}).

\textit{Feasibility of $\bar{\mathbf{x}}^*$.}
Next, we consider an accumulation point 
\begin{equation} \label{eq: full optimal solution} 
\bar{\mathbf{x}}^* = (\mathbf{x}^*, \mathbf{y}^*, t^*_f, t_l^*)
\end{equation} of the sequence~$\big\{\bar{\mathbf{x}}^*(k_f) \big\}_{k_f \in \mathbb{N}}$, as $k_f$ tends to infinity; this point is essentially an accumulation point of the double-indexed sequence~$\big\{\bar{\mathbf{x}}^*(k_l, k_f)\big\}_{k_l,k_f \in \mathbb{N}}$ by~construction. Our goal is to show that
\begin{equation} \label{eq: leader's problem true 2}
\limsup_{k_f \to \infty} \hat{z}_{ l }^*(\hat{\mathbf{C}}_f) = z^* = \min_{ \bar{\mathbf{x}} \in \bar{X}^*} \mathbb{E}_{\mathbb{Q}^*} \left\{ h_l( \bar{\mathbf{x}} , \mathbf{c}) \right\}, 
\end{equation}
where $\hat{z}_{ l }^*(\hat{\mathbf{C}}_f)$ is defined by equation (\ref{eq: leader's problem convergence 3}), and that $\bar{\mathbf{x}}^*$ is an optimal solution of the full information problem in the right-hand side of (\ref{eq: leader's problem true 2}). 

We first prove that $\bar{\mathbf{x}}^*$ is feasible in (\ref{eq: leader's problem true 2}), i.e., $\bar{\mathbf{x}}^* \in \bar{X}^*$, where $\bar{X}^*$ is defined by equation~(\ref{eq: leader's feasible set convergence}). By passing to a subsequence, if necessary, let 
\begin{equation} \label{eq: accumulation point follower} \nonumber
	\lim_{k_f \to \infty} \bar{\mathbf{x}}^*(k_f) = \bar{\mathbf{x}}^*,
\end{equation} 
where we recall that $\bar{\mathbf{x}}^*(k_f) \in \bar{X}(\hat{\mathbf{C}}_f)$ is an optimal solution of (\ref{eq: leader's problem convergence 3}). 
By similar reasoning as for the closedness of $\bar{X}(\hat{\mathbf{C}}_f)$, we conclude that $t_l^* \in [\underline{t}_l, \overline{t}_l]$ and there exists $k_0 \in \mathbb{N}$ such that, for all $ k_f \geq k_0$, we have $\tilde{\mathbf{x}}^*(k_f) = \mathbf{x}^* \in X$; recall definitions (\ref{eq: partial full optimal solution}) and (\ref{eq: full optimal solution}). 
Assuming that $k_f \geq k_0$ and, therefore, $\tilde{\mathbf{x}}^*(k_f) = \mathbf{x}^*$, Lemma \ref{lemma 2} implies that 
\[\bar{\mathbf{y}}^* := (\mathbf{y}^*, t^*_f) = 	\lim_{k_f \to \infty} \big(\tilde{\mathbf{y}}^*(k_f), \tilde{t}^*_f(k_f) \big)\] is an optimal solution of the full information follower's problem 
\begin{equation} \label{eq: follower's problem true 2} \nonumber
\max_{ \bar{\mathbf{y}} \in \bar{Y}(\mathbf{x}^*)\,} \mathbb{E}_{\mathbb{Q}^*} \left\{ h_f( \bar{\mathbf{y}} , \mathbf{c})\right\},
\end{equation}
and, therefore, $\bar{\mathbf{x}}^*$ is feasible in (\ref{eq: leader's problem true 2}), i.e., $\bar{\mathbf{x}}^* \in \bar{X}^*$.

\textit{Upper bound.}
In order to prove the equality in (\ref{eq: leader's problem true 2}), we first prove that
\[\limsup_{k_f \to \infty} \hat{z}_{ l }^*(\hat{\mathbf{C}}_f) \geq z^*.\]
Indeed, given a sequence $\big\{\bar{\mathbf{x}}^*(k_f)\big\}_{k_f \in \mathbb{N}}$ and its accumulation point~$\bar{\mathbf{x}}^* \in \bar{X}^*$ defined by equations~(\ref{eq: partial full optimal solution}) and (\ref{eq: full optimal solution}), respectively, we observe that 
\begin{equation} \label{eq: upper bound convergence}
\begin{gathered}
\lim_{k_f \to \infty} \hat{z}_{ l }^*(\hat{\mathbf{C}}_f) = \lim_{k_f \to \infty}	\mathbb{E}_{\mathbb{Q}^*}\left\{h_l\left( \bar{\mathbf{x}}^*(k_f), \mathbf{c}\right) \right\} = \\ \mathbb{E}_{\mathbb{Q}^*}\left\{ \lim_{k_f \to \infty} h_l\left( \bar{\mathbf{x}}^*(k_f), \mathbf{c}\right) \right\} = 
\mathbb{E}_{\mathbb{Q}^*}\left\{ h_l\left( \bar{\mathbf{x}}^*, \mathbf{c}\right) \right\} \geq z^*,
\end{gathered}
\end{equation}
Here, the first equality follows from the fact that $\bar{\mathbf{x}}^*(k_f)$ is an optimal solution of~(\ref{eq: leader's problem convergence 3}). The second equality is due to the bounded convergence theorem \cite{Bartle2014}, which holds since $h_l( \bar{\mathbf{x}}, \mathbf{c})$ is continuous and uniformly bounded. Furthermore, the third equality and the last inequality hold since $h_l( \bar{\mathbf{x}}, \mathbf{c})$ is continuous and $\bar{\mathbf{x}}^* \in \bar{X}^*$, respectively. 

\textit{Lower bound.}
As a result, it remains to show that 
\begin{equation} \nonumber
\limsup_{k_f \to \infty} \hat{z}_{ l }^*(\hat{\mathbf{C}}_f) \leq z^*.
\end{equation}
By definition of the minimum in~(\ref{eq: leader's problem convergence 3}), we have
\begin{equation} \label{eq: optimal partial solution convergence}
\hat{z}_{ l }^*(\hat{\mathbf{C}}_f) \leq \mathbb{E}_{\mathbb{Q}^*}\left\{ h_l( \bar{\mathbf{x}}, \mathbf{c}) \right\} \quad \forall \, \bar{\mathbf{x}} \in \bar{X}(\hat{\mathbf{C}}_f).
\end{equation}
Let $\bar{\mathbf{x}}^{\text{\tiny 0}} = (\mathbf{x}^{\text{\tiny 0}}, \mathbf{y}^{\text{\tiny 0}}, t^{\text{\tiny 0}}_f, t^{\text{\tiny 0}}_l)$ be an optimal solution of the full information problem in the right hand side of~(\ref{eq: leader's problem true}). The main challenge is that there is no guarantee that $\bar{\mathbf{x}}^{\text{\tiny 0}} \in \bar{X}(\hat{\mathbf{C}}_f)$. 
For this reason, we define $\bar{\mathbf{x}} = (\mathbf{x}, \bar{\mathbf{y}}, t_l)$ in the right-hand side of (\ref{eq: optimal partial solution convergence}) with $\mathbf{x} := \mathbf{x}^{\text{\tiny 0}}$, $t_l := t^{\text{\tiny 0}}_l$ and an \textit{arbitrary}~point
\begin{equation} \label{eq: partial full optimal solution 2}
 \bar{\mathbf{y}} := \bar{\mathbf{y}}^{\text{\tiny 0}}(k_f) = (\tilde{\mathbf{y}}^{\text{\tiny 0}}, \tilde{t}^{\,\mbox{\tiny 0}}_f) \in \hat{V}_f(\mathbf{x}^{\text{\tiny 0}}, \hat{\mathbf{C}}_f),
\end{equation} where $\hat{V}_f(\mathbf{x}^{\text{\tiny 0}}, \hat{\mathbf{C}}_f)$ is non-empty by Lemma \ref{lemma 3}. 

 Next, by design, we have 
\begin{equation} \nonumber
\bar{\mathbf{x}}^{\text{\tiny 0}}(k_f) := (\mathbf{x}^{\text{\tiny 0}}, \tilde{\mathbf{y}}^{\text{\tiny 0}}, \tilde{t}^{\,\mbox{\tiny 0}}_f, t^{\text{\tiny 0}}_l) \in \bar{X}(\hat{\mathbf{C}}_f)
\end{equation} and, therefore, (\ref{eq: optimal partial solution convergence}) yields that
\begin{equation} \label{eq: lower bound convergence}
\hat{z}_{ l }^*(\hat{\mathbf{C}}_f) \leq \mathbb{E}_{\mathbb{Q}^*}\left\{ h_l\left( \bar{\mathbf{x}}^{\text{\tiny 0}}(k_f),  \mathbf{c}\right) \right\}.
\end{equation}
Also, by Lemma \ref{lemma 2} and Assumption \textbf{A5}, any accumulation point of the sequence $\big\{\bar{\mathbf{y}}^{\text{\tiny 0}}(k_f)\big\}_{k_f \in \mathbb{N}}$ defined by equation (\ref{eq: partial full optimal solution 2}) is given~by~$(\mathbf{y}^{\text{\tiny 0}}, t_f^{\text{\tiny 0}})$. Given that the latter sequence is bounded, we conclude that 
\[\lim_{k_f \to \infty} \; \bar{\mathbf{y}}^{\text{\tiny 0}}(k_f) = (\mathbf{y}^{\text{\tiny 0}}, t_f^{\text{\tiny 0}}),\]
where the limit is understood with respect to all convergent subsequences. 

By using the same argument as in (\ref{eq: upper bound convergence}), we observe that
\[\lim_{k_f \to \infty} \; \mathbb{E}_{\mathbb{Q}^*}\left\{ h_l\left( \bar{\mathbf{x}}^{\text{\tiny 0}}(k_f), \mathbf{c}\right) \right\} = \mathbb{E}_{\mathbb{Q}^*}\left\{\lim_{k_f \to \infty} h_l\left( \bar{\mathbf{x}}^{\text{\tiny 0}}(k_f), \mathbf{c}\right) \right\} = \mathbb{E}_{\mathbb{Q}^*}\left\{h_l( \bar{\mathbf{x}}^{\text{\tiny 0}}, \mathbf{c}) \right\} = z^*.\]
By applying the limit superior to both sides of (\ref{eq: lower bound convergence}), we have: 
\begin{equation} \label{eq: lower bound convergence 2} 
	\limsup_{k_f \to \infty} \hat{z}_{ l }^*(\hat{\mathbf{C}}_f) \leq \limsup_{k_f \to \infty} \mathbb{E}_{\mathbb{Q}^*}\left\{ h_l( \bar{\mathbf{x}}^{\text{\tiny 0}}(k_f), \mathbf{c}) \right\} = z^*.
\end{equation}
Therefore, taking into account (\ref{eq: upper bound convergence}), we conclude that $\limsup_{k_f \to \infty} \hat{z}_{ l }^*(\hat{\mathbf{C}}_f) = z^*$. Finally, equations~(\ref{eq: upper bound convergence}) and (\ref{eq: lower bound convergence 2}) imply that 
\[z^* \geq	\limsup_{k_f \to \infty} \hat{z}_{ l }^*(\hat{\mathbf{C}}_f) \geq \lim_{k_f \to \infty}\hat{z}_{ l }^*(\hat{\mathbf{C}}_f) = \mathbb{E}_{\mathbb{Q}^*}\left\{ h_l(\bar{\mathbf{x}}^*, \mathbf{c}) \right\} \geq z^*,\]
which means that $\bar{\mathbf{x}}^* = (\mathbf{x}^*, \mathbf{y}^*, t^*_f, t^*_l)$ is an optimal solution of (\ref{eq: leader's problem true}). This observation concludes the~proof. 
\end{proof}
\end{theorem}

Notably, in the formulation of Theorem \ref{theorem 2}, our focus is on the limit superior, as the leader minimizes its objective function in (\ref{eq: leader's problem convergence}) with respect to all feasible values of $\mathbf{x}$ and $\mathbf{y}$ and, therefore, convergence in terms of the worst-case upper bound is of main interest. Furthermore, the proof of Theorem~\ref{theorem 2} exploits discreteness of the leader's feasible set $X$ and uniqueness of a solution to the full information follower's problem (\ref{eq: full information follower's problem}); recall Assumption \textbf{A5}. 

 Admittedly, \textbf{A5} cannot be verified in practice, since the true distribution~$\mathbb{Q}^*$ is not known. However, even if Assumption \textbf{A5} is violated, equation~(\ref{eq: upper bound convergence}) shows that the leader's optimal objective function value in [\textbf{DRI}] converges to an upper bound of the true optimal objective function value,~$z^*$. 
This observation suggests that, even though we consider an optimistic version of [\textbf{DRI}], it may offer a more conservative approach from the leader's perspective compared to [\textbf{SI}].

\section{Pessimistic approximation} \label{sec: pessimistic}
\subsection{Problem formulation and complexity}
In line with Example \ref{ex: example 1}, to relax Assumption \textbf{A3}, we consider a \textit{true basic model} given by:
\begin{subequations} \label{bilevel distributionally robust problem true} 
	\begin{align}
		\mbox{[\textbf{DRI}$^*$]: \quad } z^*_b := & \min_{\mathbf{x}, \mathbf{y}} \; \Big\{ \max_{\mathbb{Q}_l \in \mathcal{Q}_l} \;\rho_l(\mathbf{y}, \mathbb{Q}_l) \Big\}\label{obj: bilevel DR true}\\
		\mbox{s.t. } & \mathbf{x} \in X \label{cons: bilevel DR 1 true} \\
		& \mathbf{y} \in \argmax_{\,\tilde{\mathbf{y}} \in Y(\mathbf{x})} \; \Big\{ \min_{\mathbb{Q}_f \in \mathcal{Q}_f(\hat{\mathbf{C}}^*_f)} \rho_f(\tilde{\mathbf{y}}, \mathbb{Q}_f) \Big\}, \label{cons: bilevel DR 2 true}
	\end{align}
\end{subequations}
where $\hat{\mathbf{C}}^*_f$ denotes the true follower's data set. 
In essence, [\textbf{DRI}] and [\textbf{DRI}$^*$] are identical problems. However, in~[\textbf{DRI}$^*$], we do not assume that the follower's data set is a subset of the leader's, i.e., the second part of Assumption~\textbf{A3} does not hold. We emphasize that [\textbf{DRI}$^*$] is a \textit{hypothetical} model, as whenever the leader has access to $\hat{\mathbf{C}}^*_f$, it should also incorporate this data into its own data set.

We first propose a \textit{pessimistic approximation} of the true basic model~[\textbf{DRI}$^*$]; recall case~(\textit{i}) in Example \ref{ex: example 1}.
Formally, the pessimistic approximation can be formulated as: 
\begin{equation} 
	\mbox{[\textbf{DRI-P}]:} \quad \hat{z}^*_{p} := \min_{\mathbf{x} \in X} \; \max_{\mathbf{y} \in Y(\mathbf{x})} \; \max_{\mathbb{Q}_l \in \mathcal{Q}_l} \; \rho_l(\mathbf{y}, \mathbb{Q}_l), \label{pessimistic approximation} 
\end{equation}
where $\rho_l(\mathbf{y}, \mathbb{Q}_l)$ is defined by equation (\ref{eq: right-tail CVaR basic}). That is, the leader in [\textbf{DRI-P}] disregards any available information about the true follower's data set $\hat{\mathbf{C}}^*_f$ and selects the worst-case \textit{feasible} follower's policy in terms of the leader's objective function value. 

 Next, assuming that the support set $S$ is not necessarily compact, we prove that [\textbf{DRI-P}] is $\Sigma_2^p$-hard. Notably, while Assumption~\textbf{A2} requires $S$ to be compact, more general Wasserstein DRO models allow unbounded support sets; see, e.g., \cite{Esfahani2018, Gao2023}. The following result holds.

\begin{theorem} \label{theorem 3} Let Assumption \textbf{A1} hold and assume that the support set $S$ defined by equation~\upshape(\ref{eq: support set}) \itshape is not necessarily compact. Then, the pessimistic approximation \upshape [\textbf{DRI-P}] \itshape is $\Sigma_2^p$-hard. 
\begin{proof}

We show $\Sigma_2^p$-hardness by a reduction from the dominating set interdiction problem \cite{Grune2025}.
Given an undirected graph $G = (V, E)$, a \emph{dominating set} is a subset of vertices $D \subseteq V$ such that every vertex $v \in V$ is either in $D$ or adjacent to at least one vertex in $D$. That is, for every $v \in V$, either $v \in D$ or there exists $u \in D$ such that $(u, v) \in E$. A decision version of the dominating set interdiction problem can be formulated as follows:
 \begin{itemize}
 	\item[$ $] [\textbf{DSI-D}]: Given an undirected graph $G = (V, E)$, a threshold $t \in \mathbb{Z}_+$, and a budget~$b \in \mathbb{Z}_+$, is there a subset of vertices $B \subseteq V$ with $|B| \leq b$ such that there is no dominating set $D \subseteq V \setminus B$ with size $\leq t$?
 \end{itemize}
Importantly, [\textbf{DSI-D}] is proved to be $\Sigma_2^p$-complete; see Theorems 15 and 38 in \cite{Grune2025}. 

Next, we introduce the following three-level problem:
\begin{equation} \label{eq: complexity three-level problem}
	\tilde{z}^* = \min_{\mathbf{x} \in \tilde{X}} \; \max_{\mathbf{y} \in \tilde{Y}(\mathbf{x})} \; \min_{\mathbf{z} \in \tilde{Z}(\mathbf{y})} \; -\sum_{i \in V} z_i,
\end{equation}
where $\tilde{X} = \left\{\mathbf{x} \in \{0, 1\}^{|V|}: \; \sum_{i \in V} x_i \leq b \right\}$,
\begin{subequations}
	\begin{align}
		& \tilde{Y}(\mathbf{x}) = \Big\{\mathbf{y} \in \{0, 1\}^{|V|}: \; y_i \leq 1 - x_i \quad \forall i \in V, \quad \sum_{i \in V} y_i \leq t \Big\} \label{eq: follower's set kernel 1} \; \, \mbox{ and } \\
		& \tilde{Z}(\mathbf{y}) = \Big\{\mathbf{z} \in \mathbb{R}^{|V|}_+: \; z_i \leq 1 - y_i \quad \forall i \in V, \quad z_i \leq 1 - y_j \quad \forall \, (i, j) \in E \Big\}. \label{eq: follower's set kernel 2} 
\end{align}\end{subequations}

We demonstrate that [\textbf{DSI-D}] admits a ``yes''-instance if and only if $\tilde{z}^* < 0$. Indeed, if the answer to [\textbf{DSI-D}] is ``yes'', then we can select $x_i = 1$ for~$i \in B$ and $x_i = 0$ otherwise. Then, since there are no dominating sets of size $\leq t$ in the residual graph, for any $\mathbf{y} \in \tilde{Y}(\mathbf{x})$, we can identify at least one non-dominated vertex $i \in V$ such that $z^*_i = 1$. Hence, $\tilde{z}^* < 0$.

Conversely, assume that $\tilde{z}^* < 0$ but the answer to [\textbf{DSI-D}] is ``no''. Then, for any $B \subseteq V$, $|B| \leq b$, there exists a dominating set of size $\leq t$ in the residual graph. In other words, for any $\mathbf{x} \in \tilde{X}$, there exists $\mathbf{y} \in \tilde{Y}(\mathbf{x})$ such that $\tilde{Z}(\mathbf{y}) = \{\mathbf{0}\}$. This contradicts the assumption that $\tilde{z}^* < 0$.


Finally, we show that (\ref{eq: complexity three-level problem}) reduces to an instance of [\textbf{DRI-P}] with $\alpha_l = 0$ and $\varepsilon_l = \infty$, i.e., \begin{equation} \label{eq: follower's problem risk-neutral robust optimistic}
	\hat{z}^*_p = \min_{\mathbf{x} \in X} \; \max_{\mathbf{y} \in Y(\mathbf{x})} \; \max_{\mathbf{c} \in S} \; \mathbf{c}^\top \mathbf{y}.
\end{equation}
Let $\boldsymbol{\beta} \in \mathbb{R}^{|V|}_+$ and $\boldsymbol{\gamma} \in \mathbb{R}^{|E|}_+$ be dual variables corresponding to the vertex-based and the edge-based constraints in (\ref{eq: follower's set kernel 2}). Then, (\ref{eq: complexity three-level problem}) admits a single-level dual reformulation given by:
\begin{subequations} \label{eq: complexity three-level problem dual}
	\begin{align}
	 \min_{\mathbf{x} \in \tilde{X}} \; & \max_{\mathbf{y}, \boldsymbol{\beta}, \boldsymbol{\gamma}} \;\sum_{j \in V} \left(\beta_j + \sum_{i: \, (i, j) \in E} \gamma_{ij} \right)\left(y_j - 1\right) \\
		\mbox{s.t. } & \mathbf{y} \in \tilde{Y}(\mathbf{x}) \\
		& \beta_{j} + \sum_{i: \, (j, i) \in E} \gamma_{ji} \geq 1 \quad \forall j \in V \label{cons: complexity three-level problem dual 1}\\
		& \boldsymbol{\gamma}, \boldsymbol{\beta} \geq \mathbf{0}. \label{cons: complexity three-level problem dual 2}
\end{align}\end{subequations}
Next, for each $j \in V$, we set $y'_j := 1 - y_j$ and
\[\delta_j := -\beta_j - \sum_{i: \, (i, j) \in E} \gamma_{ij}. 
\]
As a result, (\ref{eq: complexity three-level problem dual}) reads as:
\begin{subequations} \label{eq: kernel dual 2}
	\begin{align}
		 \min_{\mathbf{x} \in \tilde{X}} \; & \max_{\mathbf{y}', \boldsymbol{\beta}, \boldsymbol{\gamma}, \boldsymbol{\delta}} \; \boldsymbol{\delta}^\top \mathbf{y}' \\
		\mbox{s.t. } & \text{(\ref{cons: complexity three-level problem dual 1})--(\ref{cons: complexity three-level problem dual 2})} \label{cons: kernel dual 2 1} \\
		& \delta_j = -\beta_j - \sum_{i: \, (i, j) \in E} \gamma_{ij} \quad \forall j \in V \label{cons: kernel dual 2 2} \\
		& \mathbf{x} \leq \mathbf{y}' \leq \mathbf{1} \label{cons: kernel dual 2 3}\\
		& \sum_{i \in V} y'_i \geq |V| - t \label{cons: kernel dual 2 4}.
\end{align}\end{subequations}
In particular, the integrality constraints for $\mathbf{y}'$ can be relaxed, since the objective function in (\ref{eq: kernel dual 2}) is bilinear and constraints (\ref{cons: kernel dual 2 3})-(\ref{cons: kernel dual 2 4}) define an integral polytope. 
 
We conclude that (\ref{eq: kernel dual 2}) is a special case of (\ref{eq: follower's problem risk-neutral robust optimistic}). Specifically, constraints~(\ref{cons: kernel dual 2 3})-(\ref{cons: kernel dual 2 4}) are associated with the follower's feasible set in (\ref{eq: follower's problem risk-neutral robust optimistic}), which is non-empty and compact by design. On the other hand, after lifting to a higher-dimensional space, constraints (\ref{cons: kernel dual 2 1})-(\ref{cons: kernel dual 2 2}) describe an unbounded support set. This observation concludes the~proof.
\end{proof}
\end{theorem}

Based on Theorem \ref{theorem 3}, we conclude that, even for the risk-neutral leader, the pessimistic approximation~[\textbf{DRI-P}] is $\Sigma_2^p$-hard. We therefore focus on two special cases of the problem, where the leader is either ambiguity-free ($\varepsilon_l = 0$) or risk-neutral ($\alpha_l = 0$). In these cases, the second-level problem in~[\textbf{DRI-P}] reduces to a bilinear problem with disjoint constraints. Furthermore, under mild technical  assumptions, the specific structure of the latter problem enables a \textit{direct} MILP reformulation. Both~versions of [\textbf{DRI-P}] are then solved using a standard Benders-type decomposition algorithm; 
see, e.g.,~\cite{Zeng2014}.
The following result provides a general reformulation of [\textbf{DRI-P}] for a fixed $\mathbf{x} \in X$. 

\begin{theorem} \label{theorem 4}
Let Assumptions \textbf{A1}, \textbf{A2} and the first part of Assumption \textbf{A3} hold. If the leader's ambiguity set $\mathcal{Q}_l$ is defined in terms of $p$-norm with $p \in \{1, \infty\}$, then, for any $\mathbf{x} \in X$, the second-level problem in\upshape~[\textbf{DRI-P}] \itshape admits the following non-convex quadratic reformulation: \upshape
	\begin{subequations} \label{second-level quadratic reformulation pessimistic} 
		\begin{align} 
		& \max_{\mathbf{y}, \, \boldsymbol{\xi}_{l}, \, \boldsymbol{\eta}_{l}, \, \boldsymbol{\gamma}_l} \Big\{\frac{1}{k_l(1 - \alpha_l)} \sum_{k \in K_l} \big( \gamma^{\text{\tiny (\textit{k})}}_{l} \hat{\mathbf{c}}_l^{\text{\tiny (\textit{k})}} - \boldsymbol{\xi}_l^{\text{\tiny (\textit{k})}} \big)^\top \mathbf{y} \Big\} \label{obj: second-level quadratic reformulation pessimistic} \\
			\mbox{\upshape s.t. } 	& \mathbf{y} \in Y(\mathbf{x}) \label{cons: second-level quadratic reformulation pessimistic 1}\\
			& \begin{rcases} 
				\mathbf{B}(\gamma^{\text{\tiny (\textit{k})}}_{l} \hat{\mathbf{c}}_l^{\text{\tiny (\textit{k})}} - \boldsymbol{\xi}_l^{\text{\tiny (\textit{k})}}) \leq \gamma^{\text{\tiny (\textit{k})}}_{l} \mathbf{b} \quad \\
				0 \leq \gamma^{\text{\tiny (\textit{k})}}_{l} \leq 1 \\
				-\boldsymbol{\eta}_l^{\text{\tiny (\textit{k})}} \leq \boldsymbol{\xi}_l^{\text{\tiny (\textit{k})}} \leq \boldsymbol{\eta}_l^{\text{\tiny (\textit{k})}}\\
			\end{rcases} \forall k \in K_l \label{cons: second-level quadratic reformulation pessimistic 2-4} \\
			& \sum_{k \in K_l} \gamma^{\text{\tiny (\textit{k})}}_{l} = k_l (1 - \alpha_l) \label{cons: second-level quadratic reformulation pessimistic 5}\\
			& \Big\Vert \frac{1}{k_l}\sum_{k \in K_l} \boldsymbol{\eta}_l^{\text{\tiny (\textit{k})}} \Big\Vert_p \leq \varepsilon_l. \label{cons: second-level quadratic reformulation pessimistic 6} 
	\end{align}
	\end{subequations} 
	\begin{proof} See Supplementary Material A.
	\end{proof}		
\end{theorem}

 
\subsection{Benders decomposition-based algorithm for the ambiguity-free problem}
First, we provide a reformulation of the pessimistic approximation [\textbf{DRI-P}] for the case of $\varepsilon_l = 0$. 
\begin{corollary} \label{corollary 1}
	Assume that the leader is ambiguity-free, i.e., $\varepsilon_l = 0$. Then, \upshape [\textbf{DRI-P}] \itshape reads as: \upshape
	\begin{subequations} \label{pessimistic approximation SAA} 
		\begin{align} 
			\min_{\mathbf{x} \in X} \;	& \max_{\mathbf{y}, \, \boldsymbol{\gamma}_l}\Big\{\frac{1}{k_l(1 - \alpha_l)}\sum_{k \in K_l} \big( \gamma^{\text{\tiny (\textit{k})}}_{l} \hat{\mathbf{c}}_l^{\text{\tiny (\textit{k})}}\big)^\top \mathbf{y} \Big\} \label{obj: pessimistic approximation SAA} \\
			\mbox{\upshape s.t. } 	& \mathbf{y} \in Y(\mathbf{x}) \label{cons: pessimistic approximation SAA 1}\\
			& 0 \leq \gamma^{\text{\tiny (\textit{k})}}_{l} \leq 1 
			\quad \forall k \in K_l \label{cons: pessimistic approximation SAA 2}\\
			& \sum_{k \in K_l} \gamma^{\text{\tiny (\textit{k})}}_{l} = k_l (1 - \alpha_l). \label{cons: pessimistic approximation SAA 3} \end{align}
	\end{subequations} 
\begin{proof} The result readily follows by setting $\varepsilon_{l} = 0$ in the second-level problem reformulation (\ref{second-level quadratic reformulation pessimistic}).
\end{proof}		
\end{corollary}

Next, we make the following additional assumption:
\begin{itemize}
	\item[\textbf{A6}.] The leader's confidence level $\alpha_l \in (0, 1)$ is such that $k_l(1 - \alpha_l)$ is integer.
\end{itemize}
In fact, Assumption \textbf{A6} does not impose a significant restriction. Indeed, when $k_l(1 - \alpha_l) \notin \mathbb{Z} $, one may slightly perturb $\alpha_l$ or $k_l$ to ensure that the integrality condition holds. 

 Notably, under Assumption~\textbf{A6}, the ambiguity-free problem~(\ref{pessimistic approximation SAA}) can be expressed as:
\begin{subequations} \label{pessimistic approximation SAA 2}
	\begin{align} 
	 & \min_{\mathbf{x},\, z} \; z \\
		\mbox{\upshape s.t. } 
		& \mathbf{x} \in X \\	
		& z \geq \max_{\mathbf{y} \in Y(\mathbf{x})} \Big\{ \bar{\mathbf{c}}(\boldsymbol{\gamma}_l)^\top \mathbf{y} \Big\} \quad \forall \, \boldsymbol{\gamma}_l \in \Gamma, \label{cons: pessimistic approximation SAA max}
	\end{align}
\end{subequations} 
where $\bar{\mathbf{c}}(\boldsymbol{\gamma}_l) := \frac{1}{k_l(1 - \alpha_l)}\sum_{k \in K_l} \gamma_l^{\text{\tiny (\textit{k})}} \hat{\mathbf{c}}_l^{\text{\tiny (\textit{k})}}$ and
\begin{equation} \label{eq: gamma}
	\Gamma := \Big\{\boldsymbol{\gamma}_l \in \{0, 1\}^{k_l}: \sum_{k \in K_l} \gamma_l^{\text{\tiny (\textit{k})}} = k_l (1 - \alpha_l)\Big\}.
\end{equation}
Indeed, the objective function (\ref{obj: pessimistic approximation SAA}) is bilinear and constraints (\ref{cons: pessimistic approximation SAA 2})-(\ref{cons: pessimistic approximation SAA 3}) define an integral polytope, whose extreme points are given by $\Gamma$. 

 
\begin{algorithm}[h]
	\DontPrintSemicolon \onehalfspacing 
	\textbf{Input:} the ambiguity-free problem (\ref{pessimistic approximation SAA}), a feasible leader's decision $\mathbf{x}^{\text{\tiny 0}} \in X$, and a tolerance level~$\delta \in \mathbb{R}_+$. \\
	\textbf{Output:} an optimal solution $\hat{\mathbf{x}}^*_p \in X$ and the optimal objective function value $\hat{z}^*_p$. \\
	$LB \longleftarrow -\infty$, $UB \longleftarrow \infty$ \; 
	$\hat{\mathbf{x}}^* \longleftarrow \mathbf{x}^{\text{\tiny 0}}$, \, $\hat{\Gamma} \longleftarrow \emptyset$ \;
	\While{$UB - LB > \delta$} 
	{$(\hat{\mathbf{y}}^*, \hat{\boldsymbol{\gamma}}_l^*) \longleftarrow$ an optimal solution of [\textbf{Sub}($\hat{\mathbf{x}}^*$)] \;
		$UB' \longleftarrow$ the optimal objective function value of [\textbf{Sub}($\hat{\mathbf{x}}^*$)] \;	
		\If{$UB' < UB$} { $(\hat{\mathbf{x}}^*_p, \hat{z}^*_p) \longleftarrow (\hat{\mathbf{x}}^*, UB')$ \; $UB \longleftarrow UB'$ }
		$\hat{\Gamma} \longleftarrow \hat{\Gamma} \cup \big\{\hat{\boldsymbol{\gamma}}_l^*\big\}$ \;
		$(\hat{\mathbf{x}}^*, \hat{z}^*)$ $\longleftarrow$ an optimal solution of [\textbf{MP}($\hat{\Gamma}$)] \; $LB \longleftarrow \hat{z}^*$}
	\Return{$(\hat{\mathbf{x}}^*_p, \hat{z}^*_p)$}
	\caption{A basic Benders decomposition-based algorithm for (\ref{pessimistic approximation SAA}).}
	\label{algorithm 1}
\end{algorithm}

Then, our algorithm closely follows the reformulation and linearization approach of Zeng and An~\cite{Zeng2014} and the standard Benders decomposition approach of Israeli and Wood~\cite{Israeli2002}; see Algorithm~\ref{algorithm 1} for its pseudocode. 
First, by leveraging a dual reformulation of constraints (\ref{cons: pessimistic approximation SAA max}) for a given subset~$\hat{\Gamma} \subseteq \Gamma$, we define a \textit{master problem}
\begin{subequations} \label{master problem ambiguity-free}
	\begin{align} 
		[\textbf{MP}(\hat{\Gamma})]: \quad & \min_{\mathbf{x}, \boldsymbol{\beta}, z} \; z \\
		\mbox{\upshape s.t. } 
		& \mathbf{x} \in X \\	
		& \begin{rcases} 
			z \geq \boldsymbol{\beta}(\boldsymbol{\gamma}_l)^\top(\mathbf{f} -\mathbf{L} \mathbf{x}) \\
			\boldsymbol{\beta}(\boldsymbol{\gamma}_l) \geq \mathbf{0} \\
			\mathbf{F}^\top \boldsymbol{\beta}(\boldsymbol{\gamma}_l) - \bar{\mathbf{c}}(\boldsymbol{\gamma}_l) \geq \mathbf{0} \quad
		\end{rcases} \, \forall \, \boldsymbol{\gamma}_l \in \hat{\Gamma}
	\end{align}
\end{subequations}
and a \textit{subproblem}
\begin{align} \label{pessimistic approximation SAA second-level}
	[\textbf{Sub}(\mathbf{x})]: \quad & \max_{\mathbf{y}, \, \boldsymbol{\gamma}_l} \Big\{ \bar{\mathbf{c}}(\boldsymbol{\gamma}_l)^\top \mathbf{y}: \; \mathbf{y} \in Y(\mathbf{x}), \, \boldsymbol{\gamma}_l \in \Gamma \Big\}. \end{align}
By Assumption \textbf{A1}, the primal problem in the right-hand side of (\ref{cons: pessimistic approximation SAA max}) has a finite optimum, ensuring that the optimal objective function value in [\textbf{MP}($\hat{\Gamma}$)] is finite. Furthermore, both [\textbf{MP}($\hat{\Gamma}$)] and~[\textbf{Sub}($\mathbf{x}$)] admit single-level MILP reformulations, after applying standard linearization techniques to the products of continuous and binary decision variables; see Section \ref{sec: comp study} for further details. 

Next, in lines 3--4 of Algorithm \ref{algorithm 1}, we initialize the lower and the upper bounds, $LB$ and $UB$, respectively, along with the initial set $\hat{\Gamma}$ and the leader's interdiction decision $\hat{\mathbf{x}}^*$. 
 In lines~5--14, the upper and the lower bounds are updated sequentially by solving the subproblem~[\textbf{Sub}($\hat{\mathbf{x}}^*$)] and the master problem [\textbf{MP}($\hat{\Gamma}$)], respectively. This process is repeated until $LB$ and $UB$~ converge within a predefined tolerance $\delta$. 
Overall, Algorithm \ref{algorithm 1} converges in a finite number of iterations since $\Gamma$ is finite; see Proposition 8~in~\cite{Zeng2014}. 

\subsection{Benders decomposition-based algorithm for the risk-neutral problem}
In this section, under some additional assumptions, we provide a similar algorithm for the risk-neutral pessimistic approximation [\textbf{DRI-P}]. The following result presents a reformulation of [\textbf{DRI-P}] for the case of $\alpha_l = 0$. 
\begin{corollary} \label{corollary 2}
	Assume that the leader is risk-neutral, i.e., $\alpha_l = 0$. Then, \upshape [\textbf{DRI-P}] \itshape reads as: \upshape
	\begin{subequations} \label{pessimistic approximation risk-neutral} 
		\begin{align} 
			\min_{\mathbf{x} \in X} \; & \max_{\mathbf{y},\, \boldsymbol{\xi}_{l},\, \boldsymbol{\eta}_{l}} \Big\{\frac{1}{k_l} \sum_{k \in K_l} \big( \hat{\mathbf{c}}_l^{\text{\tiny (\textit{k})}} - \boldsymbol{\xi}_l^{\text{\tiny (\textit{k})}} \big)^\top \mathbf{y} \Big\} \label{obj: pessimistic approximation risk-neutral} \\
			\mbox{\upshape s.t. } 	& \mathbf{y} \in Y(\mathbf{x}) \label{cons: pessimistic approximation risk-neutral 1}\\
			& \begin{rcases} \mathbf{B}(\hat{\mathbf{c}}_l^{\text{\tiny (\textit{k})}} - \boldsymbol{\xi}^{\text{\tiny (\textit{k})}}_l) \leq \mathbf{b} \quad \\
				-\boldsymbol{\eta}_l^{\text{\tiny (\textit{k})}} \leq \boldsymbol{\xi}_l^{\text{\tiny (\textit{k})}} \leq \boldsymbol{\eta}_l^{\text{\tiny (\textit{k})}} \label{cons: pessimistic approximation risk-neutral 2-3} \end{rcases} \, \forall k \in K_l \\
			& \Big\Vert \frac{1}{k_l}\sum_{k \in K_l} \boldsymbol{\eta}_l^{\text{\tiny (\textit{k})}} \Big\Vert_p \leq \varepsilon_l. \label{cons: pessimistic approximation risk-neutral 4}
		\end{align}
	\end{subequations} 
	\begin{proof}
		The result follows by setting $\alpha_l = 0$, which makes $\boldsymbol{\gamma}^*_l = \mathbf{1}$ optimal in the second-level problem reformulation (\ref{second-level quadratic reformulation pessimistic}). 
	\end{proof}
\end{corollary}

We make the following additional assumption:
\begin{itemize}
	\item[\textbf{A6$'$}.] For any fixed $\mathbf{x} \in X$, the follower's feasible set is defined as:
	\[Y(\mathbf{x}) = \big\{\mathbf{y} \in \mathbb{R}^n_+: \, \tilde{\mathbf{F}}\mathbf{y} \leq \tilde{\mathbf{f}}, \; \mathbf{y} \leq \mathbf{U}(\mathbf{1} - \mathbf{x}) \big\} \]
	where $\tilde{\mathbf{F}}$ is totally unimodular, $\tilde{\mathbf{f}}$ is integral and $\mathbf{U}$ is a diagonal matrix inducing upper bounds~on~$\mathbf{y}$.
\end{itemize}

Although Assumption \textbf{A6$'$} further restricts the class of interdiction problems in [\textbf{DP}], in our pessimistic risk-neutral setting, it is still satisfied by the min-cost flow interdiction problem \cite{Smith2008} and the shortest-path interdiction problem \cite{Israeli2002}.
On a positive note, similarly to \textbf{A6}, Assumption \textbf{A6$'$}~allows to treat $\mathbf{y}$ as integer variables and, therefore, simplify the second-level problem in~(\ref{pessimistic approximation risk-neutral}). 

The subsequent steps closely follow those for Algorithm~\ref{algorithm 1}. Consequently, a detailed discussion of the Benders decomposition algorithm applied to the risk-neutral problem (\ref{pessimistic approximation risk-neutral}) is provided in Supplementary Material~B~(see Algorithm~2). 

\section{Semi-pessimistic approximation} \label{sec: semi-pessimistic}
\subsection{Problem formulation and complexity}
 As illustrated by case (\textit{ii}) in Example \ref{ex: example 1}, the pessimistic approximation [\textbf{DRI-P}] can be rather conservative and may therefore yield a poor approximation of the true follower’s policy. 
To address the information gap between [\textbf{DRI}$^*$] and [\textbf{DRI-P}], we then introduce a \textit{semi-pessimistic approximation} of the true basic model [\textbf{DRI}$^*$]; recall case (\textit{ii}) of Example~\ref{ex: example 1}. Specifically,
we relax Assumption~\textbf{A3} by introducing the following alternative assumption:
\begin{itemize}
	\item[\textbf{A3$'$.}] The leader and the follower have access to two i.i.d. training data sets given by equation (\ref{eq: data sets}).
	Furthermore, for each $k \in K_f$, the leader either knows that 
	\begin{equation} \label{eq: sample-wise uncertainty sets} \nonumber
		\hat{\mathbf{c}}_f^{\text{\tiny (\textit{k})}} \in 	\hat{S}_l^{\text{\tiny (\textit{k})}} := \Big\{\mathbf{c} \in \mathbb{R}^n: \; \underline{\mathbf{b}}^{\text{\tiny (\textit{k})}} \leq \mathbf{c} \leq \overline{\mathbf{b}}^{\text{\tiny (\textit{k})}} \Big\} \subseteq S
	\end{equation} 
	or that $\hat{\mathbf{c}}_f^{\text{\tiny (\textit{k})}} \in \hat{\mathbf{C}}_l$.
\end{itemize}

According to Assumption \textbf{A3$'$}, the leader is aware about a part of the follower's data set $\hat{\mathbf{C}}_f$, whereas the remaining part is known to satisfy some component-wise interval constraints. Formally, based on Assumption \textbf{A3$'$}, we can also define $\hat{S}_l^{\text{\tiny (\textit{k})}} := \big\{\hat{\mathbf{c}}_f^{\text{\tiny (\textit{k})}}\big\}$ whenever $\hat{\mathbf{c}}_f^{\text{\tiny (\textit{k})}} \in \hat{\mathbf{C}}_l$ and introduce an uncertainty set 
\begin{equation} \label{eq: uncertainty set}
\hat{S}_l := \hat{S}_l^{\text{\tiny (1)}} \times \ldots \times \hat{S}_l^{\text{\tiny (\itshape k}_f \text{\tiny)}}
\end{equation}
that incorporates the leader's initial information about the follower's data. 

Taking into account (\ref{eq: uncertainty set}), the semi-pessimistic approximation of [\textbf{DRI}$^*$] reads as:
\begin{subequations} \label{bilevel distributionally robust problem semi-pessimistic} 
	\begin{align}
		\mbox{[\textbf{DRI-SP}]: \quad } \hat{z}_{sp}^{*} := \min_{\mathbf{x} \in X} \; \max_{\hat{\mathbf{C}}_f \in \hat{S}_l} \; & \min_{ \mathbf{y}} \Big\{ \max_{\mathbb{Q}_l \in \mathcal{Q}_l} \;\rho_l(\mathbf{y}, \mathbb{Q}_l) \Big\}\label{obj: bilevel DR semi-pessimistic}\\
		\mbox{s.t. } 
		& \mathbf{y} \in \argmax_{\; \tilde{\mathbf{y}} \in Y(\mathbf{x})} \; \Big\{ \min_{\mathbb{Q}_f \in \mathcal{Q}_f(\hat{\mathbf{C}}_f)} \rho_f(\tilde{\mathbf{y}}, \mathbb{Q}_f) \Big\}. \label{cons: bilevel DR semi-pessimistic 2}
	\end{align}
\end{subequations} 
That is, the leader in [\textbf{DRI-SP}] assumes the worst-case possible realization of the follower's data set with respect to the uncertainty set $\hat{S}_l$. Similar to the pessimistic approximation [\textbf{DRI-P}], we first establish the computational complexity of [\textbf{DRI-SP}]. The following result holds.

\begin{theorem} \label{theorem 5} \upshape [$\textbf{DRI-SP}$] \itshape is $\Sigma^p_2$-hard. 
	\begin{proof} 
		To establish the result, we use a reduction from the robust optimistic bilevel problem with interval uncertainty~\cite{Buchheim2021} given by:
		\begin{subequations} \label{robust bilevel problem with interval uncertainty} 
		\begin{align} \max_{\mathbf{x} \in \tilde{X}} \, \min_{\mathbf{c}_f \in \tilde{S}} & \max_{\mathbf{y}}
				\; \mathbf{c}_l^\top \mathbf{y} \\ 
				\mbox{s.t. } & \mathbf{y} \in \argmax_{\,\tilde{\mathbf{y}} \in \tilde{Y}(\mathbf{x})} \; \mathbf{c}_f^{\top} \tilde{\mathbf{y}}, 
		\end{align}
	 \end{subequations}
	 where $\tilde{X} = \{0, 1\}^m$,
	 \[\tilde{Y}(\mathbf{x}) = \{\mathbf{y} \in [0, 1]^n: \; \tilde{\mathbf{F}}\mathbf{y} + \tilde{\mathbf{L}} \mathbf{x} 
	 \leq \tilde{\mathbf{f}} \} \; \mbox{ and } \; \tilde{S} = \{\mathbf{c} \in \mathbb{R}^n: \; -\mathbf{1} \leq \mathbf{c} \leq \mathbf{1}\}.\]
	 The latter problem 
	 is proved to be $\Sigma^p_2$-hard; see Theorem~1 in~\cite{Buchheim2021}.

 To construct an instance of [\textbf{DRI-SP}], we set $X = \tilde{X}$, $Y(\mathbf{x}) = \tilde{Y}(\mathbf{x})$, $S = \hat{S}_l = \tilde{S}$ and assume that both decision-makers are risk-neutral and ambiguity-free, i.e., \(\alpha_l = \alpha_f = 0\) and \(\varepsilon_l = \varepsilon_f = 0\). In addition, let the leader and the follower each observe unique samples,~$\hat{\mathbf{c}}_l := -\mathbf{c}_l \in S$ and~$\hat{\mathbf{c}}_f := \mathbf{c}_f \in S$, respectively, with \(\hat{\mathbf{c}}_f\) being unknown to the leader.
		
		Given these prerequisites, (\ref{robust bilevel problem with interval uncertainty}) reduces to the following instance of [\textbf{DRI-SP}]: 
		\begin{subequations} \label{DRI-SP instance} 
			\begin{align} -\min_{\mathbf{x} \in X} \, \max_{\hat{\mathbf{c}}_f \in \hat{S}_l} & \min_{\mathbf{y}}
				\; \hat{\mathbf{c}}_l^\top \mathbf{y} \\ 
				\mbox{s.t. } & \mathbf{y} \in \argmax_{\,\tilde{\mathbf{y}} \in Y(\mathbf{x})} \; \hat{\mathbf{c}}_f^{\top} \tilde{\mathbf{y}}, 
			\end{align}
		\end{subequations}
		where we additionally reverse the sign of the leader's objective function. 
		 Without loss of generality, we may assume that $\hat{\mathbf{c}}_l \in S = \tilde{S}$, since scaling $\hat{\mathbf{c}}_l$ does not affect the set of optimal solutions in (\ref{DRI-SP instance}). Hence, (\ref{robust bilevel problem with interval uncertainty}) reduces to an instance of [\textbf{DRI-SP}], and the result follows. 	
	\end{proof}
\end{theorem}

\subsection{Scenario-based discretization}
Since the semi-pessimistic approximation is $\Sigma_2^p$-hard, 
we propose to use a discretization of the uncertainty set (\ref{eq: uncertainty set}), based on a finite number of scenarios for the follower's data set $\hat{\mathbf{C}}_f$; recall our discussion in Section \ref{subsec: approach and contributions}. That is, motivated by the sampling-based approach for robust optimization problems proposed by Calafiore and Campi~\cite{Calafiore2005, Calafiore2006}, we define a set of $r_l \in \mathbb{N}$ scenarios 
\begin{equation} \label{eq: scenarios}
\hat{S}'_l = \Big\{\hat{\mathbf{C}} ^{\text{\tiny (\textit{r})}}_f \in \hat{S}_l, \; r \in R_l = \{1, \ldots, r_l\}\Big\},
\end{equation}
 generated from a probability distribution $\mathbb{P} \in \mathcal{Q}_0(\hat{S}_l)$ that is fully supported on $\hat{S}_l$. 
Then, instead of directly solving~[\textbf{DRI-SP}], we address its \textit{scenario-based discretization} given~by:
\begin{subequations} \label{bilevel distributionally robust problem semi-pessimistic approximation} 
	\begin{align}
		\mbox{[\textbf{DRI-SP$'$}]: \quad } \hat{z}_{sp}' = \min_{\mathbf{x} \in X} \; \max_{\hat{\mathbf{C}}_f \in \hat{S}'_l } \; & \min_{ \mathbf{y}} \Big\{ \max_{\mathbb{Q}_l \in \mathcal{Q}_l} \;\rho_l(\mathbf{y}, \mathbb{Q}_l) \Big\}\label{obj: bilevel DR semi-pessimistic approximation}\\
		\mbox{s.t. } 
		& \mathbf{y} \in \argmax_{\; \tilde{\mathbf{y}} \in Y(\mathbf{x})} \; \Big\{ \min_{\mathbb{Q}_f \in \mathcal{Q}_f(\hat{\mathbf{C}}_f)} \rho_f(\tilde{\mathbf{y}}, \mathbb{Q}_f) \Big\}. \label{cons: bilevel DR semi-pessimistic approximation 2}
	\end{align}
\end{subequations} 


While the relation between [\textbf{DRI-SP$'$}] and [\textbf{DRI-SP}] is analyzed in detail later, we first demonstrate that [\textbf{DRI-SP$'$}] admits a single-level MILP reformulation of polynomial size. 
\begin{theorem} \label{theorem 6}
Let Assumptions \textbf{A1}, \textbf{A2} and \textbf{A3$\,'$} hold, and the ambiguity sets in \upshape (\ref{ambiguity set}) \itshape be defined in terms of $p$-norm with $p \in \{1, \infty\}$.
Then, the discretized semi-pessimistic approximation \upshape[\textbf{DRI-SP$'$}] \itshape admits the following MILP reformulation: \upshape
\begin{subequations} \label{semi-pessimistic reformulation MILP} 
	\begin{align}
		\hat{z}_{sp}' = \min_{\mathbf{x},\mathbf{y}, \boldsymbol{\nu}, \boldsymbol{s}, \boldsymbol{\mu} ,\boldsymbol{\beta}, \boldsymbol{\gamma}, \lambda, t, z} \; & z \\
		\mbox{\upshape s.t. } 	& \mathbf{x} \in X \\
		& \begin{rcases} 
			\hat{\mathbf{c}}_l^{\text{\tiny (\textit{k})}\top} \mathbf{y}^{\text{\tiny (\textit{r})}} - t^{\text{\tiny (\textit{r})}}_l + \boldsymbol{\Delta}_l^{\text{\tiny (\textit{k})}\top} \boldsymbol{\nu}_l^{\text{\tiny (\textit{k,r})}} \leq s^{\text{\tiny (\textit{k,r})}}_{l} \quad\\	
			\Vert \mathbf{B}^\top \boldsymbol{\nu}_l^{\text{\tiny (\textit{k,r})}} - \mathbf{y}^{\text{\tiny (\textit{r})}} \Vert_{q} \leq \lambda^{\text{\tiny (\textit{r})}}_l \\
			\boldsymbol{\nu}_l^{\text{\tiny (\textit{k,r})}} \geq \mathbf{0}, \;
			s^{\text{\tiny (\textit{k,r})}}_{l} \geq 0
		\end{rcases} \forall{k} \in K_l, \; \forall r \in R_l \\
		& \begin{rcases}
				z \geq t^{\text{\tiny (\textit{r})}}_l + \frac{1}{1 - \alpha_l} \big(\lambda^{\text{\tiny (\textit{r})}}_l \varepsilon_l + \frac{1}{k_l} \sum_{k \in K_l} s^{\text{\tiny (\textit{k,r})}}_{l} \big) \\
		\mbox{\upshape (\ref{cons: follower's LP 1})--(\ref{cons: follower's LP 6}), (\ref{cons: follower's LP dual 1})--(\ref{cons: follower's LP dual 7})} \\
			(\mathbf{f} - \mathbf{L}\mathbf{x})^\top \boldsymbol{\beta}^{\text{\tiny (\textit{r})}}_f = t^{\text{\tiny (\textit{r})}}_f - \frac{1}{1 - \alpha_f}\big(\varepsilon_f \lambda^{\text{\tiny (\textit{r})}}_f + \frac{1}{k_f} \; \sum_{k \in K_f} s^{\text{\tiny (\textit{k,r})}}_{f}\big) \quad \end{rcases} \forall r \in R_l,
	\end{align}
\end{subequations} 
\itshape where \upshape $\boldsymbol{\Delta}_l^{\text{\tiny (\textit{k})}} = \mathbf{b} - \mathbf{B} \hat{\mathbf{c}}_l^{\text{\tiny (\textit{k})}} \geq \mathbf{0}$ \itshape for each $k \in K_l$ and $\frac{1}{p} + \frac{1}{q} = 1$. 
\begin{proof} See Supplementary Material A.
\end{proof}
\end{theorem} 


 Next, we observe that [\textbf{DRI-SP}$'$] is an \textit{inner} \textit{approximation} of [\textbf{DRI-SP}]. That is, for a fixed number of scenarios $r_l$, 
it is possible for the optimal objective function value of~[\textbf{DRI}$^*$] to exceed that of [\textbf{DRI-SP$'$}], i.e.,~$z_b^* > \hat{z}'_{sp}$. 
 In this regard, Calafiore and Campi~\cite{Calafiore2005, Calafiore2006} establish finite-sample probabilistic guarantees for scenario-based approximations of convex optimization problems under uncertainty. However, the probabilistic guarantees from \cite{Calafiore2005, Calafiore2006} are not applicable to [\textbf{DRI-SP}], as the feasible set $X$ is discrete and, moreover, the objective function in [\textbf{DRI-SP}] is not convex in $\mathbf{x}$ for a fixed realization of uncertainty~$\hat{\mathbf{C}}_f \in \hat{S}_l$. We illustrate this non-convexity with the following~example.

\begin{example} \upshape 
Consider the instance of [\textbf{DRI-SP}] used in the proof of Theorem \ref{theorem 6}; see equation~(\ref{DRI-SP instance}). Then, for a given $\hat{\mathbf{c}}_f \in \hat{S}_l$, let 
	\begin{align*} \varphi(\mathbf{x}):= & \min_{\mathbf{y}}
		\; \hat{\mathbf{c}}_l^\top \mathbf{y} \\ 
		\mbox{s.t. } & \mathbf{y} \in \argmax_{\,\tilde{\mathbf{y}} \in Y(\mathbf{x})} \; \hat{\mathbf{c}}_f^{\top} \tilde{\mathbf{y}}.
	\end{align*}
Additionally, let $\hat{\mathbf{c}}_l = (-1, 1)^\top$, $\hat{\mathbf{c}}_f = (1, 1)^\top$ and 
\[Y(\mathbf{x}) = \{\mathbf{y} \in [0, 1]^2: \; y_1 + y_2 \leq 1, \quad y_i \leq 1 - x_i \quad \forall \, i \in \{1, 2\}\}.\]
We observe that $\varphi(0, 0) = -1$, $\varphi(1, 1) = 0$ and $\varphi(\frac{1}{2}, \frac{1}{2}) = 0$. Therefore,
\[\varphi(\frac{1}{2}, \frac{1}{2}) = 0 > \frac{1}{2} \varphi(0, 0) + \frac{1}{2} \varphi(1, 1) = -\frac{1}{2},\]
and the function $\varphi(\mathbf{x})$ is not convex in $\mathbf{x}$. \hfill$\square$
\end{example}

On a positive note, even in the absence of convexity, we establish that our discretization [\textbf{DRI-SP}$'$] is \textit{asymptotically robust} with respect to the true basic model [\textbf{DRI}$^*$]. In other words, it is proved that, as the number of scenarios, $r_l$, approaches infinity, 
the inequality $\hat{z}'_{sp} \geq z_b^*$ holds with probability one. 


\begin{theorem} \label{theorem 7}
Let Assumptions \textbf{A1}, \textbf{A2} and \textbf{A3$\,'$} hold. Also, let $\hat{z}'_{sp}$ and $z_b^*$ denote the optimal objective function values of \upshape [\textbf{DRI-SP$'$}] \itshape and \upshape [\textbf{DRI}$^*$], \itshape respectively. 
Then, as $r_l$ tends to infinity, for any scenario-generating distribution $\mathbb{P} \in \mathcal{Q}_0(\hat{S}_l)$ and $\delta' > 0$, the~inequality 
\begin{equation} \label{eq: scenario-based robustness}
\hat{z}'_{sp} \geq z_b^* - \delta'
\end{equation} holds almost surely.

\begin{proof} See Supplementary Material A. 
\end{proof}
\end{theorem}

 Overall, we note that the proposed scenario-based discretization~[\textbf{DRI-SP$'$}] is a \textit{heuristic} that enjoys only asymptotic performance guarantees. Nevertheless, as we demonstrate later in our numerical experiments, [\textbf{DRI-SP$'$}] substantially outperforms existing benchmark approaches and demonstrates strong out-of-sample performance even when only a small number of scenarios is employed.

Alternatively, future research could explore the use of existing cutting-set methods for robust convex optimization problems; see, e.g., \cite{Mutapcic2009}. Rather than relying on randomly generated scenarios, these methods iteratively augment the set of scenarios with the \textit{worst-case} scenarios under a given decision. In the context of [\textbf{DRI-SP}], such an approach would require \textit{repeatedly} solving the second-level problem for a fixed leader's decision $\mathbf{x}$. The resulting subproblems, however, involve bilinear constraints and are nonconvex, which poses significant computational challenges. 


\section{Computational study} \label{sec: comp study}
Building upon Example~\ref{ex: example 1}, our computational study examines the basic model [\textbf{DRI}], along with the scenario-based and pessimistic approximations, [\textbf{DRI-SP$'$}] and [\textbf{DRI-P}]. The remainder of this section is organized as follows. First, Section \ref{subsec: test instances and performance measures} describes our test instances and performance measures. 
In~Section~\ref{subsec: analysis basic}, we consider the basic model [\textbf{DRI}] and analyze out-of-sample errors arising from the decision-makers’ limited knowledge of the true distribution $\mathbb{Q}^*$, under various risk preferences. In~Section~\ref{subsec: analysis pessimistic and semi-pessimistic}, by using approximations [\textbf{DRI-SP$'$}] and [\textbf{DRI-P}] of the true basic model~[\textbf{DRI$^*$}], we assess both in-sample and out-of-sample errors of the leader due to its incorrect perception of the true follower's policy. Finally, a brief analysis of running times for all the aforementioned models is provided in~Section~\ref{subsec: running times}, followed by a brief summary in Section \ref{subsec: summary}. 

\subsection{Test instances and performance measures} \label{subsec: test instances and performance measures}
\textbf{Test instances.}
Similar to Zare~et~al.~\cite{Zare2019}, we consider a class of general interdiction problems defined as:
\begin{equation} \label{packing interdiction problem}
	\min_{\mathbf{x} \in X} \; \max_{\mathbf{y} \in Y(\mathbf{x}) } \;	\mathbf{c}^\top \mathbf{y},
\end{equation}	
where
\begin{align} \label{eq: feasible sets packing}
	& X = \big\{\mathbf{x} \in \{0, 1\}^n: \mathbf{H} \mathbf{x} \leq \mathbf{h}\ \big\} \, \, \mbox{ and } \, \, Y(\mathbf{x}) = \big\{\mathbf{y} \in \mathbb{R}_+^n: \tilde{\mathbf{F}}\mathbf{y} \leq \tilde{\mathbf{f}}, \; \mathbf{y} \leq \mathbf{U}(\mathbf{1} - \mathbf{x}) \big\}.
\end{align}
To facilitate qualitative analysis, we set $n = 10$, $d_l = \dim(\mathbf{h}) = 1$, and $\tilde{d}_f = \dim(\tilde{\mathbf{f}}) = 10$. Furthermore, all elements of $\mathbf{H}$ and $\tilde{\mathbf{F}}$ in (\ref{eq: feasible sets packing}) are generated uniformly at random from the interval $[0.01, 1]$, whereas
\begin{align*}
h_j = 0.4 \sum_{i = 1}^n H_{ji} \quad \forall j \in \{1, \ldots, d_l\}, \, \mbox{ } \,
\tilde{f}_j = 0.4 \sum_{i = 1}^n \tilde{F}_{ji} \quad \forall j \in \{1, \ldots, \tilde{d}_f\}
\end{align*}
and $\mathbf{U} = \mathbf{I}$. Given that the profit vector $\mathbf{c}$ is also nonnegative, (\ref{packing interdiction problem}) can be viewed as an instance of the packing interdiction problem \cite{Dinitz2013}. 
In particular, it can be verified that $Y(\mathbf{x})$ is non-empty and bounded for every fixed $\mathbf{x} \in X$, which aligns with Assumption \textbf{A1}.

Next, in line with Assumption \textbf{A2}, we introduce an interval support set 
\begin{equation} \label{eq: support packing}
	S := \{\mathbf{c} \in \mathbb{R}^n: c_i \in [0.01, 1] \quad \forall i \in N := \{1, \ldots, n\}\},
\end{equation}
 and assume that each $c_i$ follows a scaled \textit{truncated normal distribution}. This distribution is obtained by drawing samples from a normal distribution with mean $\mu_i = \frac{0.5 i}{n + 1}$ and standard deviation 
\begin{equation} \label{eq: variance}
	\sigma_i = 0.05 + \frac{0.4 i}{n + 1} \quad \forall i \in N, 
\end{equation}
 truncating them to the interval $[0,1]$, and scaling to the range $[0.01,1]$. The true distribution~$\mathbb{Q}^*$ is then defined as a product of the respective marginal distributions.  Notably,  our data-generating distribution~$\mathbb{Q}^*$ is similar to the portfolio return distribution used in \cite{Esfahani2018}, where larger returns are associated with~higher~risk. 

 In our experiments, we generate $10$ random test instances of the packing interdiction problem~(\ref{packing interdiction problem}) using a fixed random seed; each instance is then evaluated across $10$ independently generated data sets sampled from~$\mathbb{Q}^*$. The core parameters of our basic model [\textbf{DRI}] are then defined as follows. By default, $\alpha_l = 0.95$ and $\alpha_f = 0.95$ represent risk-averse scenarios for the leader and the follower, respectively; risk-neutral scenarios are attained by setting~$\alpha_l = 0$ and~$\alpha_f = 0$. 
 Furthermore, we set~$p = 1$ and use the following analytical expressions for the Wasserstein radii: 
\begin{equation} \label{eq: Wassertein radii}
\varepsilon_l := \frac{\delta_l}{\sqrt{k_l}} \quad \text{and} \quad \varepsilon_f := \frac{\delta_f}{\sqrt{k_f}},
\end{equation}
where $\delta_l \in \mathbb{R}_+$ and $\delta_f \in \mathbb{R}_+$ are appropriate constants that may depend on the problem's dimension. This approach leverages the convergence rate proportional to $k_l^{-\frac{1}{2}}$ (or, respectively, $k_f^{-\frac{1}{2}} $), which is considered theoretically optimal and achievable in practice; see, e.g., \cite{Esfahani2018, Kuhn2019}. 

 Similar to the one-stage DRO model in \cite{Esfahani2018}, in the absence of the true distribution $\mathbb{Q}^*$, the parameters $\delta_l$ and $\delta_f$ in (\ref{eq: Wassertein radii}) can potentially be tuned via \textit{cross-validation} to improve the model's out-of-sample~performance. In [\textbf{DRI}], 
this would involve splitting the leader’s and the follower’s data sets into training and validation subsets, with the follower’s subsets contained in the leader’s under Assumption \textbf{A3}. The MILP reformulation (\ref{MILP reformulation}) is then solved using the training data, and the resulting leader's decision is evaluated on the validation data. Importantly, the Wasserstein radii $\varepsilon_l$ and~$\varepsilon_f$ in~[\textbf{DRI}] are interdependent, with the
optimal choice of $\varepsilon_f$ depending on the leader’s decision and vice versa.  Since addressing this interdependence would require additional modeling assumptions, our experimental evaluation simply relies on a sufficiently large data set sampled from~$\mathbb{Q}^*$. 

 Finally, to generate scenarios for the scenario-based approximation [\textbf{DRI-SP$'$}], we define the uncertainty set $\hat{S}_l$ given by equation (\ref{eq: uncertainty set}). The first $k_{lf} \leq k_f$ samples in the true follower's data~set~$\hat{\mathbf{C}}^*_f$ are assumed to be known to the leader. 
 In contrast, for each $k \in \{k_{lf} + 1, \ldots, k_f\}$ and $i \in N$,~let 
\begin{equation} \label{eq: uncertainty set packing} \nonumber
\big(\hat{\mathbf{c}}_f^{\text{\tiny (\textit{k})}}\big)_i \in \left[ c^*_{ki} - \kappa \Delta_{ki}, c^*_{ki} + \kappa (1 - \Delta_{ki}) \right] \cap [0.01, 1].
\end{equation}
Here, $c^*_{ki}$ is the nominal value used by the follower, $\kappa$ represents a noise level, and $\Delta_{ki} \in [0, 1]$ is a shift parameter. For each of the 10 test instances, $\kappa$ is set to 0 with probability $0.5$ and $\kappa = \tilde{\kappa} > 0$, otherwise. In other words, approximately 50\% of the data entries in the remaining $k_f - k_{lf}$ samples are subject to interval uncertainty.  Moreover, $\Delta_{ki}$ is selected uniformly at random from the interval~$[0, 1]$, and the sample-generating distribution $\mathbb{P}$ is supposed to be uniform on $\hat{S}_l$; see, e.g., \cite{Calafiore2005,Calafiore2006}. 

For comparison, 
we also consider the \textit{augmented basic model}; recall case (\textit{iii}) of Example~\ref{ex: example 1}. For simplicity, within the comparison, we set $k_l = k_f$, and consequently this model reduces to [\textbf{DRI}], with the follower's data set replaced by the leader's. With a slight abuse of notation, the augmented basic model is also referred to~as~[\textbf{DRI}]. 

\textbf{Computational settings.}
All experiments are performed on a PC with CPU i9-12900U and RAM 32 GB. The MILP reformulations (\ref{MILP reformulation}) and (\ref{semi-pessimistic reformulation MILP}), as well Algorithms \ref{algorithm 1} and~2, are implemented in Java using CPLEX 22.1. For both algorithms, we set $\delta = 10^{-3}$ and~$\mathbf{x}^{\text{\tiny 0}} = \mathbf{0}$. 
We also discuss our approach to linearization of (\ref{MILP reformulation}) and other related MILP reformulations in Supplementary Material B. The respective code can be found at \url{https://github.com/sk19941995/Data-driven-interdiction}.

\textbf{Performance measures.}
First, we define the \textit{relative out-of-sample loss} of the leader and the follower in the basic model [\textbf{DRI}]. Let $(\hat{\mathbf{x}}^*, \hat{\mathbf{y}}^*)$ be an optimal solution of the MILP reformulation~(\ref{MILP reformulation}). Then, the relative out-of-sample loss of the follower is defined as:
\begin{equation} \label{eq: relative loss follower}
	\mbox{RL}_f^{\text{\tiny (\textit{out})}}(\hat{\mathbf{x}}^*, \hat{\mathbf{y}}^*, \mathbb{Q}^*) = \frac{\rho_f(\hat{\mathbf{y}}^*, \mathbb{Q}^*)}{\,\max_{\mathbf{y} \in Y(\hat{\mathbf{x}}^*)} \rho_f(\mathbf{y}, \mathbb{Q}^*)}.
\end{equation}
It follows that $\mbox{RL}_f^{\text{\tiny (\textit{out})}} \leq 1$ and $\mbox{RL}_f^{\text{\tiny (\textit{out})}} = 1$ if and only if $\hat{\mathbf{y}}^*$ is an optimal solution of the full-information follower's problem with $\mathbf{x} = \hat{\mathbf{x}}^*$. In a similar way, the leader's relative out-of-sample loss is defined as:
\begin{equation} \label{eq: relative loss leader}
	\mbox{RL}_l^{\text{\tiny (\textit{out})}}(\hat{\mathbf{x}}^*, \hat{\mathbf{y}}^*, \mathbb{Q}^*) = \frac{\rho_l(\hat{\mathbf{y}}^*, \mathbb{Q}^*)}{\rho_l(\tilde{\mathbf{y}}^*, \mathbb{Q}^*)},
\end{equation}
where $\tilde{\mathbf{y}}^*$ solves the following full information leader's problem: 
\begin{subequations} \label{full information leader's problem}
	\begin{align}
		& \min_{\mathbf{x}, \mathbf{y}} \; \rho_l(\mathbf{y}, \mathbb{Q}^*) \\
		\mbox{s.t. } & \mathbf{x} \in X \\
		& \mathbf{y} \in \argmax_{\,\tilde{\mathbf{y}} \in Y(\mathbf{x})} \; \Big\{ \min_{\mathbb{Q}_f \in \mathcal{Q}_f(\hat{\mathbf{C}}^*_f)} \rho_f(\tilde{\mathbf{y}}, \mathbb{Q}_f) \Big\}. 
	\end{align}
\end{subequations}
In other words, (\ref{full information leader's problem}) is the true basic model [\textbf{DRI$^*$}], where the leader also has full information about the true distribution $\mathbb{Q}^*$. 
Since $(\hat{\mathbf{x}}^*, \hat{\mathbf{y}}^*)$ is also feasible in (\ref{full information leader's problem}), we have $\mbox{RL}_l^{\text{\tiny (\textit{out})}} \geq 1$ and $\mbox{RL}_l^{\text{\tiny (\textit{out})}} = 1$ if and only if $(\hat{\mathbf{x}}^*, \hat{\mathbf{y}}^*)$ is an optimal solution of the full information problem (\ref{full information leader's problem}). Following the earlier discussion, instead of computing the ratios in (\ref{eq: relative loss follower}) and~(\ref{eq: relative loss leader}) explicitly, we use a sample average approximation based on $1000$ samples generated from $\mathbb{Q}^*$.

Next, for [\textbf{DRI-SP}$'$] and~[\textbf{DRI-P}], the relative out-of-sample loss of the leader is defined as follows. Given that~$\hat{\mathbf{x}}^* \in X$ is the associated optimal leader's decision, we introduce
\begin{equation} \label{eq: relative loss leader pessimistic}
	\mbox{RL}^{\text{\tiny (\textit{out})}}_{l}(\hat{\mathbf{x}}^*, \mathbb{Q}^*) = \frac{\rho_l(\hat{\mathbf{y}}^*, \mathbb{Q}^*)}{\rho_l(\tilde{\mathbf{y}}^*, \mathbb{Q}^*)},
\end{equation}
where $\tilde{\mathbf{y}}^*$ is an optimal solution of (\ref{full information leader's problem}) and $\hat{\mathbf{y}}^*$ is an optimal solution of (\ref{full information leader's problem}) with an additional constraint $\mathbf{x} = \hat{\mathbf{x}}^*$. 
Finally, the \textit{relative in-sample loss} for any approximation of the true basic model~[\textbf{DRI}$^*$] is defined as:
\begin{equation} \label{eq: relative in-sample loss} 	\mbox{RL}^{\text{\tiny (\textit{in})}}_{l} = \frac{\hat{z}^*}{z^*_{b}}.
\end{equation}
Here, $z^*_{b}$ is the optimal objective function value of [\textbf{DRI$^*$}], whereas $\hat{z}^*$ is that of its approximation,~[\textbf{DRI-SP$'$}] or [\textbf{DRI-P}]. The in-sample and the out-of-sample relative losses for the augmented basic model [\textbf{DRI}] are defined precisely in the same way as for the abovementioned~approximations.

\subsection{Analysis of basic model} \label{subsec: analysis basic}
\begin{figure}[h]
	\centering	\scalebox{0.8}{
		\begin{subfigure}{0.49\textwidth}
			\centering
			\begin{tikzpicture}
				\begin{axis}[axis lines=box,
					xlabel={$\delta_f$},
					ylabel={Relative out-of-sample loss},
					xmode=log,
					log basis x=10,
					xmin=0.0001, xmax=10, 
					ymin=0.5, ymax=1.2,
					grid=major,
					legend pos=north east,
					xtick = {0.0001, 0.001, 0.01, 0.1, 1, 10}
					]
					
					\addplot[color=black, thick, mark=triangle*, mark size = 2pt] coordinates {
						(0.0001,0.834) (0.0002,0.836) (0.0003,0.836) (0.0004,0.839) (0.0007,0.841)
						(0.0011,0.840) (0.0018,0.843) (0.0029,0.844) (0.0046,0.853) (0.0075,0.864)
						(0.0121,0.873) (0.0196,0.879) (0.0316,0.881) (0.0511,0.891) (0.0825,0.877)
						(0.1,0.883) (0.1334,0.882) (0.2154,0.846) (0.3481,0.765) (0.5623,0.733)
						(0.9085,0.738) (1.4678,0.735) (2.3714,0.738) (3.8312,0.731)
						(6.1897,0.736) (10.0,0.732)
					};
					\addlegendentry{Risk-averse follower}
					
					\addplot[color=darkbluegray, thick, mark=square*, mark size=1.5pt] coordinates {
						(0.0001,0.994) (0.0002,0.994) (0.0003,0.994) (0.0004,0.994) (0.0007,0.994)
						(0.0011,0.994) (0.0018,0.994) (0.0029,0.994) (0.0046,0.994) (0.0075,0.994)
						(0.0121,0.994) (0.0196,0.994) (0.0316,0.994) (0.0511,0.994) (0.0825,0.994)
						(0.1,0.994) (0.1334,0.994) (0.2154,0.992) (0.3481,0.988) (0.5623,0.986)
						(0.9085,0.968) (1.4678,0.932) (2.3714,0.892) (3.8312,0.834)
						(6.1897,0.739) (10.0,0.696)
					};
					\addlegendentry{Risk-neutral follower}
					
					\path[name path=upperA] plot coordinates {
						(0.0001,0.834+0.086) (0.0002,0.836+0.086) (0.0003,0.836+0.087) (0.0004,0.839+0.084) (0.0007,0.841+0.084)
						(0.0011,0.840+0.086) (0.0018,0.843+0.083) (0.0029,0.844+0.085) (0.0046,0.853+0.085) (0.0075,0.864+0.079)
						(0.0121,0.873+0.072) (0.0196,0.879+0.070) (0.0316,0.881+0.067) (0.0511,0.891+0.062) (0.0825,0.877+0.050)
						(0.1,0.883+0.046) (0.1334,0.882+0.043) (0.2154,0.846+0.061) (0.3481,0.765+0.115) (0.5623,0.733+0.118)
						(0.9085,0.738+0.117) (1.4678,0.735+0.116) (2.3714,0.738+0.120) (3.8312,0.731+0.116)
						(6.1897,0.736+0.118) (10.0,0.732+0.117)
					};
					\path[name path=lowerA] plot coordinates {
						(0.0001,0.834-0.086) (0.0002,0.836-0.086) (0.0003,0.836-0.087) (0.0004,0.839-0.084) (0.0007,0.841-0.084)
						(0.0011,0.840-0.086) (0.0018,0.843-0.083) (0.0029,0.844-0.085) (0.0046,0.853-0.085) (0.0075,0.864-0.079)
						(0.0121,0.873-0.072) (0.0196,0.879-0.070) (0.0316,0.881-0.067) (0.0511,0.891-0.062) (0.0825,0.877-0.050)
						(0.1,0.883-0.046) (0.1334,0.882-0.043) (0.2154,0.846-0.061) (0.3481,0.765-0.115) (0.5623,0.733-0.118)
						(0.9085,0.738-0.117) (1.4678,0.735-0.116) (2.3714,0.738-0.120) (3.8312,0.731-0.116)
						(6.1897,0.736-0.118) (10.0,0.732-0.117)
					};
					\addplot[black!10, fill opacity=0.5] fill between[of=upperA and lowerA];
					
					\path[name path=upperB] plot coordinates {
						(0.0001,0.994+0.009) (0.0002,0.994+0.009) (0.0003,0.994+0.009) (0.0004,0.994+0.009) (0.0007,0.994+0.009)
						(0.0011,0.994+0.009) (0.0018,0.994+0.009) (0.0029,0.994+0.009) (0.0046,0.994+0.009) (0.0075,0.994+0.009)
						(0.0121,0.994+0.009) (0.0196,0.994+0.009) (0.0316,0.994+0.009) (0.0511,0.994+0.008) (0.0825,0.994+0.008)
						(0.1,0.994+0.008) (0.1334,0.994+0.009) (0.2154,0.992+0.011) (0.3481,0.988+0.015) (0.5623,0.986+0.015)
						(0.9085,0.968+0.024) (1.4678,0.932+0.024) (2.3714,0.892+0.028) (3.8312,0.834+0.044)
						(6.1897,0.739+0.034) (10.0,0.696+0.057)
					};
					\path[name path=lowerB] plot coordinates {
						(0.0001,0.994-0.009) (0.0002,0.994-0.009) (0.0003,0.994-0.009) (0.0004,0.994-0.009) (0.0007,0.994-0.009)
						(0.0011,0.994-0.009) (0.0018,0.994-0.009) (0.0029,0.994-0.009) (0.0046,0.994-0.009) (0.0075,0.994-0.009)
						(0.0121,0.994-0.009) (0.0196,0.994-0.009) (0.0316,0.994-0.009) (0.0511,0.994-0.008) (0.0825,0.994-0.008)
						(0.1,0.994-0.008) (0.1334,0.994-0.009) (0.2154,0.992-0.011) (0.3481,0.988-0.015) (0.5623,0.986-0.015)
						(0.9085,0.968-0.024) (1.4678,0.932-0.024) (2.3714,0.892-0.028) (3.8312,0.834-0.044)
						(6.1897,0.739-0.034) (10.0,0.696-0.057)
					};
					\addplot[darkbluegray!40, fill opacity=0.5] fill between[of=upperB and lowerB];
				\end{axis}
			\end{tikzpicture}
			\caption{$k_l = k_f = 30$, $\delta_l = 0.1$ and the leader is risk-averse.}
			\label{fig: experiment 3}
		\end{subfigure}
		\hfill
		\begin{subfigure}{0.49\textwidth}
			\centering
			\begin{tikzpicture}
				\begin{axis}[axis lines=box,
					xlabel={$\delta_l$},
					ylabel={Relative out-of-sample loss},
					xmode=log,
					log basis x=10,
					xmin=0.0001, xmax=10, 
					ymin=0.9, ymax=1.6,
					grid=major,
					legend pos=north east,
					xtick = {0.0001, 0.001, 0.01, 0.1, 1, 10}
					]
					
					\addplot[color=black, thick, mark=triangle*, mark size = 2pt] coordinates {
						(0.0001,1.031) (0.0002,1.031) (0.0003,1.031) (0.0004,1.031) (0.0007,1.029)
						(0.0011,1.029) (0.0018,1.029) (0.0029,1.029) (0.0046,1.029) (0.0075,1.029)
						(0.0121,1.026) (0.0196,1.024) (0.0316,1.024) (0.0511,1.026) (0.0825,1.026)
						(0.1,1.024) (0.1334,1.024) (0.2154,1.031) (0.3481,1.036) (0.5623,1.044)
						(0.9085,1.181) (1.4678,1.229) (2.3714,1.229) (3.8312,1.229)
						(6.1897,1.229) (10.0,1.229)
					};
					\addlegendentry{Risk-averse leader}
					
					\addplot[color=darkbluegray, thick, mark=square*, mark size=1.5pt] coordinates {
						(0.0001,1.009) (0.0002,1.009) (0.0003,1.009) (0.0004,1.009) (0.0007,1.009)
						(0.0011,1.009) (0.0018,1.009) (0.0029,1.009) (0.0046,1.009) (0.0075,1.009)
						(0.0121,1.009) (0.0196,1.009) (0.0316,1.009) (0.0511,1.009) (0.0825,1.009)
						(0.1,1.009) (0.1334,1.010) (0.2154,1.010) (0.3481,1.011) (0.5623,1.012)
						(0.9085,1.013) (1.4678,1.013) (2.3714,1.017) (3.8312,1.026)
						(6.1897,1.039) (10.0,1.062)
					};
					\addlegendentry{Risk-neutral leader}
					
					\path[name path=upperA] plot coordinates {
						(0.0001,1.031+0.035) (0.0002,1.031+0.035) (0.0003,1.031+0.035) (0.0004,1.031+0.035) (0.0007,1.029+0.033)
						(0.0011,1.029+0.033) (0.0018,1.029+0.033) (0.0029,1.029+0.033) (0.0046,1.029+0.033) (0.0075,1.029+0.033)
						(0.0121,1.026+0.030) (0.0196,1.024+0.027) (0.0316,1.024+0.028) (0.0511,1.026+0.029) (0.0825,1.026+0.029)
						(0.1,1.024+0.027) (0.1334,1.024+0.027) (0.2154,1.031+0.035) (0.3481,1.036+0.038) (0.5623,1.044+0.043)
						(0.9085,1.181+0.114) (1.4678,1.229+0.119) (2.3714,1.229+0.119) (3.8312,1.229+0.119)
						(6.1897,1.229+0.119) (10.0,1.229+0.119)
					};
					\path[name path=lowerA] plot coordinates {
						(0.0001,1.031-0.035) (0.0002,1.031-0.035) (0.0003,1.031-0.035) (0.0004,1.031-0.035) (0.0007,1.029-0.033)
						(0.0011,1.029-0.033) (0.0018,1.029-0.033) (0.0029,1.029-0.033) (0.0046,1.029-0.033) (0.0075,1.029-0.033)
						(0.0121,1.026-0.030) (0.0196,1.024-0.027) (0.0316,1.024-0.028) (0.0511,1.026-0.029) (0.0825,1.026-0.029)
						(0.1,1.024-0.027) (0.1334,1.024-0.027) (0.2154,1.031-0.035) (0.3481,1.036-0.038) (0.5623,1.044-0.043)
						(0.9085,1.181-0.114) (1.4678,1.229-0.119) (2.3714,1.229-0.119) (3.8312,1.229-0.119)
						(6.1897,1.229-0.119) (10.0,1.229-0.119)
					};
					\addplot[black!10, fill opacity=0.5] fill between[of=upperA and lowerA];
					
					\path[name path=upperB] plot coordinates {
						(0.0001,1.009+0.013) (0.0002,1.009+0.013) (0.0003,1.009+0.013) (0.0004,1.009+0.013) (0.0007,1.009+0.013)
						(0.0011,1.009+0.013) (0.0018,1.009+0.013) (0.0029,1.009+0.013) (0.0046,1.009+0.013) (0.0075,1.009+0.013)
						(0.0121,1.009+0.013) (0.0196,1.009+0.013) (0.0316,1.009+0.013) (0.0511,1.009+0.013) (0.0825,1.009+0.013)
						(0.1,1.009+0.013) (0.1334,1.010+0.014) (0.2154,1.010+0.013) (0.3481,1.011+0.015) (0.5623,1.012+0.016)
						(0.9085,1.013+0.018) (1.4678,1.013+0.018) (2.3714,1.017+0.022) (3.8312,1.026+0.031)
						(6.1897,1.039+0.042) (10.0,1.062+0.056)
					};
					\path[name path=lowerB] plot coordinates {
						(0.0001,1.009-0.013) (0.0002,1.009-0.013) (0.0003,1.009-0.013) (0.0004,1.009-0.013) (0.0007,1.009-0.013)
						(0.0011,1.009-0.013) (0.0018,1.009-0.013) (0.0029,1.009-0.013) (0.0046,1.009-0.013) (0.0075,1.009-0.013)
						(0.0121,1.009-0.013) (0.0196,1.009-0.013) (0.0316,1.009-0.013) (0.0511,1.009-0.013) (0.0825,1.009-0.013)
						(0.1,1.009-0.013) (0.1334,1.010-0.014) (0.2154,1.010-0.013) (0.3481,1.011-0.015) (0.5623,1.012-0.016)
						(0.9085,1.013-0.018) (1.4678,1.013-0.018) (2.3714,1.017-0.022) (3.8312,1.026-0.031)
						(6.1897,1.039-0.042) (10.0,1.062-0.056)
					};
					\addplot[darkbluegray!40, fill opacity=0.5] fill between[of=upperB and lowerB];
				\end{axis}
			\end{tikzpicture}
			\caption{$k_l = k_f = 30$, $\delta_f = 0.1$ and the follower is risk-averse.}
			\label{fig: experiment 4}
	\end{subfigure}}
	\caption{The average relative out-of-sample loss (with MADs) of the follower (\ref{eq: relative loss follower}) (a) and the leader (\ref{eq: relative loss leader}) (b) as a function of $\delta_f$ and~$\delta_l$, respectively, evaluated over 100 random test instances.}
	\label{fig: experiments 3 and 4}
\end{figure}
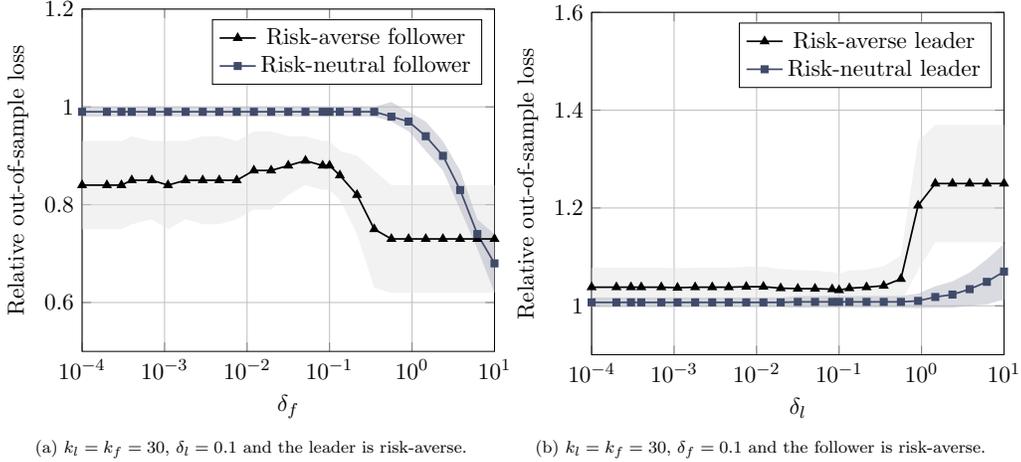

 In the basic model [\textbf{DRI}], we define the follower's data set $\hat{\mathbf{C}}_f$~as a subset of the leader's data set~$\hat{\mathbf{C}}_l$; recall Assumption~\textbf{A3}. 
Initially, we set $k_l = k_f = 30$ and examine the follower's relative loss~(\ref{eq: relative loss follower}) and the leader's relative loss (\ref{eq: relative loss leader}) as functions of $\delta_f$ and~$\delta_l$, respectively; recall equation (\ref{eq: Wassertein radii}). Thus, Figures \ref{fig: experiment 3} and \ref{fig: experiment 4} present the average relative loss of the follower and the leader, respectively, along with mean absolute deviations (MADs), computed over 100~random test instances; recall that we generate $10$ random test instances of the packing interdiction problem (\ref{packing interdiction problem}), with $10$ different data sets for each instance.

We make the following observations:
\begin{itemize} \setlength{\itemsep}{0.05cm} 
	\setlength{\parskip}{0.05cm} 
	\setlength{\topsep}{2pt} 
	\item Similar to the one-stage DRO model in \cite{Esfahani2018}, for both the risk-averse leader and the risk-averse follower, the relative out-of sample loss first tends to improve and then deteriorates with the increase of the Wasserstein radius. In other words, both decision-makers alternate between the sample average ($\varepsilon_l = 0$ or $\varepsilon_f = 0$) and the robust regimes ($\varepsilon_l = \infty$ or $\varepsilon_f = \infty$), with a potential improvement in-between these two regimes. 
	\item For both the risk-averse leader and the risk-averse follower, a \textit{critical radius} can be identified at approximately $\varepsilon_l \approx \frac{0.1}{\sqrt{30}}$ and $\varepsilon_f \approx \frac{0.1}{\sqrt{30}}$, corresponding to the best achievable out-of-sample performance. Meanwhile, the DRO approach yields a relative gain of about $5\%$ over the sample average approximation (SAA) for the risk-averse follower, but only about $1\%$ for the risk-averse leader. One possible explanation for this observation is that the leader's decision variables are binary, while the follower's decision variables are continuous.
	\item For the risk-neutral leader and the risk-neutral follower, no relative gain over the SAA is achieved. 
	This outcome is somewhat expected, given our choice of the $\ell_1$-norm; see, e.g., Lemma~2 in \cite{Ketkov2024}. 
	\item Finally, we observe that the robust regime for the risk-averse decision-makers is achieved at smaller values of $\delta_l$ and $\delta_f$ compared to the risk-neutral decision-makers. In this regard, we note that the conditional value-at-risk can be viewed as an additional layer of robustness with respect to the risk-neutral objective~function. 
\end{itemize}

\begin{figure}[h]
	\centering	\scalebox{0.8}{
		\begin{subfigure}{0.49\textwidth}
			\centering
			\begin{tikzpicture}
				\begin{axis}[axis lines=box,
					xlabel={$k_f$},
					ylabel={Relative out-of-sample loss},
					xmode=log,
					log basis x=10,
					legend pos=north east,
					grid=major,
					xmin = 10, xmax = 300,
					ymin=0.5, ymax=1.2,
					xtick={10,30,100,300},
					xticklabels={10, 30, 100, 300}
					]
					
					\addplot[color=black, thick, mark=triangle*, mark size = 2pt] coordinates {
						(5,0.711) (10,0.712) (15,0.709) (20,0.713) (30,0.712)
						(40,0.709) (50,0.716) (75,0.729) (100,0.742) (125,0.756)
						(150,0.762) (175,0.781) (200,0.835) (250,0.859) (300,0.874)
					};
					\addlegendentry{Risk-averse follower}
					
					\addplot[color=darkbluegray, thick, mark=square*, mark size=1.5pt] coordinates {
						(5,0.923) (10,0.958) (15,0.969) (20,0.981) (30,0.985)
						(40,0.990) (50,0.990) (75,0.995) (100,0.996) (125,0.997)
						(150,0.998) (175,0.998) (200,0.999) (250,0.999) (300,0.999)		
					};
					\addlegendentry{Risk-neutral follower}
					
					\path[name path=upperA] plot coordinates {
						(5,0.711+0.110) (10,0.712+0.111) (15,0.709+0.109) (20,0.713+0.111) (30,0.712+0.110)
						(40,0.709+0.111) (50,0.716+0.113) (75,0.729+0.119) (100,0.742+0.122) (125,0.756+0.130)
						(150,0.762+0.125) (175,0.781+0.108) (200,0.835+0.068) (250,0.859+0.051) (300,0.874+0.038)
					};
					\path[name path=lowerA] plot coordinates {
						(5,0.711-0.110) (10,0.712-0.111) (15,0.709-0.109) (20,0.713-0.111) (30,0.712-0.110)
						(40,0.709-0.111) (50,0.716-0.113) (75,0.729-0.119) (100,0.742-0.122) (125,0.756-0.130)
						(150,0.762-0.125) (175,0.781-0.108) (200,0.835-0.068) (250,0.859-0.051) (300,0.874-0.038)
					};
					\addplot[black!10, fill opacity=0.5] fill between[of=upperA and lowerA];
					
					\path[name path=upperB] plot coordinates {
						(5,0.923+0.042) (10,0.958+0.030) (15,0.969+0.025) (20,0.981+0.018) (30,0.985+0.014)
						(40,0.990+0.011) (50,0.990+0.011) (75,0.995+0.007) (100,0.996+0.006) (125,0.997+0.005)
						(150,0.998+0.004) (175,0.998+0.003) (200,0.999+0.002) (250,0.999+0.001) (300,0.999+0.001)
					};
					\path[name path=lowerB] plot coordinates {
						(5,0.923-0.042) (10,0.958-0.030) (15,0.969-0.025) (20,0.981-0.018) (30,0.985-0.014)
						(40,0.990-0.011) (50,0.990-0.011) (75,0.995-0.007) (100,0.996-0.006) (125,0.997-0.005)
						(150,0.998-0.004) (175,0.998-0.003) (200,0.999-0.002) (250,0.999-0.001) (300,0.999-0.001)
					};
					\addplot[darkbluegray!40, fill opacity=0.5] fill between[of=upperB and lowerB];
				\end{axis}
			\end{tikzpicture}
			\caption{$k_l = 300$, $\delta _f = 0.5$, $\delta _l = 0.1$ and the leader is risk-averse.}
			\label{fig: experiment 1}
		\end{subfigure}
		\hfill	
		\begin{subfigure}{0.49\textwidth}
			\centering
			\begin{tikzpicture}
				\begin{axis}[axis lines=box,
					xlabel={$k_l$},
					ylabel={\text{Relative out-of-sample loss}},
					xmode=log,
					log basis x=10,
					legend pos=north east,
					grid=major,
					xmin = 10, xmax = 300,
					ymin=0.9, ymax=1.3,
					xtick={10,30,100,300},
					xticklabels={10, 30, 100, 300}
					]
					
					\addplot[color=black, thick, mark=triangle*, mark size = 2pt] coordinates {
						(10,1.127) (15,1.088) (20,1.070) (30,1.067) (40,1.063)
						(50,1.048) (75,1.047) (100,1.043) (125,1.035) (150,1.034)
						(175,1.034) (200,1.028) (250,1.024) (300,1.021)
					};
					\addlegendentry{Risk-averse leader}
					
					\addplot[color=darkbluegray, thick, mark=square*, mark size=1.5pt] coordinates {
						(10,1.026) (15,1.021) (20,1.013) (30,1.013) (40,1.013)
						(50,1.010) (75,1.006) (100,1.004) (125,1.003) (150,1.004)
						(175,1.004) (200,1.004) (250,1.002) (300,1.002)
					};
					\addlegendentry{Risk-neutral leader}
					
					\path[name path=upperA] plot coordinates {
						(10,1.127+0.104) (15,1.088+0.076) (20,1.070+0.060) (30,1.067+0.063) (40,1.063+0.060)
						(50,1.048+0.048) (75,1.047+0.045) (100,1.043+0.045) (125,1.035+0.038) (150,1.034+0.038)
						(175,1.034+0.039) (200,1.028+0.034) (250,1.024+0.029) (300,1.021+0.027)
					};
					\path[name path=lowerA] plot coordinates {
						(10,1.127-0.104) (15,1.088-0.076) (20,1.070-0.060) (30,1.067-0.063) (40,1.063-0.060)
						(50,1.048-0.048) (75,1.047-0.045) (100,1.043-0.045) (125,1.035-0.038) (150,1.034-0.038)
						(175,1.034-0.039) (200,1.028-0.034) (250,1.024-0.029) (300,1.021-0.027)
					};
					\addplot[black!10, fill opacity=0.5] fill between[of=upperA and lowerA];
					
					\path[name path=upperB] plot coordinates {
						(10,1.026+0.029) (15,1.021+0.025) (20,1.013+0.017) (30,1.013+0.017) (40,1.013+0.017)
						(50,1.010+0.014) (75,1.006+0.009) (100,1.004+0.007) (125,1.003+0.005) (150,1.004+0.006)
						(175,1.004+0.006) (200,1.004+0.006) (250,1.002+0.004) (300,1.002+0.003)
					};
					\path[name path=lowerB] plot coordinates {
						(10,1.026-0.029) (15,1.021-0.025) (20,1.013-0.017) (30,1.013-0.017) (40,1.013-0.017)
						(50,1.010-0.014) (75,1.006-0.009) (100,1.004-0.007) (125,1.003-0.005) (150,1.004-0.006)
						(175,1.004-0.006) (200,1.004-0.006) (250,1.002-0.004) (300,1.002-0.003)
					};
					\addplot[darkbluegray!40, fill opacity=0.5] fill between[of=upperB and lowerB];
				\end{axis}
			\end{tikzpicture}
			\caption{$k_f = 10$, $\delta _l = 0.5$, $\delta _f = 0.1$ and the follower is risk-averse.}
			\label{fig: experiment 2}
	\end{subfigure}}
	\caption{The average relative out-of-sample loss (with MADs) of the follower (\ref{eq: relative loss follower}) (a) and the leader (\ref{eq: relative loss leader}) (b) as a function of their respective sample sizes, evaluated over 100 random test instances.}
	\label{fig: experiments 1 and 2}
\end{figure}
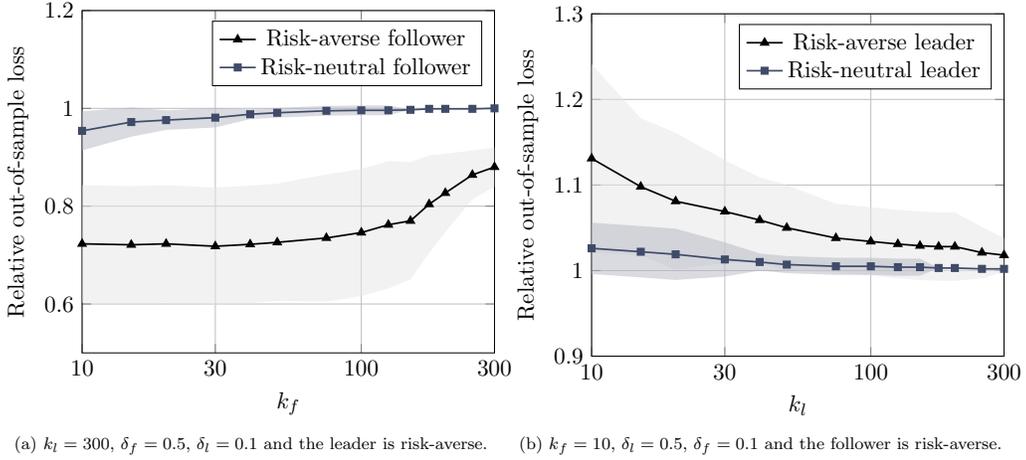

Next, we use the relative out-of-sample losses (\ref{eq: relative loss follower}) and (\ref{eq: relative loss leader}) to analyze convergence, respectively, for the leader and the follower; see Figures \ref{fig: experiment 1} and \ref{fig: experiment 2}. 
We make the following observations:
\begin{itemize}
	\setlength{\itemsep}{0.05cm} 
	\setlength{\parskip}{0.05cm} 
	\setlength{\topsep}{2pt} 
	\item In line with Lemma~\ref{lemma 2}, 
	the follower's relative loss in (\ref{eq: relative loss follower}) tends to $1$ as $k_f$ increases, irrespective of the leader's decision~$\mathbf{x} \in X$.
	\item In line with  Theorem~\ref{theorem 2}, the leader's relative loss in (\ref{eq: relative loss leader}) tends to 1 as~$k_l$ increases, even when~$k_f$ remains fixed; see equation~(\ref{eq: leader's problem convergence 3}).
	\item Based on Figures~\ref{fig: experiment 1} and~\ref{fig: experiment 2}, we conclude that the convergence is faster for the case of the risk-neutral decision-makers. This observation is rather natural, given that the selected parameters,~$\delta_f$ and~$\delta_l$, correspond to different regimes for the risk-averse and the risk-neutral decision-makers; recall Figures~\ref{fig: experiment 3}~and~\ref{fig: experiment 4}. 
\end{itemize}

 Finally, we conduct a detailed analysis of the decision-makers' risk preferences by examining different combinations of risk-neutral (N) and risk-averse (A) regimes for both the leader and the follower. For each combination of regimes and the associated optimal solution $(\hat{\mathbf{x}}^*, \hat{\mathbf{y}}^*)$ of the MILP reformulation~(\ref{MILP reformulation}), we estimate the true right-tail and the true left-tail CVaRs, $\rho_l(\hat{\mathbf{y}}^*, \mathbb{Q}^*)$ and $\rho_f(\hat{\mathbf{y}}^*, \mathbb{Q}^*)$, as well as the true mean, $\mathbb{E}_{\mathbb{Q}^*} \{\mathbf{c}^\top \hat{\mathbf{y}}^*\}$, based on 1000 samples generated from $\mathbb{Q}^*$; 
 see Table~\ref{fig: experiment 5}. 

\begin{table}[h]
	\begin{subtable}{\textwidth} \centering \onehalfspacing 
	\tiny
	\begin{tabular}{c|c|c ccc}
		\hline
		\multirow{2}{*}{Sample sizes} & \multirow{2}{*}{Risk preferences} & \multicolumn{3}{c}{True risk measure} \\
		\cline{3-5} 
		& & Left-tail CVaR & Mean & Right-tail CVaR \\ \hline
		\multirow{4}{*}{$k_l = k_f = 30$} 
		& NN & \textbf{0.37 (0.04)} & \textbf{0.94 (0.09)} & \textbf{1.61 (0.15)} \\ 
		& NA & \textbf{0.37 (0.03)} & 0.76 (0.04) & 1.20 (0.06) \\ 
		& AN & 0.40 (0.05) & 0.98 (0.10) & 1.65 (0.14) \\ 
		& AA & 0.39 (0.03) & 0.78 (0.05) & 1.22 (0.07) \\ \hline
		\multirow{4}{*}{$k_l = k_f = 300$} 
		& NN & 0.36 (0.04) & \textbf{0.94 (0.09)} & \textbf{1.61 (0.15)} \\ 
		& NA & \textbf{0.39 (0.04)} & 0.81 (0.07) & 1.29 (0.10) \\
		& AN & 0.38 (0.05) & 0.95 (0.10) & \textbf{1.61 (0.14)} \\ 
		& AA & 0.40 (0.05) & 0.82 (0.07) & 1.29 (0.09) \\ \hline	
	\end{tabular}
	\caption{}
	\label{tab: 1}
 \end{subtable}
 \vspace{3pt}
 \begin{subtable}{\textwidth} \centering \onehalfspacing 
	\tiny
	\begin{tabular}{c|c|c ccc}
		\hline
		\multirow{2}{*}{Sample sizes} & \multirow{2}{*}{Risk preferences} & \multicolumn{3}{c}{True risk measure} \\
		\cline{3-5} 
		& & Left-tail CVaR & Mean & Right-tail CVaR \\ \hline
		\multirow{4}{*}{$k_l = k_f = 30$} 
		& NN & 0.51 (0.06) & \textbf{1.06 (0.05)} & \textbf{1.73 (0.08)} \\ 
		& NA & \textbf{0.52 (0.05)} & 0.96 (0.05) & 1.48 (0.08) \\ 
		& AN & 0.66 (0.11) & 1.15 (0.08) & 1.74 (0.09) \\ 
		& AA & 0.59 (0.06) & 1.01 (0.06) & 1.49 (0.08) \\ \hline
		\multirow{4}{*}{$k_l = k_f = 300$} 
		& NN & 0.49 (0.06) & \textbf{1.06 (0.04)} & 1.75 (0.07) \\ 
		& NA & \textbf{0.53 (0.04)} & 0.99 (0.04) & 1.53 (0.08) \\
		& AN & 0.67 (0.10) & 1.14 (0.07) & \textbf{1.68 (0.06)} \\ 
		& AA & 0.63 (0.08) & 1.02 (0.06) & 1.46 (0.06) \\ \hline
	\end{tabular}
	\caption{}
	\label{tab: 2}
 \end{subtable}
	\caption{\footnotesize The average out-of-sample performance (with MADs) of the decision-makers with varying risk preferences, for $\delta_l = \delta_f = 0.1$, evaluated over 100 random test instances. The true data-generating distribution has variances either proportional (a) or \textit{inversely} proportional (b) to the mean. Symbols
	``N'' and ``A'' in column~$2$ represent risk-neutral and risk-averse policies, respectively, while the first letter refers to the leader's policy and the second letter refers to the follower's policy. The policy that is beneficial for both the leader and the follower is highlighted in bold.}
	\label{fig: experiment 5}
\end{table}

 In this experiment, we consider two classes of data-generating distributions, where $\sigma_i$, $i \in N$, is either proportional to the mean, as defined by equation (\ref{eq: variance}), or inversely proportional to the mean,~i.e., 
\begin{equation} \label{eq: inverse variance} \nonumber
	\sigma_i = 0.05 + (0.4 - \frac{0.4 i}{n + 1}) \quad \forall i \in N; 
\end{equation}
see Tables \ref{tab: 1} and \ref{tab: 2}, respectively. %
The second class of distributions 
serves to clearly illustrate the contrast between the risk-neutral and the risk-averse regimes for the leader; see, e.g., the values in rows~NN and AN. More specifically, when the variance is proportional to the mean, 
the leader's optimal solution in the risk-averse case (A) often coincides with the one in the risk-neutral case~(N), leading to similar out-of-sample~performance. For both distribution classes, we identify policy combinations that are \textit{mutually beneficial} for the leader and the follower; in Table \ref{fig: experiment 5}, these combinations are highlighted in bold. 
Furthermore, given that the risk-averse regime is shown to exhibit slower convergence, we analyze both small ($k_l = k_f = 30$) and large ($k_l = k_f = 300$) sample~sizes. 

We observe that, for small sample sizes, the risk-neutral policy 
is typically beneficial for both decision-makers. However, for large sample sizes, the neutral-averse (NA) and the averse-neutral~(AN) policies become beneficial under the left and the right-tail CVaRs, respectively. This trend can be explained by the fact that these policies are, in a sense, the ``closest'' approximations of the respective true risk measures. What is even more interesting is that the averse-averse (AA) policy, while theoretically justified, is generally not beneficial. Specifically, for large sample sizes, it is only preferable for the follower under the left-tail CVaR, and only for the leader under the right-tail~CVaR. Therefore, adopting the AA policy makes sense only if avoiding both extremely small and large profits is deemed more important than optimizing average gains.

\subsection{Analysis of semi-pessimistic and pessimistic approximations} \label{subsec: analysis pessimistic and semi-pessimistic}
\begin{figure}[h]\centering	\scalebox{0.8}{
	\begin{subfigure}{0.5\textwidth} 
		\centering
		\begin{tikzpicture}
			\begin{axis}[
				axis lines=box,
				xlabel={$r_l$},
				ylabel={Relative in-sample loss},
				xmin=1, xmax=10, 
				ymin=0.9, ymax=1.2,
				grid=major,
				legend pos=north east,
				mark size=3pt,
				xtick={2,4,6,8,10},
				ytick={0.9,0.95,1,1.05,1.1}
				]
				
			\addplot[color=black, mark=triangle*, mark size=2pt] coordinates {
			 (1,1.003) (2,1.003) (3,1.003) (4,1.003) (5,1.003)
			 (6,1.003) (7,1.003) (8,1.003) (9,1.003) (10,1.003)
			};
			\addlegendentry{\footnotesize Augmented basic model}
				
			\addplot[color=limegreen, mark=*, mark size=2pt] coordinates {
			(1,0.999) (2,1.012) (3,1.016) (4,1.019) (5,1.020)
			(6,1.023) (7,1.026) (8,1.027) (9,1.028) (10,1.032)
			};
			\addlegendentry{\footnotesize Scenario-based approximation ($\tilde{\kappa} = 0.2$)}
			
			\addplot[color=forestgreen, mark=square*, mark size=1.5pt] coordinates {
				(1,1.002) (2,1.016) (3,1.023) (4,1.027) (5,1.029)
				(6,1.031) (7,1.035) (8,1.037) (9,1.038) (10,1.041)
			};
			\addlegendentry{\footnotesize Scenario-based approximation ($\tilde{\kappa} = 0.5$)}
			
			\addplot[mark=*, mark size=2pt, dashed, color=limegreen] coordinates {
				(1,0.937) (2,0.966) (3,0.976) (4,0.989) (5,0.989)
				(6,0.989) (7,0.999) (8,0.999) (9,0.999) (10,1.000)
			};
			
			\addplot[mark=square*, mark size=1.5pt, dashed, color=forestgreen] coordinates {
				(1,0.959) (2,0.988) (3,0.998) (4,0.999) (5,1.000)
				(6,1.000) (7,1.000) (8,1.000) (9,1.000) (10,1.000)
			};
				
			\path[name path=upperC] plot coordinates {
				(1,1.003+0.030) (2,1.003+0.030) (3,1.003+0.030) (4,1.003+0.030) (5,1.003+0.030)
				(6,1.003+0.030) (7,1.003+0.030) (8,1.003+0.030) (9,1.003+0.030) (10,1.003+0.030)
				};
				\path[name path=lowerC] plot coordinates {
				(1,1.003-0.030) (2,1.003-0.030) (3,1.003-0.030) (4,1.003-0.030) (5,1.003-0.030)
				(6,1.003-0.030) (7,1.003-0.030) (8,1.003-0.030) (9,1.003-0.030) (10,1.003-0.030)
				};
				\addplot[black!10, fill opacity=0.5] fill between[of=upperC and lowerC];
				
			\path[name path=upperA] plot coordinates {
				(1,0.999+0.020) (2,1.012+0.025) (3,1.016+0.025) (4,1.019+0.027) (5,1.020+0.027)
				(6,1.023+0.030) (7,1.026+0.031) (8,1.027+0.031) (9,1.028+0.031) (10,1.032+0.035)
			};
			\path[name path=lowerA] plot coordinates {
				(1,0.999-0.020) (2,1.012-0.025) (3,1.016-0.025) (4,1.019-0.027) (5,1.020-0.027)
				(6,1.023-0.030) (7,1.026-0.031) (8,1.027-0.031) (9,1.028-0.031) (10,1.032-0.035)
			};
			\addplot[limegreen!20, fill opacity=0.5] fill between[of=upperA and lowerA];
			
			\path[name path=upperB] plot coordinates {
				(1,1.002+0.017) (2,1.016+0.026) (3,1.023+0.029) (4,1.027+0.031) (5,1.029+0.034)
				(6,1.031+0.034) (7,1.035+0.036) (8,1.037+0.037) (9,1.038+0.038) (10,1.041+0.040)
			};
			\path[name path=lowerB] plot coordinates {
				(1,1.002-0.017) (2,1.016-0.026) (3,1.023-0.029) (4,1.027-0.031) (5,1.029-0.034)
				(6,1.031-0.034) (7,1.035-0.036) (8,1.037-0.037) (9,1.038-0.038) (10,1.041-0.040)
			};
			\addplot[forestgreen!20, fill opacity=0.5] fill between[of=upperB and lowerB];	
			
			\addplot[black, thin, gray, domain=1:10, samples=2] {1};
			\end{axis}
		\end{tikzpicture}
		\caption{}
		\label{fig: experiment 6a}
	\end{subfigure}
	\hfill
	\begin{subfigure}{0.5\textwidth}
		\centering 
		\begin{tikzpicture}
			\begin{axis}[
				axis lines=box,
				xlabel={$r_l$},
				ylabel={Relative out-of-sample loss},
				xmin=1, xmax=10, 
				ymin=0.97, ymax=1.2, 
				grid=major,
				legend pos=north east,
				mark size=3pt,
				xtick={2,4,6,8,10},
				ytick={1,1.05,1.1,1.15}
				]
				
				\addplot[color=black, mark=triangle*, mark size = 2pt, dashed] coordinates {
					(1,1.039) (2,1.039) (3,1.039) (4,1.039) (5,1.039)
					(6,1.039) (7,1.039) (8,1.039) (9,1.039) (10,1.039)
				};
				\addlegendentry{\footnotesize True basic model}
				
				\addplot[color=black, mark=triangle*, mark size = 2pt] coordinates {
				(1,1.101) (2,1.101) (3,1.101) (4,1.101) (5,1.101)
				(6,1.101) (7,1.101) (8,1.101) (9,1.101) (10,1.101)
				};
				\addlegendentry{\footnotesize Augmented basic model}
				
				\addplot[color=limegreen, mark=*, mark size=2pt] coordinates {
				(1,1.052) (2,1.052) (3,1.052) (4,1.049) (5,1.048)
				(6,1.046) (7,1.049) (8,1.049) (9,1.049) (10,1.052)
				};
				\addlegendentry{\footnotesize Scenario-based approximation ($\tilde{\kappa} = 0.2$)}
				
				\addplot[color=forestgreen, mark=square*, mark size=1.5pt] coordinates {
					(1,1.069) (2,1.057) (3,1.057) (4,1.058) (5,1.063)
					(6,1.061) (7,1.057) (8,1.055) (9,1.051) (10,1.050)
				};
				\addlegendentry{\footnotesize Scenario-based approximation ($\tilde{\kappa} = 0.5$)}
				
				\path[name path=upperA] plot coordinates {
					(1,1.039+0.039) (2,1.039+0.039) (3,1.039+0.039) (4,1.039+0.039) (5,1.039+0.039)
					(6,1.039+0.039) (7,1.039+0.039) (8,1.039+0.039) (9,1.039+0.039) (10,1.039+0.039)
				};
				\path[name path=lowerA] plot coordinates {
				 (1,1.039-0.039) (2,1.039-0.039) (3,1.039-0.039) (4,1.039-0.039) (5,1.039-0.039)
				 (6,1.039-0.039) (7,1.039-0.039) (8,1.039-0.039) (9,1.039-0.039) (10,1.039-0.039)
				};
				\addplot[black!30, fill opacity=0.5] fill between[of=upperA and lowerA];
				
				\path[name path=upperB] plot coordinates {
					(1,1.101+0.086) (2,1.101+0.086) (3,1.101+0.086) (4,1.101+0.086) (5,1.101+0.086)
					(6,1.101+0.086) (7,1.101+0.086) (8,1.101+0.086) (9,1.101+0.086) (10,1.101+0.086)
				};
				\path[name path=lowerB] plot coordinates {
					(1,1.101-0.086) (2,1.101-0.086) (3,1.101-0.086) (4,1.101-0.086) (5,1.101-0.086)
					(6,1.101-0.086) (7,1.101-0.086) (8,1.101-0.086) (9,1.101-0.086) (10,1.101-0.086)
				};
				\addplot[black!10, fill opacity=0.5] fill between[of=upperB and lowerB];
				
				\path[name path=upperC] plot coordinates {
					(1,1.052+0.049) (2,1.052+0.049) (3,1.052+0.049) (4,1.049+0.048) (5,1.048+0.047)
					(6,1.046+0.046) (7,1.049+0.047) (8,1.049+0.047) (9,1.049+0.047) (10,1.052+0.050)
				};
				\path[name path=lowerC] plot coordinates {
					(1,1.052-0.049) (2,1.052-0.049) (3,1.052-0.049) (4,1.049-0.048) (5,1.048-0.047)
					(6,1.046-0.046) (7,1.049-0.047) (8,1.049-0.047) (9,1.049-0.047) (10,1.052-0.050)
				};
				\addplot[limegreen!20, fill opacity=0.5] fill between[of=upperC and lowerC];
				
				\path[name path=upperD] plot coordinates {
					(1,1.069+0.067) (2,1.057+0.057) (3,1.057+0.055) (4,1.058+0.057) (5,1.063+0.060)
					(6,1.061+0.059) (7,1.057+0.056) (8,1.055+0.056) (9,1.051+0.052) (10,1.050+0.050)
				};
				\path[name path=lowerD] plot coordinates {
					(1,1.069-0.067) (2,1.057-0.057) (3,1.057-0.055) (4,1.058-0.057) (5,1.063-0.060)
					(6,1.061-0.059) (7,1.057-0.056) (8,1.055-0.056) (9,1.051-0.052) (10,1.050-0.050)
				};
				\addplot[forestgreen!10, fill opacity=0.5] fill between[of=upperD and lowerD];
				
				\addplot[black, thin, gray, domain=1:10, samples=2] {1};
			\end{axis}
		\end{tikzpicture}
		\caption{}
		\label{fig: experiment 6b}
	\end{subfigure}}
	\caption{The average relative in-sample (\ref{eq: relative loss leader pessimistic}) (a) and out-of-sample (\ref{eq: relative in-sample loss}) (b) loss (with MADs) as a function of the number of scenarios,~$r_l$, for $k_l = k_f = 30$, $\delta_l = \delta_f = 0.1$ and $k_{lf} = 20$, evaluated over 100 random test instances. The dashed lines in (a) correspond to the empirical $5\%$ percentile of the relative in-sample loss.}
	\label{fig: experiment 6}
\end{figure}
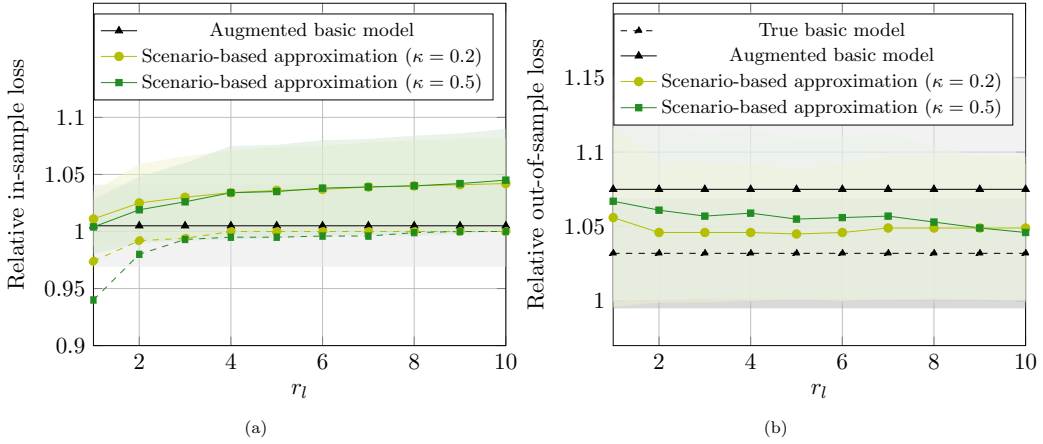

We first analyze the scenario-based approximation [\textbf{DRI-SP$'$}] in terms of the leader's relative out-of-sample and in-sample losses, (\ref{eq: relative loss leader pessimistic}) and (\ref{eq: relative in-sample loss}), as functions of the number of scenarios $r_l \in~\{1, \ldots, 10\}$; see Figure \ref{fig: experiment 6}. In particular, we set $k_{lf} = 20$ and~$k_f = 30$, i.e., one-third of the follower's data set is subject to uncertainty, and examine medium and large noise levels, $\tilde{\kappa} = 0.2$ and $\tilde{\kappa} = 0.5$. To verify robustness of [\textbf{DRI-SP$'$}], in addition to the average in-sample loss, for each $r_l$, we provide its empirical~$5\%$ percentile. For comparison, we also indicate the results for the underlying augmented basic model [\textbf{DRI}] and the true basic model [\textbf{DRI$^*$}] that are not affected by the number of~scenarios.


Based on Figure \ref{fig: experiment 6}, we make the following observations:
\begin{itemize}
	\setlength{\itemsep}{0.05cm} 
	\setlength{\parskip}{0.05cm} 
	\setlength{\topsep}{2pt} 
	\item In line with Theorem \ref{theorem 7}, [\textbf{DRI-SP$'$}] tends to become more robust in relation to [\textbf{DRI$^*$}], with the increase of $r_l$; see Figure \ref{fig: experiment 6a}. While robustness is achieved on average, i.e., the average relative in-sample loss is larger than~$1$, about~5$\%$ of the test instances for~$r_l = 10$ may violate robustness.
	\item Increasing $\tilde{\kappa}$ tends to provide more robust decisions at the cost of a slightly worse out-of-sample performance. 
	Indeed, larger uncertainty sets cause greater variability in the follower's data, which makes adverse follower data sets more plausible even when $r_l$ is small. Overall, 
	the additional out-of-sample error of~[\textbf{DRI-SP$'$}] relative to~[\textbf{DRI$^*$}] is uniformly bounded by 3--4$\%$ on~average. 
	\item The augmented basic model [\textbf{DRI}] shows no significant relationship with the true model~[\textbf{DRI$^*$}], resulting in overly optimistic in-sample and relatively poor out-of-sample performance.
\end{itemize}

Notably, for larger-scale problems, increasing the number of scenarios significantly raises the computational complexity of the MILP reformulation~(\ref{semi-pessimistic reformulation MILP}). While this may also enhance robustness, as illustrated in Figure~\ref{fig: experiment 6a}, we observe that, after some point, the marginal benefit of additional scenarios may not translate into substantial gains in the leader's out-of-sample performance. For example, based on Figure \ref{fig: experiment 6a}, we fix $r_l = 5$ in the remainder of the analysis.

\begin{figure}[h] \centering	\scalebox{0.8}{
	\begin{subfigure}{0.5\textwidth}
		\centering
		\begin{tikzpicture}
			\begin{axis}[
				axis lines=box,
				xlabel={$k_{lf}$},
				ylabel={Relative in-sample loss},
				xmin = 0, xmax = 30,
				ymin=0.9, ymax=1.7,
				grid=major,
				legend pos=north east,
				mark size=3pt,
				ytick = {0.9, 1, 1.1, 1.2, 1.3, 1.4, 1.5}
				]
				
				\addplot[color=black, mark=triangle*, mark size = 2pt] coordinates {
					(0,1.006) (2,1.006) (4,0.999) (6,1.001) (8,1.002) (10,0.999)
					(12,0.997) (14,0.997) (16,0.996) (18,0.998) (20,0.994) (22,1.003)
					(24,1.002) (26,1.000) (28,1.001) (30,1.000)
				};
				\addlegendentry{\footnotesize Augmented basic model }

				\addplot[color=forestgreen, mark=square*, mark size=1.5pt] coordinates { 
					(0,1.086) (2,1.079) (4,1.067) (6,1.067) (8,1.060) (10,1.051)
					(12,1.056) (14,1.050) (16,1.048) (18,1.042) (20,1.026) (22,1.031)
					(24,1.019) (26,1.014) (28,1.013) (30,1.000)
				};
				\addlegendentry{\footnotesize Scenario-based approximation ($\tilde{\kappa} = 0.5$)}

				\addplot[color=violet, mark=x, mark size=2pt] coordinates {
					(0,1.357) (2,1.357) (4,1.346) (6,1.349) (8,1.351) (10,1.348)
					(12,1.346) (14,1.346) (16,1.344) (18,1.347) (20,1.342) (22,1.355)
					(24,1.354) (26,1.351) (28,1.354) (30,1.352)
				};
				
				\addlegendentry{\footnotesize Pessimistic approximation}
				
				\addplot[mark=square*, mark size=1.5pt, dashed, color=forestgreen] coordinates {
					(0,1.000) (2,1.000) (4,0.985) (6,0.985) (8,0.985) (10,0.988)
					(12,0.980) (14,0.974) (16,0.978) (18,0.989) (20,0.978) (22,0.998)
					(24,0.995) (26,0.993) (28,1.000) (30,1.000)
				};
				
				\path[name path=upperA] plot coordinates {
					(0,1.006+0.056) (2,1.006+0.056) (4,0.999+0.058) (6,1.001+0.057) (8,1.002+0.054) (10,0.999+0.050)
					(12,0.997+0.052) (14,0.997+0.048) (16,0.996+0.043) (18,0.998+0.038) (20,0.994+0.037) (22,1.003+0.030)
					(24,1.002+0.025) (26,1.000+0.017) (28,1.001+0.011) (30,1.000+0.000)
				};
				
				\path[name path=lowerA] plot coordinates {
					(0,1.006-0.056) (2,1.006-0.056) (4,0.999-0.058) (6,1.001-0.057) (8,1.002-0.054) (10,0.999-0.050)
					(12,0.997-0.052) (14,0.997-0.048) (16,0.996-0.043) (18,0.998-0.038) (20,0.994-0.037) (22,1.003-0.030)
					(24,1.002-0.025) (26,1.000-0.017) (28,1.001-0.011) (30,1.000-0.000)
				};
				\addplot[black!10, fill opacity=0.5] fill between[of=upperA and lowerA];
				
				\addplot[name path=upperB, draw = none] coordinates {
					(0,1.086+0.056) (2,1.079+0.057) (4,1.067+0.049) (6,1.067+0.052) (8,1.060+0.049) (10,1.051+0.049)
					(12,1.056+0.053) (14,1.050+0.051) (16,1.048+0.048) (18,1.042+0.047) (20,1.026+0.034) (22,1.031+0.035)
					(24,1.019+0.024) (26,1.014+0.020) (28,1.013+0.021) (30,1.000+0.000)
				};
				
				\addplot[name path=lowerB, draw = none] coordinates {
					(0,1.086-0.056) (2,1.079-0.057) (4,1.067-0.049) (6,1.067-0.052) (8,1.060-0.049) (10,1.051-0.049)
					(12,1.056-0.053) (14,1.050-0.051) (16,1.048-0.048) (18,1.042-0.047) (20,1.026-0.034) (22,1.031-0.035)
					(24,1.019-0.024) (26,1.014-0.020) (28,1.013-0.021) (30,1.000-0.000)
				};
				
				\addplot[forestgreen!20, fill opacity=0.5] fill between[of=upperB and lowerB];
				
				\path[name path=upperC] plot coordinates {
					(0,1.357+0.090) (2,1.357+0.085) (4,1.346+0.081) (6,1.349+0.084) (8,1.351+0.088) (10,1.348+0.088)
					(12,1.346+0.090) (14,1.346+0.089) (16,1.344+0.091) (18,1.347+0.090) (20,1.342+0.087) (22,1.355+0.088)
					(24,1.354+0.087) (26,1.351+0.087) (28,1.354+0.091) (30,1.352+0.094)
				};
				
				\path[name path=lowerC] plot coordinates {
					(0,1.357-0.090) (2,1.357-0.085) (4,1.346-0.081) (6,1.349-0.084) (8,1.351-0.088) (10,1.348-0.088)
					(12,1.346-0.090) (14,1.346-0.089) (16,1.344-0.091) (18,1.347-0.090) (20,1.342-0.087) (22,1.355-0.088)
					(24,1.354-0.087) (26,1.351-0.087) (28,1.354-0.091) (30,1.352-0.094)
				};
				
				\addplot[violet!10, fill opacity=0.5] fill between[of=upperC and lowerC];
				\addplot[black, thin, gray, domain=0:30, samples=2] {1};
			\end{axis}
		\end{tikzpicture}
		\caption{}
		\label{fig: experiment 7a}
	\end{subfigure}
	\hfill
	\begin{subfigure}{0.5\textwidth}
		\centering 
		\begin{tikzpicture}
			\begin{axis}[
				axis lines=box,
				xlabel={$k_{lf}$},
				ylabel={Relative out-of-sample loss},
				xmin = 0,
				xmax = 30,
				ymin=0.95, ymax=1.4,
				grid=major,
				legend pos=north east,
				mark size=3pt, 	ytick = {1, 1.05, 1.1, 1.15, 1.2, 1.25}
				]
				
				\addplot[color=black, mark=triangle*, mark size=2pt, dashed] coordinates { 
				 (0,1.022) (2,1.021) (4,1.025) (6,1.026) (8,1.028) (10,1.031)
				 (12,1.032) (14,1.026) (16,1.025) (18,1.025) (20,1.027) (22,1.026)
				 (24,1.028) (26,1.028) (28,1.028) (30,1.026)
				};
				
				\addlegendentry{\footnotesize True basic model}

				\addplot[color=black, mark=triangle*, mark size = 2pt] coordinates {
					(0,1.116) (2,1.115) (4,1.115) (6,1.111) (8,1.111) (10,1.112)
					(12,1.105) (14,1.094) (16,1.090) (18,1.081) (20,1.074) (22,1.069)
					(24,1.052) (26,1.045) (28,1.033) (30,1.026)
				};
				\addlegendentry{\footnotesize Augmented basic model}
				
				\addplot[color=forestgreen, mark=square*, mark size=1.5pt] coordinates { 
					(0,1.061) (2,1.071) (4,1.082) (6,1.068) (8,1.068) (10,1.055)
					(12,1.060) (14,1.052) (16,1.054) (18,1.050) (20,1.047) (22,1.040)
					(24,1.048) (26,1.040) (28,1.038) (30,1.026)
				};
			
				\addlegendentry{\footnotesize Scenario-based approximation ($\tilde{\kappa} = 0.5$)}	
				
				\addplot[color=violet, mark=x, mark size=2pt] coordinates {
					(0,1.127) (2,1.128) (4,1.124) (6,1.124) (8,1.125) (10,1.124)
					(12,1.124) (14,1.126) (16,1.125) (18,1.133) (20,1.135) (22,1.138)
					(24,1.136) (26,1.135) (28,1.136) (30,1.131)
				};
				\addlegendentry{\footnotesize Pessimistic approximation}
				
				\path[name path=upperA] plot coordinates {
					(0,1.116+0.100) (2,1.115+0.096) (4,1.115+0.091) (6,1.111+0.089) (8,1.111+0.090) (10,1.112+0.093)
					(12,1.105+0.087) (14,1.094+0.078) (16,1.090+0.080) (18,1.081+0.073) (20,1.074+0.066) (22,1.069+0.064)
					(24,1.052+0.050) (26,1.045+0.047) (28,1.033+0.034) (30,1.026+0.031)
				};
				\path[name path=lowerA] plot coordinates {
					(0,1.116-0.100) (2,1.115-0.096) (4,1.115-0.091) (6,1.111-0.089) (8,1.111-0.090) (10,1.112-0.093)
					(12,1.105-0.087) (14,1.094-0.078) (16,1.090-0.080) (18,1.081-0.073) (20,1.074-0.066) (22,1.069-0.064)
					(24,1.052-0.050) (26,1.045-0.047) (28,1.033-0.034) (30,1.026-0.031)
				};
				\addplot[black!10, fill opacity=0.5] fill between[of=upperA and lowerA];
				
				\path[name path=upperC] plot coordinates {
					(0,1.022+0.027) (2,1.021+0.025) (4,1.025+0.028) (6,1.026+0.029) (8,1.028+0.030) (10,1.031+0.032)
					(12,1.032+0.032) (14,1.026+0.029) (16,1.025+0.030) (18,1.025+0.028) (20,1.027+0.032) (22,1.026+0.032)
					(24,1.028+0.032) (26,1.028+0.032) (28,1.028+0.032) (30,1.026+0.031)
				};
				
				\path[name path=lowerC] plot coordinates {
					(0,1.022-0.027) (2,1.021-0.025) (4,1.025-0.028) (6,1.026-0.029) (8,1.028-0.030) (10,1.031-0.032)
					(12,1.032-0.032) (14,1.026-0.029) (16,1.025-0.030) (18,1.025-0.028) (20,1.027-0.032) (22,1.026-0.032)
					(24,1.028-0.032) (26,1.028-0.032) (28,1.028-0.032) (30,1.026-0.031)
				};
				
				\addplot[black!30, fill opacity=0.5] fill between[of=upperC and lowerC];	
				
				\path[name path=upperD] plot coordinates {
					(0,1.127+0.095) (2,1.128+0.095) (4,1.124+0.089) (6,1.124+0.086) (8,1.125+0.085) (10,1.124+0.086)
					(12,1.124+0.090) (14,1.126+0.091) (16,1.125+0.096) (18,1.133+0.095) (20,1.135+0.094) (22,1.138+0.102)
					(24,1.136+0.100) (26,1.135+0.098) (28,1.136+0.100) (30,1.131+0.106)
				};
				
				\path[name path=lowerD] plot coordinates {
				 (0,1.127-0.095) (2,1.128-0.095) (4,1.124-0.089) (6,1.124-0.086) (8,1.125-0.085) (10,1.124-0.086)
				 (12,1.124-0.090) (14,1.126-0.091) (16,1.125-0.096) (18,1.133-0.095) (20,1.135-0.094) (22,1.138-0.102)
				 (24,1.136-0.100) (26,1.135-0.098) (28,1.136-0.100) (30,1.131-0.106)
				};
			
				\addplot[violet!10, fill opacity=0.5] fill between[of=upperD and lowerD];
				
				\path[name path=upperB] plot coordinates {
					(0,1.061+0.058) (2,1.071+0.061) (4,1.082+0.076) (6,1.068+0.057) (8,1.068+0.056) (10,1.055+0.049)
					(12,1.060+0.051) (14,1.052+0.048) (16,1.054+0.056) (18,1.050+0.052) (20,1.047+0.049) (22,1.040+0.043)
					(24,1.048+0.049) (26,1.040+0.040) (28,1.038+0.038) (30,1.026+0.031)
				};
				
				\path[name path=lowerB] plot coordinates {
					(0,1.061-0.058) (2,1.071-0.061) (4,1.082-0.076) (6,1.068-0.057) (8,1.068-0.056) (10,1.055-0.049)
					(12,1.060-0.051) (14,1.052-0.048) (16,1.054-0.056) (18,1.050-0.052) (20,1.047-0.049) (22,1.040-0.043)
					(24,1.048-0.049) (26,1.040-0.040) (28,1.038-0.038) (30,1.026-0.031)
				};
				
				\addplot[forestgreen!20, fill opacity=0.5] fill between[of=upperB and lowerB];
				\addplot[black, thin, gray, domain=0:30, samples=2] {1};
			\end{axis}
		\end{tikzpicture}
		\caption{}
		\label{fig: experiment 7b}
	\end{subfigure}}
	\caption{The average relative in-sample (\ref{eq: relative loss leader pessimistic}) (a) and out-of-sample (\ref{eq: relative in-sample loss}) (b) loss of the \textit{ambiguity-free} leader (with MADs) as a function of the common sample size, $k_{lf}$, for $k_l = k_f = 30$, $\delta_l = \delta_f = 0.1$ and $\alpha_l = 0.9$, evaluated over 100 random test instances. The follower is assumed to be risk-averse. The dashed line corresponds to the empirical $5\%$ percentile of the relative in-sample loss for the scenario-based semi-pessimistic approximation.}
	\label{fig: experiment 7}
\end{figure}
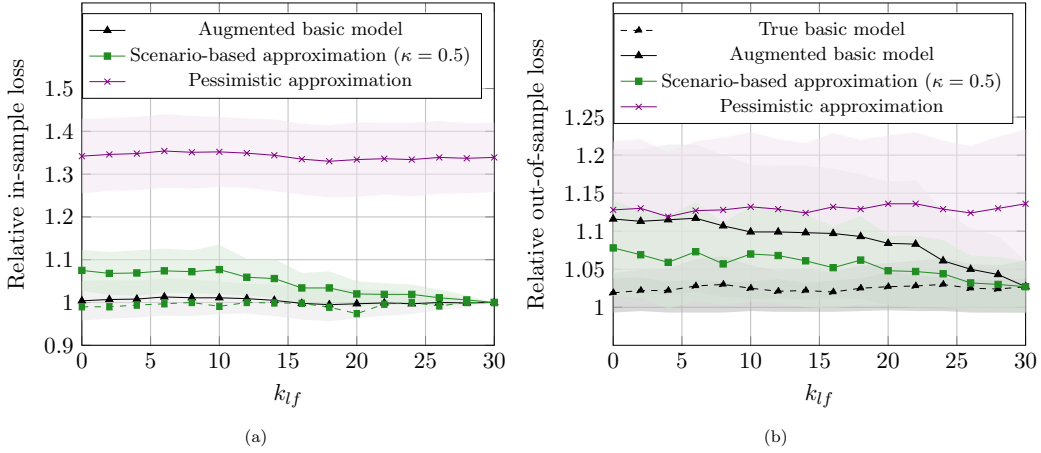

 In the following, we compare [\textbf{DRI-SP$'$}] with the pessimistic approximation [\textbf{DRI-P}] and the augmented basic model [\textbf{DRI}] for the ambiguity-free leader ($\varepsilon_l = 0$); see Figure \ref{fig: experiment 7}. Specifically, the common number of samples, $k_{lf}$, is increased from~$0$ to $30$, with the true follower's data set being~fixed and $k_l = k_f = 30$. 
Finally, we set $\alpha_l = 0.9$ to satisfy Assumption \textbf{A6}. 

From Figure \ref{fig: experiment 7}, we observe that the relative in-sample loss of the leader in [\textbf{DRI-SP$'$}] and [\textbf{DRI}] converges to 1 as the leader acquires more data from the follower. A similar trend is observed in terms of the relative out-of-sample performance, i.e., both models converge to the true basic model~[\textbf{DRI$^*$}]. Importantly, both [\textbf{DRI}] and~[\textbf{DRI-P}] exhibit poor out-of-sample performance. Even when the leader is entirely uncertain about the follower's data ($k_{lf} = 0$), meaning that approximately half of all data entries are subject to severe noise ($\tilde{\kappa} = 0.5$), the scenario-based approximation [\textbf{DRI-SP$'$}] still significantly outperforms both outlined models. It can be explained by the fact that the augmented basic model [\textbf{DRI}] substantially alters the follower's data set. At the same time, the pessimistic approximation [\textbf{DRI-P}] not only disregards the follower's data but also overlooks the structure of the follower's objective function.

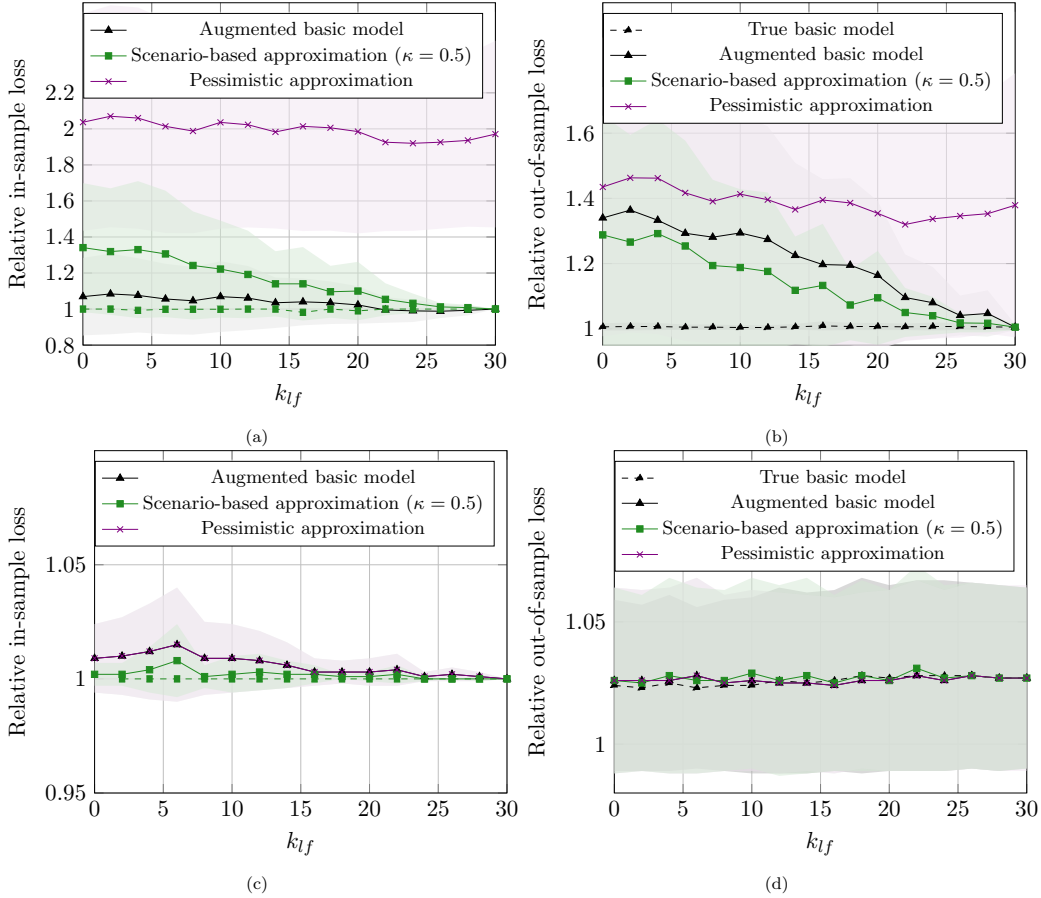
\begin{figure}[h] 
	 \begin{minipage}{\textwidth} \centering
	 \scalebox{0.8}{\begin{subfigure}{0.5\textwidth}
		\centering
		\begin{tikzpicture}
			\begin{axis}[
				axis lines=box,
				xlabel={$k_{lf}$},
				ylabel={Relative in-sample loss},
				xmin = 0, xmax = 30,
				ymin=0.8, ymax=2.2,
				grid=major,
				legend pos=north east,
				mark size=3pt,
				ytick = {0.8, 1,1.2, 1.4, 1.6,1.8, 2}
				]
				
				\addplot[color=black, mark=triangle*, mark size = 2pt] coordinates {
					(0,1.029) (2,1.036) (4,1.038) (6,1.029) (8,1.030) (10,1.028)
					(12,1.006) (14,1.003) (16,1.014) (18,1.017) (20,0.995) (22,0.980)
					(24,0.984) (26,0.989) (28,0.997) (30,1.000)
				};
				\addlegendentry{\footnotesize Augmented basic model}

				\addplot[color=forestgreen, mark=square*, mark size=1.5pt] coordinates { 
					(0,1.169) (2,1.181) (4,1.173) (6,1.164) (8,1.158) (10,1.119)
					(12,1.121) (14,1.099) (16,1.119) (18,1.087) (20,1.058) (22,1.042)
					(24,1.020) (26,1.013) (28,1.009) (30,1.000)
				};
				\addlegendentry{\footnotesize Scenario-based approximation ($\tilde{\kappa} = 0.5$)}

				\addplot[color=violet, mark=x, mark size=2pt] coordinates {
					(0,1.683) (2,1.699) (4,1.715) (6,1.699) (8,1.703) (10,1.710)
					(12,1.683) (14,1.670) (16,1.699) (18,1.710) (20,1.671) (22,1.655)
					(24,1.670) (26,1.678) (28,1.708) (30,1.724)
				};
				
				\addlegendentry{\footnotesize Pessimistic approximation}
				
				\addplot[mark=square*, mark size=1.5pt, dashed, color=forestgreen] coordinates {
					(0,1.000) (2,1.000) (4,0.985) (6,1.000) (8,0.990) (10,1.000)
					(12,0.995) (14,1.000) (16,0.994) (18,0.988) (20,1.000) (22,0.989)
					(24,1.000) (26,1.000) (28,1.000) (30,1.000)
				};
				
				\path[name path=upperA] plot coordinates {
					(0,1.029+0.198) (2,1.036+0.198) (4,1.038+0.183) (6,1.029+0.172) (8,1.030+0.169) (10,1.028+0.152)
					(12,1.006+0.116) (14,1.003+0.123) (16,1.014+0.115) (18,1.017+0.109) (20,0.995+0.092) (22,0.980+0.075)
					(24,0.984+0.059) (26,0.989+0.055) (28,0.997+0.026) (30,1.000+0.000)
				};
				\path[name path=lowerA] plot coordinates {
					(0,1.029-0.198) (2,1.036-0.198) (4,1.038-0.183) (6,1.029-0.172) (8,1.030-0.169) (10,1.028-0.152)
					(12,1.006-0.116) (14,1.003-0.123) (16,1.014-0.115) (18,1.017-0.109) (20,0.995-0.092) (22,0.980-0.075)
					(24,0.984-0.059) (26,0.989-0.055) (28,0.997-0.026) (30,1.000-0.000)
				};
				\addplot[black!10, fill opacity=0.5] fill between[of=upperA and lowerA];
				
				\path[name path=upperC] plot coordinates {
					(0,1.683+0.295) (2,1.699+0.315) (4,1.715+0.338) (6,1.699+0.318) (8,1.703+0.328) (10,1.710+0.339)
					(12,1.683+0.313) (14,1.670+0.282) (16,1.699+0.320) (18,1.710+0.340) (20,1.671+0.295) (22,1.655+0.274)
					(24,1.670+0.296) (26,1.678+0.306) (28,1.708+0.345) (30,1.724+0.366)
				};
				
				\path[name path=lowerC] plot coordinates {
					(0,1.683-0.295) (2,1.699-0.315) (4,1.715-0.338) (6,1.699-0.318) (8,1.703-0.328) (10,1.710-0.339)
					(12,1.683-0.313) (14,1.670-0.282) (16,1.699-0.320) (18,1.710-0.340) (20,1.671-0.295) (22,1.655-0.274)
					(24,1.670-0.296) (26,1.678-0.306) (28,1.708-0.345) (30,1.724-0.366)
				};
				
				\addplot[violet!10, fill opacity=0.5] fill between[of=upperC and lowerC];

				\addplot[name path=upperB, draw = none] coordinates {
					(0,1.169+0.182) (2,1.181+0.197) (4,1.173+0.209) (6,1.164+0.190) (8,1.158+0.196) (10,1.119+0.137)
					(12,1.121+0.147) (14,1.099+0.130) (16,1.119+0.147) (18,1.087+0.122) (20,1.058+0.084) (22,1.042+0.064)
					(24,1.020+0.032) (26,1.013+0.022) (28,1.009+0.017) (30,1.000+0.000)
				};
				
				\addplot[name path=lowerB, draw = none] coordinates {
					(0,1.169-0.182) (2,1.181-0.197) (4,1.173-0.209) (6,1.164-0.190) (8,1.158-0.196) (10,1.119-0.137)
					(12,1.121-0.147) (14,1.099-0.130) (16,1.119-0.147) (18,1.087-0.122) (20,1.058-0.084) (22,1.042-0.064)
					(24,1.020-0.032) (26,1.013-0.022) (28,1.009-0.017) (30,1.000-0.000)
				};
				
				\addplot[forestgreen!20, fill opacity=0.5] fill between[of=upperB and lowerB];
			 \addplot[black, thin, gray, domain=0:30, samples=2] {1};
			\end{axis}
		\end{tikzpicture}
		\caption{}
		\label{fig: experiment 8a}
	\end{subfigure}
	\hfill
	\begin{subfigure}{0.5\textwidth}
		\centering 
		\begin{tikzpicture}
			\begin{axis}[
				axis lines=box,
				xlabel={$k_{lf}$},
				ylabel={Relative out-of-sample loss},
				xmin = 0,
				xmax = 30,
				ymin=0.95, ymax=1.7,
				grid=major,
				legend pos=north east,
				mark size=3pt, 	ytick = {1, 1.1, 1.2, 1.3}
				]
				
				\addplot[color=black, mark=triangle*, mark size=2pt, dashed] coordinates { 
					(0,1.007) (2,1.008) (4,1.007) (6,1.007) (8,1.006) (10,1.008)
					(12,1.006) (14,1.006) (16,1.004) (18,1.004) (20,1.007) (22,1.008)
					(24,1.009) (26,1.008) (28,1.007) (30,1.004)
				};
				
				\addlegendentry{\footnotesize True basic model}

				\addplot[color=black, mark=triangle*, mark size = 2pt] coordinates {
					(0,1.243) (2,1.245) (4,1.246) (6,1.236) (8,1.211) (10,1.178)
					(12,1.147) (14,1.140) (16,1.144) (18,1.146) (20,1.110) (22,1.073)
					(24,1.058) (26,1.040) (28,1.027) (30,1.004)
				};
				\addlegendentry{\footnotesize Augmented basic model}
				
				\addplot[color=forestgreen, mark=square*, mark size=1.5pt] coordinates { 
					(0,1.191) (2,1.183) (4,1.152) (6,1.157) (8,1.146) (10,1.116)
					(12,1.101) (14,1.102) (16,1.116) (18,1.094) (20,1.070) (22,1.047)
					(24,1.026) (26,1.022) (28,1.018) (30,1.004)
				};
				
				\addlegendentry{\footnotesize Scenario-based approximation ($\tilde{\kappa} = 0.5$)}	
				
				\addplot[color=violet, mark=x, mark size=2pt] coordinates {
					(0,1.223) (2,1.226) (4,1.232) (6,1.224) (8,1.220) (10,1.229)
					(12,1.218) (14,1.203) (16,1.229) (18,1.230) (20,1.200) (22,1.183)
					(24,1.205) (26,1.209) (28,1.234) (30,1.241)
				};
				\addlegendentry{\footnotesize Pessimistic approximation}
				
				\path[name path=upperA] plot coordinates {
					(0,1.243+0.254) (2,1.245+0.257) (4,1.246+0.260) (6,1.236+0.249) (8,1.211+0.233) (10,1.178+0.207)
					(12,1.147+0.173) (14,1.140+0.159) (16,1.144+0.175) (18,1.146+0.191) (20,1.110+0.145) (22,1.073+0.095)
					(24,1.058+0.076) (26,1.040+0.057) (28,1.027+0.040) (30,1.004+0.007)
				};
				\path[name path=lowerA] plot coordinates {
					(0,1.243-0.254) (2,1.245-0.257) (4,1.246-0.260) (6,1.236-0.249) (8,1.211-0.233) (10,1.178-0.207)
					(12,1.147-0.173) (14,1.140-0.159) (16,1.144-0.175) (18,1.146-0.191) (20,1.110-0.145) (22,1.073-0.095)
					(24,1.058-0.076) (26,1.040-0.057) (28,1.027-0.040) (30,1.004-0.007)
				};
				\addplot[black!10, fill opacity=0.5] fill between[of=upperA and lowerA];
				
				\path[name path=upperC] plot coordinates {
					(0,1.007+0.011) (2,1.008+0.013) (4,1.007+0.011) (6,1.007+0.012) (8,1.006+0.010) (10,1.008+0.011)
					(12,1.006+0.009) (14,1.006+0.009) (16,1.004+0.007) (18,1.004+0.006) (20,1.007+0.010) (22,1.008+0.012)
					(24,1.009+0.014) (26,1.008+0.012) (28,1.007+0.010) (30,1.004+0.007)
				};
				
				\path[name path=lowerC] plot coordinates {
					(0,1.007-0.011) (2,1.008-0.013) (4,1.007-0.011) (6,1.007-0.012) (8,1.006-0.010) (10,1.008-0.011)
					(12,1.006-0.009) (14,1.006-0.009) (16,1.004-0.007) (18,1.004-0.006) (20,1.007-0.010) (22,1.008-0.012)
					(24,1.009-0.014) (26,1.008-0.012) (28,1.007-0.010) (30,1.004-0.007)
				};
				
				\addplot[black!30, fill opacity=0.5] fill between[of=upperC and lowerC];	
				
				\path[name path=upperD] plot coordinates {
					(0,1.223+0.243) (2,1.226+0.247) (4,1.232+0.255) (6,1.224+0.241) (8,1.220+0.236) (10,1.229+0.254)
					(12,1.218+0.241) (14,1.203+0.219) (16,1.229+0.257) (18,1.230+0.267) (20,1.200+0.227) (22,1.183+0.211)
					(24,1.205+0.250) (26,1.209+0.249) (28,1.234+0.284) (30,1.241+0.295)
				};
				
				\path[name path=lowerD] plot coordinates {
				(0,1.223-0.243) (2,1.226-0.247) (4,1.232-0.255) (6,1.224-0.241) (8,1.220-0.236) (10,1.229-0.254)
				(12,1.218-0.241) (14,1.203-0.219) (16,1.229-0.257) (18,1.230-0.267) (20,1.200-0.227) (22,1.183-0.211)
				(24,1.205-0.250) (26,1.209-0.249) (28,1.234-0.284) (30,1.241-0.295)
				};
				
				\addplot[violet!10, fill opacity=0.5] fill between[of=upperD and lowerD];
				
				\path[name path=upperB] plot coordinates {
					(0,1.191+0.216) (2,1.183+0.222) (4,1.152+0.194) (6,1.157+0.192) (8,1.146+0.180) (10,1.116+0.148)
					(12,1.101+0.134) (14,1.102+0.130) (16,1.116+0.160) (18,1.094+0.125) (20,1.070+0.101) (22,1.047+0.064)
					(24,1.026+0.038) (26,1.022+0.033) (28,1.018+0.029) (30,1.004+0.007)
				};
				
				\path[name path=lowerB] plot coordinates {
					(0,1.191-0.216) (2,1.183-0.222) (4,1.152-0.194) (6,1.157-0.192) (8,1.146-0.180) (10,1.116-0.148)
					(12,1.101-0.134) (14,1.102-0.130) (16,1.116-0.160) (18,1.094-0.125) (20,1.070-0.101) (22,1.047-0.064)
					(24,1.026-0.038) (26,1.022-0.033) (28,1.018-0.029) (30,1.004-0.007)
				};
				
				\addplot[forestgreen!20, fill opacity=0.5] fill between[of=upperB and lowerB];
				\addplot[black, thin, gray, domain=0:30, samples=2] {1};
			\end{axis}
		\end{tikzpicture}
		\caption{}
		\label{fig: experiment 8b}
	\end{subfigure}}
 \end{minipage}
 \newline
 \begin{minipage}{\textwidth} \centering
	\scalebox{0.8}{\begin{subfigure}{0.5\textwidth}
		\centering
		\begin{tikzpicture}
			\begin{axis}[
				axis lines=box,
				xlabel={$k_{lf}$},
				ylabel={Relative in-sample loss},
				xmin = 0, xmax = 30,
				ymin=0.95, ymax=1.1,
				grid=major,
				legend pos=north east,
				mark size=3pt,
				ytick = {0.95, 1, 1.05}
				]
				
				\addplot[color=black, mark=triangle*, mark size = 2pt] coordinates {
					(0,1.020) (2,1.021) (4,1.016) (6,1.014) (8,1.014) (10,1.010)
					(12,1.006) (14,1.007) (16,1.007) (18,1.007) (20,1.004) (22,1.003)
					(24,1.004) (26,1.000) (28,1.001) (30,1.000)
				};
				\addlegendentry{\footnotesize Augmented basic model}

				\addplot[color=forestgreen, mark=square*, mark size=1.5pt] coordinates { 
					(0,1.010) (2,1.010) (4,1.006) (6,1.005) (8,1.005) (10,1.007)
					(12,1.001) (14,1.003) (16,1.003) (18,1.005) (20,1.001) (22,1.001)
					(24,1.001) (26,1.000) (28,1.001) (30,1.000)
				};
				\addlegendentry{\footnotesize Scenario-based approximation ($\tilde{\kappa} = 0.5$)}

				\addplot[color=violet, mark=x, mark size=2pt] coordinates {
					(0,1.020) (2,1.021) (4,1.016) (6,1.014) (8,1.014) (10,1.011)
					(12,1.007) (14,1.007) (16,1.007) (18,1.007) (20,1.004) (22,1.004)
					(24,1.004) (26,1.001) (28,1.001) (30,1.000)
				};
				
				\addlegendentry{\footnotesize Pessimistic approximation}
				
				\addplot[mark=square*, mark size=1.5pt, dashed, color=forestgreen] coordinates {
					(0,1.000) (2,1.000) (4,1.000) (6,1.000) (8,1.000) (10,1.000)
					(12,1.000) (14,1.000) (16,1.000) (18,1.000) (20,1.000) (22,1.000)
					(24,1.000) (26,1.000) (28,1.000) (30,1.000)
				};
				
				\path[name path=upperA] plot coordinates {
					(0,1.020+0.030) (2,1.021+0.032) (4,1.016+0.025) (6,1.014+0.023) (8,1.014+0.021) (10,1.010+0.016)
					(12,1.006+0.011) (14,1.007+0.012) (16,1.007+0.012) (18,1.007+0.012) (20,1.004+0.007) (22,1.003+0.006)
					(24,1.004+0.007) (26,1.000+0.001) (28,1.001+0.002) (30,1.000+0.000)
				};
				\path[name path=lowerA] plot coordinates {
					(0,1.020-0.030) (2,1.021-0.032) (4,1.016-0.025) (6,1.014-0.023) (8,1.014-0.021) (10,1.010-0.016)
					(12,1.006-0.011) (14,1.007-0.012) (16,1.007-0.012) (18,1.007-0.012) (20,1.004-0.007) (22,1.003-0.006)
					(24,1.004-0.007) (26,1.000-0.001) (28,1.001-0.002) (30,1.000-0.000)
				};
				\addplot[black!10, fill opacity=0.5] fill between[of=upperA and lowerA];
				
				\path[name path=upperC] plot coordinates {
				 (0,1.020+0.030) (2,1.021+0.031) (4,1.016+0.025) (6,1.014+0.023) (8,1.014+0.021) (10,1.011+0.016)
				 (12,1.007+0.011) (14,1.007+0.012) (16,1.007+0.012) (18,1.007+0.012) (20,1.004+0.007) (22,1.004+0.006)
				 (24,1.004+0.007) (26,1.001+0.001) (28,1.001+0.002) (30,1.000+0.000)
				};
				
				\path[name path=lowerC] plot coordinates {
					(0,1.020-0.030) (2,1.021-0.031) (4,1.016-0.025) (6,1.014-0.023) (8,1.014-0.021) (10,1.011-0.016)
					(12,1.007-0.011) (14,1.007-0.012) (16,1.007-0.012) (18,1.007-0.012) (20,1.004-0.007) (22,1.004-0.006)
					(24,1.004-0.007) (26,1.001-0.001) (28,1.001-0.002) (30,1.000-0.000)
				};
				
				\addplot[violet!10, fill opacity=0.5] fill between[of=upperC and lowerC];
				
				\addplot[name path=upperB, draw = none] coordinates {
					(0,1.010+0.019) (2,1.010+0.019) (4,1.006+0.013) (6,1.005+0.013) (8,1.005+0.009) (10,1.007+0.011)
					(12,1.001+0.006) (14,1.003+0.006) (16,1.003+0.006) (18,1.005+0.009) (20,1.001+0.002) (22,1.001+0.003)
					(24,1.001+0.003) (26,1.000+0.000) (28,1.001+0.001) (30,1.000+0.000)
				};
				
				\addplot[name path=lowerB, draw = none] coordinates {
				(0,1.010+0.019) (2,1.010+0.019) (4,1.006+0.013) (6,1.005+0.013) (8,1.005+0.009) (10,1.007+0.011)
				(12,1.001+0.006) (14,1.003+0.006) (16,1.003+0.006) (18,1.005+0.009) (20,1.001+0.002) (22,1.001+0.003)
				(24,1.001+0.003) (26,1.000+0.000) (28,1.001+0.001) (30,1.000+0.000)
				};
				
				\addplot[forestgreen!20, fill opacity=0.5] fill between[of=upperB and lowerB];
				\addplot[black, thin, gray, domain=0:30, samples=2] {1};
				
			\end{axis}
		\end{tikzpicture}
		\caption{}
		\label{fig: experiment 8c}
	\end{subfigure}
	\hfill
	\begin{subfigure}{0.5\textwidth}
		\centering 
		\begin{tikzpicture}
			\begin{axis}[
				axis lines=box,
				xlabel={$k_{lf}$},
				ylabel={Relative out-of-sample loss},
				xmin = 0,
				xmax = 30,
				ymin=0.98, ymax=1.12,
				grid=major,
				legend pos=north east,
				mark size=3pt, 	ytick = {1, 1.05}
				]
				
				\addplot[color=black, mark=triangle*, mark size=2pt, dashed] coordinates { 
					(0,1.026) (2,1.026) (4,1.026) (6,1.025) (8,1.026) (10,1.024)
					(12,1.025) (14,1.025) (16,1.024) (18,1.023) (20,1.022) (22,1.023)
					(24,1.022) (26,1.023) (28,1.023) (30,1.024)
				};
				
				\addlegendentry{\footnotesize True basic model}

				\addplot[color=black, mark=triangle*, mark size = 2pt] coordinates {
					(0,1.028) (2,1.027) (4,1.028) (6,1.029) (8,1.028) (10,1.027)
					(12,1.023) (14,1.025) (16,1.023) (18,1.025) (20,1.024) (22,1.024)
					(24,1.023) (26,1.023) (28,1.023) (30,1.024)
				};
				\addlegendentry{\footnotesize Augmented basic model}
				
				\addplot[color=forestgreen, mark=square*, mark size=1.5pt] coordinates { 
					(0,1.031) (2,1.029) (4,1.027) (6,1.025) (8,1.029) (10,1.024)
					(12,1.021) (14,1.026) (16,1.024) (18,1.024) (20,1.023) (22,1.023)
					(24,1.023) (26,1.023) (28,1.023) (30,1.024)
				};
				
				\addlegendentry{\footnotesize Scenario-based approximation ($\tilde{\kappa} = 0.5$)}	
				
				\addplot[color=violet, mark=x, mark size=2pt] coordinates {
					(0,1.027) (2,1.026) (4,1.027) (6,1.028) (8,1.027) (10,1.026)
					(12,1.023) (14,1.025) (16,1.023) (18,1.025) (20,1.024) (22,1.024)
					(24,1.023) (26,1.023) (28,1.023) (30,1.024)
				};
				\addlegendentry{\footnotesize Pessimistic approximation}
				
				\path[name path=upperA] plot coordinates {
					(0,1.028+0.034) (2,1.027+0.033) (4,1.028+0.033) (6,1.029+0.034) (8,1.028+0.033) (10,1.027+0.032)
					(12,1.023+0.029) (14,1.025+0.030) (16,1.023+0.029) (18,1.025+0.030) (20,1.024+0.030) (22,1.024+0.030)
					(24,1.023+0.030) (26,1.023+0.030) (28,1.023+0.030) (30,1.024+0.030)
				};
				\path[name path=lowerA] plot coordinates {
					(0,1.028-0.034) (2,1.027-0.033) (4,1.028-0.033) (6,1.029-0.034) (8,1.028-0.033) (10,1.027-0.032)
					(12,1.023-0.029) (14,1.025-0.030) (16,1.023-0.029) (18,1.025-0.030) (20,1.024-0.030) (22,1.024-0.030)
					(24,1.023-0.030) (26,1.023-0.030) (28,1.023-0.030) (30,1.024-0.030)
				};
				\addplot[black!10, fill opacity=0.5] fill between[of=upperA and lowerA];
				
				\path[name path=upperC] plot coordinates {
					(0,1.026+0.032) (2,1.026+0.032) (4,1.026+0.032) (6,1.025+0.031) (8,1.026+0.032) (10,1.024+0.031)
					(12,1.025+0.031) (14,1.025+0.031) (16,1.024+0.030) (18,1.023+0.029) (20,1.022+0.029) (22,1.023+0.029)
					(24,1.022+0.029) (26,1.023+0.030) (28,1.023+0.030) (30,1.024+0.030)
				};
				
				\path[name path=lowerC] plot coordinates {
				(0,1.026-0.032) (2,1.026-0.032) (4,1.026-0.032) (6,1.025-0.031) (8,1.026-0.032) (10,1.024-0.031)
				(12,1.025-0.031) (14,1.025-0.031) (16,1.024-0.030) (18,1.023-0.029) (20,1.022-0.029) (22,1.023-0.029)
				(24,1.022-0.029) (26,1.023-0.030) (28,1.023-0.030) (30,1.024-0.030)
				};
				
				\addplot[black!30, fill opacity=0.5] fill between[of=upperC and lowerC];	
				
				\path[name path=upperD] plot coordinates {
					(0,1.027+0.033) (2,1.026+0.032) (4,1.027+0.032) (6,1.028+0.033) (8,1.027+0.032) (10,1.026+0.032)
					(12,1.023+0.029) (14,1.025+0.030) (16,1.023+0.029) (18,1.025+0.030) (20,1.024+0.030) (22,1.024+0.030)
					(24,1.023+0.030) (26,1.023+0.030) (28,1.023+0.030) (30,1.024+0.030)
				};
				
				\path[name path=lowerD] plot coordinates {
					(0,1.027-0.033) (2,1.026-0.032) (4,1.027-0.032) (6,1.028-0.033) (8,1.027-0.032) (10,1.026-0.032)
					(12,1.023-0.029) (14,1.025-0.030) (16,1.023-0.029) (18,1.025-0.030) (20,1.024-0.030) (22,1.024-0.030)
					(24,1.023-0.030) (26,1.023-0.030) (28,1.023-0.030) (30,1.024-0.030)
				};
				
				\addplot[violet!10, fill opacity=0.5] fill between[of=upperD and lowerD];
				
				\path[name path=upperB] plot coordinates {
					(0,1.031+0.037) (2,1.029+0.036) (4,1.027+0.033) (6,1.025+0.030) (8,1.029+0.034) (10,1.024+0.031)
					(12,1.021+0.028) (14,1.026+0.031) (16,1.024+0.031) (18,1.024+0.030) (20,1.023+0.029) (22,1.023+0.029)
					(24,1.023+0.029) (26,1.023+0.030) (28,1.023+0.030) (30,1.024+0.030)
				};
				
				\path[name path=lowerB] plot coordinates {
			  (0,1.031-0.037) (2,1.029-0.036) (4,1.027-0.033) (6,1.025-0.030) (8,1.029-0.034) (10,1.024-0.031)
			  (12,1.021-0.028) (14,1.026-0.031) (16,1.024-0.031) (18,1.024-0.030) (20,1.023-0.029) (22,1.023-0.029)
			  (24,1.023-0.029) (26,1.023-0.030) (28,1.023-0.030) (30,1.024-0.030)
				};
				
				\addplot[forestgreen!20, fill opacity=0.5] fill between[of=upperB and lowerB];
				\addplot[black, thin, gray, domain=0:30, samples=2] {1};
			\end{axis}
		\end{tikzpicture}
		\caption{}
		\label{fig: experiment 8d}
	\end{subfigure}}
 \end{minipage}
	\caption{The average relative in-sample (\ref{eq: relative loss leader pessimistic}) (a, c) and out-of-sample (\ref{eq: relative in-sample loss}) (b, d) loss of the \textit{risk-neutral} leader (with MADs) as a function of the common sample size, $k_{lf}$, for $k_l = k_f = 30$ and $\delta_l = \delta_f = 0.1$, evaluated over 100 random test instances. The follower is assumed to be risk-averse (a, b) and risk-neutral (c, d). The dashed line corresponds to the empirical $5\%$ percentile of the relative in-sample loss for the scenario-based semi-pessimistic approximation.}
\label{fig: experiment 8}
\end{figure}	

In our last experiment, we provide a similar comparison of three models, [\textbf{DRI-SP$'$}],~[\textbf{DRI-P}] and [\textbf{DRI}], but with a risk-neutral leader, i.e., $\alpha_l = 0$. The follower is assumed to be either risk-averse (Figures \ref{fig: experiment 8a} and \ref{fig: experiment 8b}) or risk-neutral~(Figures \ref{fig: experiment 8c} and \ref{fig: experiment 8d}). In order to address the unimodularity Assumption \textbf{A6$'$}, we simply replace the follower's feasible set $Y(\mathbf{x})$ with 
\begin{equation} \label{eq: feasible sets packing unimodular} 
Y_{uni}(\mathbf{x}) = \big\{\mathbf{y} \in \mathbb{R}_+^n: \sum_{i = 1}^n y_i \leq \lfloor 0.2n \rfloor, \; \mathbf{y} \leq \mathbf{U}(\mathbf{1} - \mathbf{x}) \big\}.
\end{equation} 

For the risk-averse follower, our observations are similar to those in the ambiguity-free case. However, with the risk-neutral follower, [\textbf{DRI-P}] and [\textbf{DRI}] demonstrate solid in-sample and out-of-sample performance, with a small~out-of-sample error in relation to [\textbf{DRI}$^*$]. Indeed, when both decision-makers are risk-neutral and~$\delta_l = \delta_f = 0.1$, they typically operate in the sample average regime, aligning their objective functions; see Figure \ref{fig: experiments 3 and 4}. Therefore, 
 in most cases, [\textbf{DRI-P}] and~[\textbf{DRI}] provide identical optimal solutions and can be generally recommended for this particular~case.

\subsection{Analysis of running times} \label{subsec: running times}
\begin{figure}[h]\centering	\scalebox{0.8}{
		\begin{subfigure}{0.5\textwidth}
			\centering
			\begin{tikzpicture}
				\begin{axis}[
					axis lines=box,
					xlabel={$k_l$},
					ylabel={Solution times},
					xmin = 5,
					xmax = 300,
					xmode=log,
					log basis x=10,
					ymode=log,
					log basis y=10,
					ymin=0, ymax=100000,
					grid=major,
					legend pos=north east,
					mark size=3pt, ytick = {1,10, 100, 1000},
					xtick={10,30,100,300},
					xticklabels={10, 30, 100, 300}
					]
					
					\addplot[color=black, mark=triangle*, mark size=2pt] coordinates { 
						(5, 0.21) (10, 0.24) (15, 0.28) (20, 0.33)
						(30, 0.48) (40, 0.61) (50, 0.77) (75, 1.17)
						(100, 1.62) (125, 2.19) (150, 2.54) (175, 3.94)
						(200, 4.85) (250, 6.27) (300, 8.66)
					};
					\addlegendentry{\footnotesize Basic model}
					
					\addplot[color=limegreen, mark=square*, mark size=1.5pt] coordinates { 
						(5, 0.46) (10, 0.81) (15, 1.34) (20, 1.72)
						(30, 2.99) (40, 4.07) (50, 8.14) (75, 14.87)
						(100, 24.86) (125, 40.32) (150, 31.59) (175, 86.05)
						(200, 106.45) (250, 195.83) (300, 264.13)
					};
					
					\addlegendentry{\footnotesize Scenario-based approximation ($r_l = 5$)}
					
					\addplot[color=limegreen, mark=*, mark size=2pt] coordinates {
						(5, 1.08) (10, 2.42) (15, 2.86) (20, 3.95)
						(30, 6.79) (40, 9.64) (50, 21.70) (75, 49.90)
						(100, 67.65) (125, 125.70) (150, 188.83) (175, 239.53)
						(200, 353.32) (250, 554.07) (300, 818.03)
					};
					\addlegendentry{\footnotesize Scenario-based approximation ($r_l = 10$)}

					\path[name path=upperA] plot coordinates {
						(5, 0.21 + 0.04) (10, 0.24 + 0.04) (15, 0.28 + 0.02) (20, 0.33 + 0.03)
						(30, 0.48 + 0.04) (40, 0.61 + 0.06) (50, 0.77 + 0.07) (75, 1.17 + 0.08)
						(100, 1.62 + 0.19) (125, 2.19 + 0.25) (150, 2.54 + 0.20) (175, 3.94 + 0.39)
						(200, 4.85 + 0.28) (250, 6.27 + 0.60) (300, 8.66 + 0.64)
					};
					
					\path[name path=lowerA] plot coordinates {
						(5, 0.21 - 0.04) (10, 0.24 - 0.04) (15, 0.28 - 0.02) (20, 0.33 - 0.03)
						(30, 0.48 - 0.04) (40, 0.61 - 0.06) (50, 0.77 - 0.07) (75, 1.17 - 0.08)
						(100, 1.62 - 0.19) (125, 2.19 - 0.25) (150, 2.54 - 0.20) (175, 3.94 - 0.39)
						(200, 4.85 - 0.28) (250, 6.27 - 0.60) (300, 8.66 - 0.64)
					};
					
					\addplot[black!30, fill opacity=0.5] fill between[of=upperA and lowerA];

					\path[name path=upperB] plot coordinates {
						(5, 0.46 + 0.05) (10, 0.81 + 0.08) (15, 1.34 + 0.19) (20, 1.72 + 0.19)
						(30, 2.99 + 0.27) (40, 4.07 + 0.41) (50, 8.14 + 0.68) (75, 14.87 + 2.11)
						(100, 24.86 + 5.21) (125, 40.32 + 5.41) (150, 31.59 + 5.91) (175, 86.05 + 17.37)
						(200, 106.45 + 21.01) (250, 195.83 + 25.91) (300, 264.13 + 48.04)
					};
					
					\path[name path=lowerB] plot coordinates {
						(5, 0.46 - 0.05) (10, 0.81 - 0.08) (15, 1.34 - 0.19) (20, 1.72 - 0.19)
						(30, 2.99 - 0.27) (40, 4.07 - 0.41) (50, 8.14 - 0.68) (75, 14.87 - 2.11)
						(100, 24.86 - 5.21) (125, 40.32 - 5.41) (150, 31.59 - 5.91) (175, 86.05 - 17.37)
						(200, 106.45 - 21.01) (250, 195.83 - 25.91) (300, 264.13 - 48.04)
					};

					\addplot[limegreen!10, fill opacity=0.5] fill between[of=upperB and lowerB];	
					
					\path[name path=upperC] plot coordinates {
						(5, 1.08 + 0.15) (10, 2.42 + 0.18) (15, 2.86 + 0.25) (20, 3.95 + 0.41)
						(30, 6.79 + 0.54) (40, 9.64 + 0.87) (50, 21.70 + 2.90) (75, 49.90 + 6.85)
						(100, 67.65 + 11.32) (125, 125.70 + 24.70) (150, 188.83 + 44.60) (175, 239.53 + 41.47)
						(200, 353.32 + 25.35) (250, 554.07 + 103.97) (300, 818.03 + 139.40)
					};
					
					\path[name path=lowerC] plot coordinates {
						(5, 1.08 - 0.15) (10, 2.42 - 0.18) (15, 2.86 - 0.25) (20, 3.95 - 0.41)
						(30, 6.79 - 0.54) (40, 9.64 - 0.87) (50, 21.70 - 2.90) (75, 49.90 - 6.85)
						(100, 67.65 - 11.32) (125, 125.70 - 24.70) (150, 188.83 - 44.60) (175, 239.53 - 41.47)
						(200, 353.32 - 25.35) (250, 554.07 - 103.97) (300, 818.03 - 139.40)
					};
					
					\addplot[limegreen!10, fill opacity=0.5] fill between[of=upperC and lowerC];	
				\end{axis}
			\end{tikzpicture}
			\caption{$k_f = k_l, n = 10$.}
			\label{fig: experiment 9a}
		\end{subfigure}
		\hfill
		\begin{subfigure}{0.5\textwidth}
			\centering
			\begin{tikzpicture}
				\begin{axis}[
					axis lines=box,
					xlabel={$n$},
					ylabel={Solution times},
					xmin = 6,
					xmax = 50,
					xmode=log,
					log basis x=10,
					ymode=log,
					log basis y=10,
					ymin=0, ymax=1000000,
					grid=major,
					legend pos=north east,
					mark size=3pt, ytick = {1,10, 100, 1000},
					xtick={10,20,24,32,50},
					xticklabels={10, 20,24, 33,50}
					]
					
					\addplot[color=black, mark=triangle*, mark size=2pt] coordinates { 
						(6, 0.19) (7, 0.35) (8, 0.43)
						(10, 0.57) (12, 0.67) (14, 1.13) (17, 2.96)
						(20, 12.21) (24, 29.08) (28, 179.82) (32, 1077.84)
					};
				
					\addlegendentry{\footnotesize Basic model}
					
					\addplot[color=limegreen, mark=square*, mark size=1.5pt] coordinates { 
						(6, 1.27) (7, 1.66) (8, 2.04)
						(10, 3.49) (12, 5.99) (14, 11.10)
						(17, 52.43) (20, 263.70) (24, 1043.94)
					};
				
					\addlegendentry{\footnotesize Scenario-based approximation ($r_l = 5$)}
					
					\addplot[color=limegreen, mark=*, mark size=2pt] coordinates {
						(6, 2.42) (7, 3.49) (8, 4.26)
						(10, 7.64) (12, 13.50) (14, 28.09)
						(17, 191.28) (20, 1095.67)
					};
					\addlegendentry{\footnotesize Scenario-based approximation ($r_l = 10$)}

					\path[name path=upperA] plot coordinates {
						(6, 0.19 + 0.07) (7, 0.35 + 0.04) (8, 0.43 + 0.03)
						(10, 0.57 + 0.06) (12, 0.67 + 0.06) (14, 1.13 + 0.11) (17, 2.96 + 0.62)
						(20, 12.21 + 2.54) (24, 29.08 + 7.47) (28, 179.82 + 89.84) (32, 1077.84 + 466.55)
					};
					
					\path[name path=lowerA] plot coordinates {
						(6, 0.19 - 0.07) (7, 0.35 - 0.04) (8, 0.43 - 0.03)
						(10, 0.57 - 0.06) (12, 0.67 - 0.06) (14, 1.13 - 0.11) (17, 2.96 - 0.62)
						(20, 12.21 - 2.54) (24, 29.08 - 7.47) (28, 179.82 - 89.84) (32, 1077.84 - 466.55)
					};
					
					\addplot[black!30, fill opacity=0.5] fill between[of=upperA and lowerA];

					\path[name path=upperB] plot coordinates {
						(6, 1.27 + 0.10) (7, 1.66 + 0.11) (8, 2.04 + 0.22)
						(10, 3.49 + 0.32) (12, 5.99 + 0.64) (14, 11.10 + 2.12)
						(17, 52.43 + 18.38) (20, 263.70 + 39.09) (24, 1043.94 + 306.63)
					};
					
					\path[name path=lowerB] plot coordinates {
						(6, 1.27 - 0.10) (7, 1.66 - 0.11) (8, 2.04 - 0.22)
						(10, 3.49 - 0.32) (12, 5.99 - 0.64) (14, 11.10 - 2.12)
						(17, 52.43 - 18.38) (20, 263.70 - 39.09) (24, 1043.94 - 306.63)
					};
					
					\addplot[limegreen!10, fill opacity=0.5] fill between[of=upperB and lowerB];	
					
					\path[name path=upperC] plot coordinates {
						(6, 2.42 + 0.32) (7, 3.49 + 0.39) (8, 4.26 + 0.44)
						(10, 7.64 + 0.83) (12, 13.50 + 1.99) (14, 28.09 + 4.54)
						(17, 191.28 + 23.46) (20, 1095.67 + 168.83)
					};
					
					\path[name path=lowerC] plot coordinates {
						(6, 2.42 - 0.32) (7, 3.49 - 0.39) (8, 4.26 - 0.44)
						(10, 7.64 - 0.83) (12, 13.50 - 1.99) (14, 28.09 - 4.54)
						(17, 191.28 - 23.46) (20, 1095.67 - 168.83)
					};
					
					\addplot[limegreen!10, fill opacity=0.5] fill between[of=upperC and lowerC];
					
				\end{axis}
			\end{tikzpicture}
			\caption{$k_l = k_f = 30$.}
			\label{fig: experiment 9b}
	\end{subfigure}}
	\caption{Average solution times in seconds (with MADs) of the basic and semi-pessimistic formulations as a function of $k_l = k_f$ (a) and $n$ (b), for $\delta_l = \delta_f = 0.1$, $k_{lf} = \lfloor\frac{2}{3}k_f\rfloor$ and $\kappa = 0.2$, evaluated over 10 random test instances. The time limit is set to 60 minutes. If at least one out of the 10 test instances is not solved within this limit, we discard the corresponding value of parameter, along with all larger values.}
	\label{fig: experiment 9}
\end{figure}
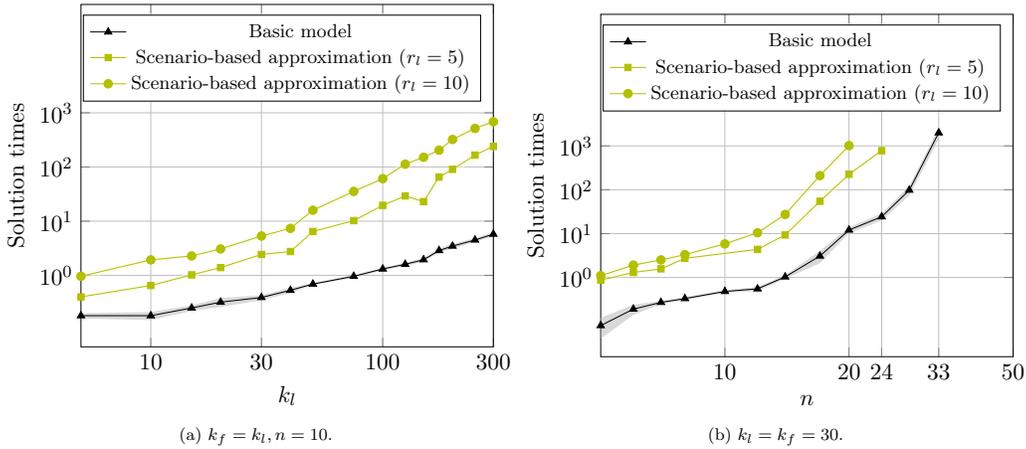

Here, we provide a brief analysis of running times for all three considered formulations. The leader's and the follower's data sets, $\hat{\mathbf{C}}_l$ and $\hat{\mathbf{C}}_f$, are assumed to be identical and, therefore, [\textbf{DRI}] is referred to as a basic model. In Figures~\ref{fig: experiment 9a} and~\ref{fig: experiment 9b}, we examine how solution times for~[\textbf{DRI}] and~[\textbf{DRI-SP$'$}] scale with the sample sizes, $k_l$ and $k_f$, and the problem size~$n$, respectively. All plots are presented on a log-log scale. For each parameter value, we report the average solution times~(with~MADs) over~10 random test instances of~(\ref{packing interdiction problem}) and a unique data set, all generated using the same random seed.  The time limit is set to~60 minutes. If at least one of the 10 test instances remains unsolved within this limit, we discard the corresponding parameter value, along with all larger~values. 

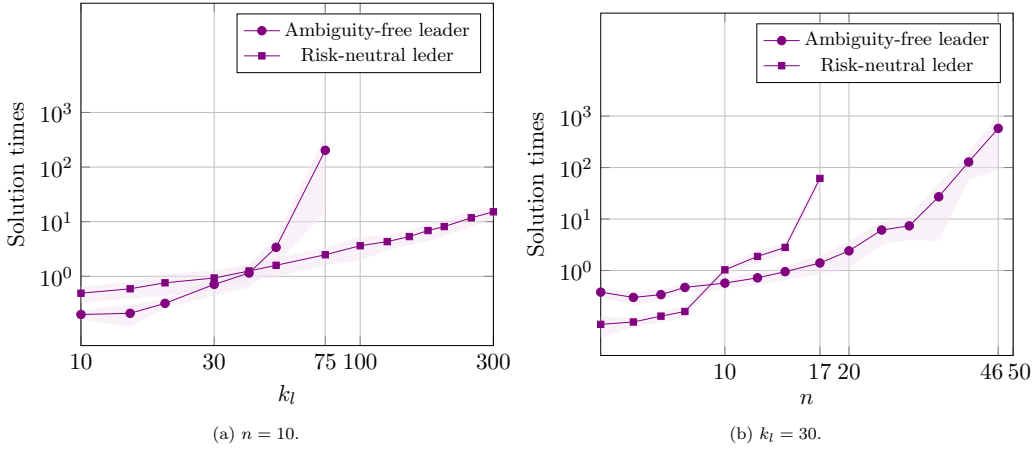
\begin{figure}[h] \centering	\scalebox{0.8}{
		\begin{subfigure}{0.5\textwidth}
			\centering
			\begin{tikzpicture}
				\begin{axis}[
					axis lines=box,
					xlabel={$k_l$},
					ylabel={Solution times},
					xmin = 10,
					xmax = 300,
					xmode=log,
					log basis x=10,
					ymode=log,
					log basis y=10,
					ymin=0, ymax=100000,
					grid=major,
					legend pos=north east,
					mark size=3pt, ytick = {1,10, 100, 1000},
					xtick={10,30,75,100,300},
					xticklabels={10, 30, 75, 100, 300}
					]

					\addplot[color=violet, mark=*, mark size=2pt] coordinates { 
						(10, 0.23) (15, 0.22) (20, 0.49)
						(30, 0.89) (40, 1.35) (50, 3.54)
						(75, 165.54)
					};
					
					\addlegendentry{\footnotesize Ambiguity-free leader}
					
					\addplot[color=violet, mark=square*, mark size=1.5pt] coordinates {
						(10, 0.72) (15, 0.88) (20, 1.07) (30, 1.33)
						(40, 1.75) (50, 2.07) (75, 3.43) (100, 4.78)
						(125, 6.06) (150, 7.57) (175, 9.20) (200, 11.13)
						(250, 16.88) (300, 21.06)
					};
					\addlegendentry{\footnotesize Risk-neutral leader}

					\path[name path=upperB] plot coordinates {
						(10, 0.23 + 0.06) (15, 0.22 + 0.07) (20, 0.49 + 0.20)
						(30, 0.89 + 0.33) (40, 1.35 + 0.44) (50, 3.54 + 1.41)
						(75, 165.54 + 113.00)
					};
					
					\path[name path=lowerB] plot coordinates {
						(10, 0.23 - 0.06) (15, 0.22 - 0.07) (20, 0.49 - 0.20)
						(30, 0.89 - 0.33) (40, 1.35 - 0.44) (50, 3.54 - 1.41)
						(75, 165.54 - 113.00)
					};
					
					\addplot[violet!10, fill opacity=0.5] fill between[of=upperB and lowerB];	
					
					\path[name path=upperC] plot coordinates {
						(10, 0.72 + 0.29) (15, 0.88 + 0.30) (20, 1.07 + 0.40) (30, 1.33 + 0.41)
						(40, 1.75 + 0.64) (50, 2.07 + 0.63) (75, 3.43 + 1.11) (100, 4.78 + 1.54)
						(125, 6.06 + 1.71) (150, 7.57 + 1.94) (175, 9.20 + 2.39) (200, 11.13 + 2.87) (250, 16.88 + 5.26) (300, 21.06 + 5.00)
					};
					
					\path[name path=lowerC] plot coordinates {
						(10, 0.72 - 0.29) (15, 0.88 - 0.30) (20, 1.07 - 0.40) (30, 1.33 - 0.41)
						(40, 1.75 - 0.64) (50, 2.07 - 0.63) (75, 3.43 - 1.11) (100, 4.78 - 1.54)
						(125, 6.06 - 1.71) (150, 7.57 - 1.94) (175, 9.20 - 2.39) (200, 11.13 - 2.87) (250, 16.88 - 5.26) (300, 21.06 - 5.00)	
					};
					
					\addplot[violet!10, fill opacity=0.5] fill between[of=upperC and lowerC];
					
				\end{axis}
			\end{tikzpicture}
			\caption{$n = 10$.}
			\label{fig: experiment 10a}
		\end{subfigure}
		\hfill
		\begin{subfigure}{0.5\textwidth}
			\centering
			\begin{tikzpicture}
				\begin{axis}[
					axis lines=box,
					xlabel={$n$},
					ylabel={Solution times},
					xmin = 5,
					xmax = 50,
					xmode=log,
					log basis x=10,
					ymode=log,
					log basis y=10,
					ymin=0, ymax=100000,
					grid=major,
					legend pos=north east,
					mark size=3pt, ytick = {1,10, 100, 1000},
					xtick={10,17,20,39,50},
					xticklabels={10,17,20,39\;\;,\;\;50}
					]

					\addplot[color=violet, mark=*, mark size=2pt] coordinates { 
						(5, 0.42) (6, 0.40) (7, 0.50) (8, 0.54)
						(10, 0.72) (12, 0.94) (14, 1.28) (17, 2.00)
						(20, 2.62) (24, 9.53) (28, 6.99) (33, 67.68)
						(39, 234.81)
					};
					
					\addlegendentry{\footnotesize Ambiguity-free leader}
					
					\addplot[color=violet, mark=square*, mark size=1.5pt] coordinates {
						(5, 0.11) (6, 0.14) (7, 0.14) (8, 0.18)
						(10, 1.23) (12, 2.42) (14, 3.98) (17, 111.03)
					};
					\addlegendentry{\footnotesize Risk-neutral leader}

					\path[name path=upperB] plot coordinates {
						(5, 0.42 + 0.08) (6, 0.40 + 0.12) (7, 0.50 + 0.15) (8, 0.54 + 0.14)
						(10, 0.72 + 0.24) (12, 0.94 + 0.28) (14, 1.28 + 0.37) (17, 2.00 + 0.69)
						(20, 2.62 + 0.75) (24, 9.53 + 3.73) (28, 6.99 + 3.70) (33, 67.68 + 62.65)
						(39, 234.81 + 218.26)
					};
					
					\path[name path=lowerB] plot coordinates {
						(5, 0.42 - 0.08) (6, 0.40 - 0.12) (7, 0.50 - 0.15) (8, 0.54 - 0.14)
						(10, 0.72 - 0.24) (12, 0.94 - 0.28) (14, 1.28 - 0.37) (17, 2.00 - 0.69)
						(20, 2.62 - 0.75) (24, 9.53 - 3.73) (28, 6.99 - 3.70) (33, 67.68 - 62.65)
						(39, 234.81 - 218.26)
					};
					
					\addplot[violet!10, fill opacity=0.5] fill between[of=upperB and lowerB];	
					
					\path[name path=upperC] plot coordinates {
						(5, 0.11 + 0.02) (6, 0.14 + 0.05) (7, 0.14 + 0.02) (8, 0.18 + 0.03)
						(10, 1.23 + 0.22) (12, 2.42 + 0.66) (14, 3.98 + 1.07) (17, 111.03 + 29.58)
					};
					
					\path[name path=lowerC] plot coordinates {
						(5, 0.11 - 0.02) (6, 0.14 - 0.05) (7, 0.14 - 0.02) (8, 0.18 - 0.03)
						(10, 1.23 - 0.22) (12, 2.42 - 0.66) (14, 3.98 - 1.07) (17, 111.03 - 29.58)
					};
					
					\addplot[violet!10, fill opacity=0.5] fill between[of=upperC and lowerC];
					
				\end{axis}
			\end{tikzpicture}
			\caption{$k_l = 30$.}
			\label{fig: experiment 10b}
	\end{subfigure}}
	\caption{Average solution times in seconds (with MADs) of the pessimistic formulation as a function of $k_l$ (a) and $n$~(b), for $\delta_l = 0.1$ and $\alpha_l = 0.9$, evaluated over 10 random test instances. The time limit is set to 60 minutes. If at least one out of the 10 test instances is not solved within this limit, we discard the corresponding value of parameter, along with all larger values.}
	\label{fig: experiment 10}
\end{figure}	

Figure \ref{fig: experiment 9} suggests that solution times scale well with $k_l$, $k_f$, and $r_l$, showing a relatively moderate growth pattern. However, intuitively, when $n$ increases, the number of binary variables in the MILP reformulations (\ref{MILP reformulation}) and (\ref{semi-pessimistic reformulation MILP}) also grows, leading to an exponential increase in the running time. In this regard, we recall that even the deterministic version of our problem, [\textbf{DP}], is $NP$-hard. As a result, we suggest using models [\textbf{DRI}] and [\textbf{DRI-SP$'$}] for moderately-sized problem instances, with up to around 30 leader's decision variables.

Then, we analyze how solution times for [\textbf{DRI-P}] scale in $k_l$ and~$n$, given that $k_l$ is adjusted to satisfy Assumption \textbf{A6}; see Figures \ref{fig: experiment 10a} and \ref{fig: experiment 10b}, respectively. 
We observe that [\textbf{DRI-P}] for the ambiguity-free leader scales more efficiently with $n$, while the approximation for the risk-neutral leader scales better with~$k_l$. This observation can be explained by the design of Algorithms \ref{algorithm 1}~and~2. Thus, Algorithm \ref{algorithm 1} for the ambiguity-free leader enumerates solutions in the discrete set (\ref{eq: gamma}), whose size grows exponentially with the increase of~$k_l$. At the same time, Algorithm~2 operates with the discrete set (\ref{eq: feasible sets packing unimodular}), whose size is exponential~in~$n$. Overall, although our main focus is on the qualitative results, we believe that Algorithms \ref{algorithm 1} and 2 could be further enhanced to handle larger problem~instances.


\subsection{Summary of numerical results} \label{subsec: summary}
Our computational results can be summarized as follows:
\begin{itemize}
	\setlength{\itemsep}{0.05cm} 
	\setlength{\parskip}{0.05cm} 
	\setlength{\topsep}{2pt} 
	\item Similar to one-stage DRO models, calibrating the Wasserstein radius is essential for optimizing the model's out-of-sample performance. That is, the optimal Wasserstein radii in [\textbf{DRI}] significantly depend on the respective sample sizes, as well as the decision-makers' risk preferences. Furthermore, 
	 risk-aversion is essential for mitigating extremely small or large profits, while risk-neutral policies may offer a more effective approach for optimizing average gains.
	
	\item 
	As shown in Figures \ref{fig: experiment 7b} and \ref{fig: experiment 8b}, the leader's incomplete knowledge of the follower's data may have a much greater impact than the distributional estimation errors. In this regard, we suggest using either the scenario-based approximation [\textbf{DRI-SP$'$}] or the pessimistic approximation [\textbf{DRI-P}]. Thus, [\textbf{DRI-SP$'$}], despite offering only asymptotic theoretical guarantees, exhibits strong in-sample and out-of-sample performance, even with a relatively small number of scenarios and significant uncertainty in the follower's data. In contrast,~[\textbf{DRI-P}] tends to yield poor out-of-sample performance, unless the leader's and the follower's objective functions are well aligned.
	
	\item As a remark, [\textbf{DRI-SP$'$}] may also be applied to the case where the follower's global parameters,~$\alpha_f$ and $\varepsilon_{f}$, are subject to interval uncertainty. However, combining both data and parameter uncertainty is likely to make the scenario-based approximation more conservative, thereby aligning it more closely with the pessimistic approximation. 
\end{itemize}

\section{Conclusion} \label{sec: conclusion}
In this study, we explore a class of general linear interdiction problems with uncertainty in the follower's profit/cost vector. 
The profit vector is assumed to follow an unknown true distribution, which is estimated independently by  the leader and the follower, based on their own historical data. Importantly, due to the conflicting objective criteria, the follower may choose to withhold some private data from the leader. To the best of our knowledge, this is the first data-driven interdiction model in the literature, where the decision-makers have asymmetric information regarding the data. 

To address the distributional uncertainty, both decision-makers in our model solve conventional one-stage Wasserstein distributionally robust optimization (DRO) problems based on the conditional value-at-risk~(CVaR).
Moreover, we propose three distinct optimization models to address the leader's incomplete knowledge about the follower's~data. In the \textit{basic model}, the leader has full knowledge of the follower's data, and therefore the follower's data set is contained in the leader's. In contrast, under the \textit{pessimistic} and the \textit{semi-pessimistic} approximations, the leader either assumes the worst-case feasible follower's policy or constructs an interval-based uncertainty set for the follower's~data.

We establish that the basic model is $NP$-hard, whereas both the pessimistic and the semi-pessimistic approximations are $\Sigma_2^p$-hard. Then, the basic and the semi-pessimistic models are addressed by leveraging a single-level mixed-integer linear programming (MILP) reformulation and discretization with a subsequent reformulation, respectively. In contrast, for two special cases, the pessimistic approximation is solved via a Benders-type decomposition algorithm. 

Notably, we prove that, under a mild assumption, the basic model is asymptotically consistent, i.e., as the decision-makers acquire more data, their optimal solutions and the respective optimal objective function values converge to those of the underlying stochastic
programming problem. Furthermore, we show that our discretization of the leader's uncertainty set in the semi-pessimistic approximation is asymptotically robust. Put differently, as the number of scenarios increases, the resulting
problem provides a robust approximation of the basic model with the true follower's data set. To the best of our knowledge, these results have not been previously introduced in the bilevel optimization context and may provide a basis for addressing more complex hierarchical optimization problems under~uncertainty. 



Finally, we conduct a numerical study of all our models using a set of randomly generated instances of the packing interdiction problem. It turns out that the leader's misspecification of the true follower's policy in our bilevel context may have a greater impact than the individual errors of the decision-makers, which arise from their incomplete knowledge of the true distribution. Overall, our numerical results support the use of the scenario-based semi-pessimistic approximation, even when the follower's data is subject to significant noise, and the pessimistic approximation when the leader's and the follower's objective functions are well aligned.


\textbf{Competing Interests:} None. 

\textbf{Availability of Data and Materials:} Data were generated through numerical experiments, with the code available at: \url{https://github.com/sk19941995/Data-driven-interdiction}.

\textbf{Declaration of generative AI and AI-assisted technologies in the manuscript preparation process:}
During the preparation of this work, the authors used ChatGPT for language refinement and clarity. After using this tool, the authors reviewed and edited the content as needed and take full responsibility for the content of the published article.

\singlespacing 
\bibliographystyle{apa}
\bibliography{bibliography}

\includepdf[pages=-]{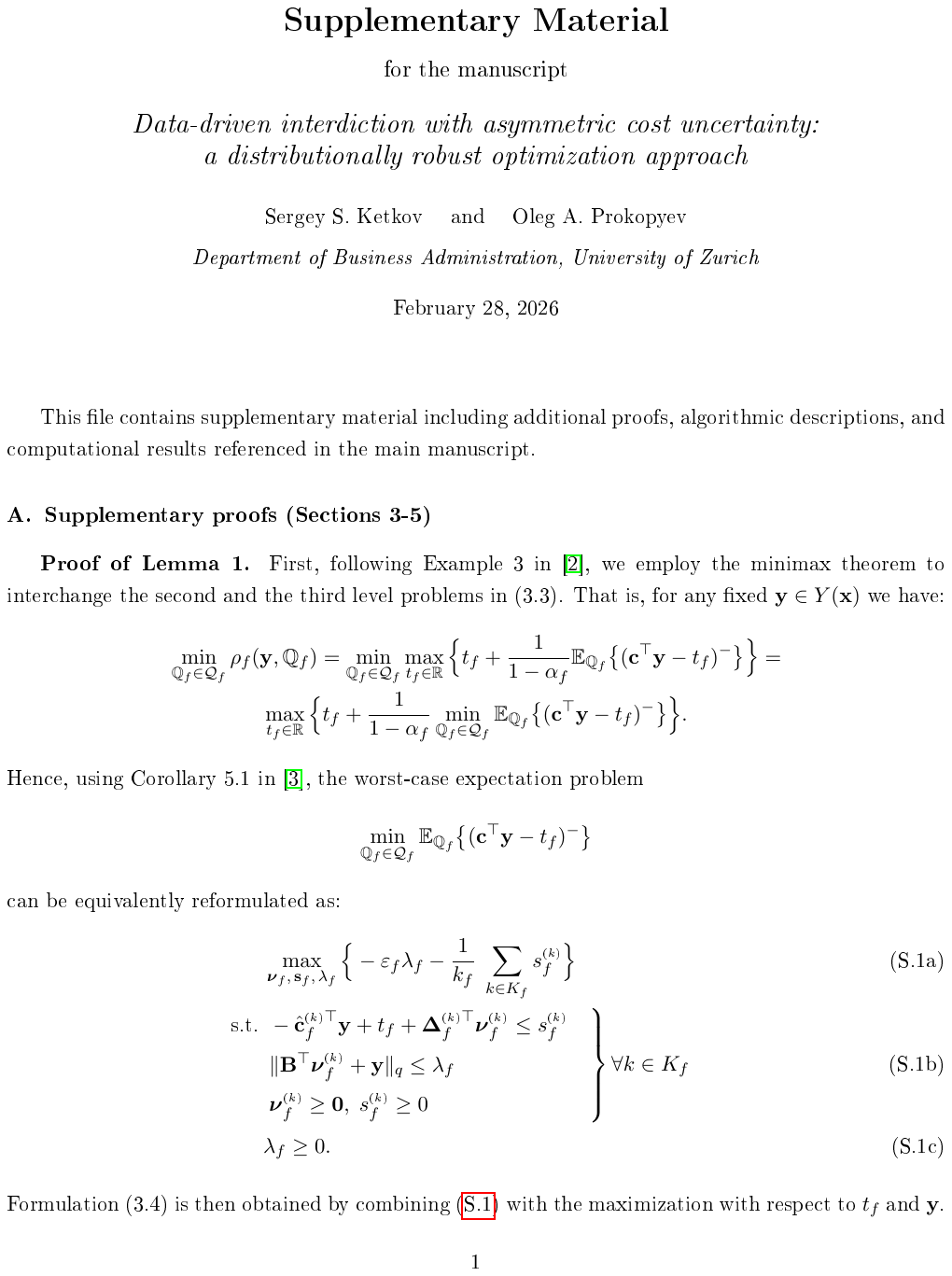}
	
\end{document}